\input amstex

\def\b1{\text{\bf 1}}

\def\CA{{\Cal A}}

\def\CC{{\Cal C}}
\def\CD{{\Cal D}}
\def\CF{{\Cal F}}

\def\CH{{\Cal H}}
\def\CI{{\Cal I}}
\def\CJ{{\Cal J}}
\def\CK{{\Cal K}}
\def\CG{{\Cal G}}
\def\CL{{\Cal L}}
\def\CM{{\Cal M}}
\def\CN{{\Cal N}}
\def\CO{{\Cal O}}

\def\CS{{\Cal S}}
\def\CT{{\Cal T}}

\def\CV{{\Cal V}}

\def\gr{\text{gr}}

\def\End{\text{End}}
\def\Hom{\text{Hom}}

\def\et{\text{\'et}}
\def\h{\text{h}}
\def\hotimes{\widehat{\otimes}}

\def\#{\,\check{}}

\def\fm{{\frak m}}

\def\ft{{\frak t}}

\def\con{{\text{con}}}

\def\ss{{\text{ss}}}

\def\id{\text{id}}

\def\Ker{\text{Ker}}

\def\Spec{\text{Spec}}

\def\pcr{\text{pcr}}

\def\dR{\text{dR}}
\def\crys{\text{crys}}

\def\eff{\text{eff}}

\def\perf{\text{perf}}
\def\nd{\text{nd}}
\def\st{\text{st}}
\def\lim{\text{lim}}

\def\holim{\text{holim}}
\def\BdR{\text{B}_{\text{dR}}}
\def\AdR{\text{A}_{\text{dR}}}
\def\Bcrys{\text{B}_{\text{crys}}}
\def\Acrys{\text{A}_{\text{crys}}}
\def\Bst{\text{B}_{\text{st}}}

\def\holim{\text{holim}}

\def\ds{\diamondsuit}
\def\limleft{\mathop{\vtop{\ialign{##\crcr
  \hfil\rm lim\hfil\crcr
  \noalign{\nointerlineskip}\leftarrowfill\crcr
  \noalign{\nointerlineskip}\crcr}}}}
\def\limright{\mathop{\vtop{\ialign{##\crcr
  \hfil\rm lim\hfil\crcr
  \noalign{\nointerlineskip}\rightarrowfill\crcr
  \noalign{\nointerlineskip}\crcr}}}}

\def\hra{\hookrightarrow}
\def\iso{\buildrel\sim\over\rightarrow} 

\def\lra{\longrightarrow}


\documentstyle{amsppt}
\NoBlackBoxes


\topmatter  \title On the  crystalline period map \endtitle  \author A.~Beilinson \endauthor\address Department of Mathematics, University of Chicago, Chicago, IL 60637 \endaddress \email sasha\@math.uchicago.edu \endemail  \thanks The author was  supported in part by NSF grant DMS-1001660. \endthanks    \subjclassyear{2010}   \subjclass Primary 14F30, 14F40; secondary 14F20  \endsubjclass    \keywords p-adic periods, log crystalline cohomology, h-topology, alterations \endkeywords      \date   
 \enddate    
 \dedicatory To Irene  
  \enddedicatory  \endtopmatter
 
\document
 
This article is a direct continuation of \cite{B}. It contains a simple proof of comparison theorems in $p$-adic Hodge theory (the Fontaine-Jannsen conjecture). Different proofs were found earlier by Faltings, Niziol, and Tsuji, the case of open varieties treated by Yamashita. An alternative approach, based on a canonical identification of the log crystalline cohomology for lci maps with  the {\it noncompleted} (for the Hodge filtration)  derived de Rham complex, was developed by Bhatt \cite{Bh2}.

 Let $\CV ar_{\bar{K}}$ be  the category of algebraic varieties over an algebraic closure 
 $\bar{K}$ of a $p$-adic field $K$.
Our principal character is a natural  h-sheaf  of dg algebras $\CA_{\crys}$ on   $\CV ar_{\bar{K}}$
equipped with a Frobenius action; let  $X\mapsto R\Gamma_{\!\crys}(X):= R\Gamma (X_{\h},\CA_{\crys})$, $X\in  \CV ar_{\bar{K}}$, be
 the corresponding cohomology theory.  We construct $\CA_{\crys}$ applying h-localization procedure, as in \cite{B} 2.6, to 
the absolute (i.e., relative to  $\Bbb Z_p$)  log crystalline cohomology  of log $O_{\bar{K}}$-schemes coming from  semi-stable  pairs. The algebra  $\CA_{\crys}(\Spec\,\bar{K})$ equals the Fontaine ring $A_{\crys}$.  The crystalline $p$-adic Poincar\'e lemma asserts that $\CA_{\crys}\otimes^L \Bbb Z/p^n$ is a constant  h-sheaf with fiber $A_{\crys}/p^n$. It yields a morphism of dg algebras
$\rho_{\crys}: 
R\Gamma_{\!\crys}(X)\to R\Gamma_{\!\et}(X,\Bbb Z_p )\otimes_{\Bbb Z_p}\! \Acrys$. This
{\it  crystalline period map}  becomes an isomorphism after   $\Bcrys$-localization, which implies the Fontaine-Jannsen conjecture.

The article runs as follows: \S 1, that takes most of the pages, is an exposition of  log crystalline cohomology and  Hyodo-Kato theory.  In \S 2 we define $\CA_{\crys}$ and prove  the crystalline $p$-adic Poincar\'e lemma. 
Like its Hodge-completed counterpart  of \cite{B} 3.3, the Poincar\'e lemma comes from  Bhatt's theorem \cite{B} 4.3, \cite{Bh1} 1.1, but to deduce it we use the conjugate spectral sequence instead of the Hodge-de Rham one. 
We also define the Hyodo-Kato cohomology\footnote{It is constructed by the same procedure as $R\Gamma_{\!\crys}(X)$ using the Hyodo-Kato cohomology of log 
schemes coming from  semi-stable  pairs as the input.}  $R\Gamma_{\!\text{H}\!\text{K}}(X)$
and  identify $R\Gamma_{\!\crys}(X)\otimes\Bbb Q$ with its twisted form; since, by Hyodo-Kato,  another twist of 
$R\Gamma_{\!\text{H}\!\text{K}}(X)$
equals $R\Gamma_{\!\dR}(X)$, we get control on $R\Gamma_{\!\crys}(X)\otimes\Bbb Q$.  In \S 3 we prove the above isomorphism assertion for $\rho_{\crys}$, deduce from it the Fontaine and  Fontaine-Jannsen  conjectures C$_{\crys}$, C$_{\text{st}}$, C$_{\text{pst}}$, and show that the Hodge filtration completion of   $\rho_{\crys}$  equals the period map $\rho_{\dR}$ from \cite{B} 3.5.

I am very grateful to Bhargav Bhatt, Volodya Drinfeld, Luc Illusie, Kazuya Kato, and the referee  for  valuable comments, suggestions, and corrections.

\head 1. Log crystalline cohomology: a review. \endhead

The section can be divided into five parts: \newline {\it (i) Log preliminaries} (1.1--1.4). The main reference for log schemes is \cite{K1}; for an expanded exposition,  see    \cite{GR} Ch.~3, 7 and \cite{Og1}. 
We show  that  the key existence result for pd-envelopes, proved in \cite{K1} 5.4 for fine log schemes, remains true  for integral quasi-coherent log schemes, which is our preferred log crystalline setting.
\newline{\it (ii) Log crystalline basics} (1.5--1.12). One finds a  concise exposition of the basics of Berthelot's theory \cite{Ber}   in  \cite{BhdJ}; we discuss in similar vein  Kato's log crystalline theory  \cite{K1} \S\S5-6, \cite{HK} \S2, including the comparison theorem of  Illusie and Olsson.\newline{\it (iii) Frobenius log crystals  and the Hyodo-Kato theory} (1.13--1.16). We deduce the key global results of Hyodo-Kato directly from the identification of nondegenerate Frobenius crystals up to isogeny on $\Spec\, (O_K /p )$ with  $(\varphi ,N)$-modules (cf.~\cite{Fa} \S2,
\cite{Ol2} 5.3). The original approach of  \cite{HK} \S 5 (see also \cite{Ill3} \S 7) is local and uses de Rham-Witt complexes. The core of both arguments goes back to Dwork (\cite{Ka} 3.1). \newline{\it (iv) The Fontaine rings and absolute
 crystalline cohomology of log schemes over $O_{\bar{K}}$} (1.17--1.18). This is an exposition of \cite{K2} \S\S 3--4. \newline{\it (v) Log de Rham complex in characteristic} 0 (1.19). This is an exposition of a theorem of Ogus \cite{Og2} Th.~1.3.

\subhead{\rm 1.1}\endsubhead  {\it Log schemes.}  A log scheme is denoted as $(Z,\CM )=(Z,\CM,\alpha )$, $Z$ is the underlying scheme, $\CM$ is the (\'etale) monoid sheaf, $\alpha :\CM \to\CO_Z$ is the log structure map.
A {\it chart}, or {\it $M$-chart}, for it is a homomorphism of monoids $c_M : M\to \Gamma (Z,\CM )$ (the chart map) which yields an isomorphism $M^a_Z \iso \CM$; here $M_Z$ is the constant sheaf of monoids on $X_{\et}$ with fibers $M$ and $M^a_Z =(M^a_Z ,\alpha )$ is the log structure generated by the  
prelog one $\alpha c_M : M_Z \to \CO_Z$. A chart is {\it integral} if 
$M$ is integral (i.e., the canonical map $M \to M^{\text{gr}}$ - the group completion of $M$ - is injective) and {\it fine} if $M$ is fine, i.e., integral and finitely generated.

We say that $(Z,\CM )$ is {\it log affine} if $Z$ is affine and 
$\Gamma (Z,\CM )$ generates $\CM$,\footnote{I.e., every local section of $\CM$ is the product of a global section and an invertible function.}  {\it quasi-coherent} (\cite{K1} 2.1) if it admits a chart \'etale locally,  {\it integral}\footnote{In present article ``integral" in the sense of ``reduced irreducible" (scheme) is never used.} (\cite{K1}  2.2) if $\CM$ is integral, and {\it fine} (\cite{K1}  2.3) if it admits a fine chart \'etale locally. We identify schemes with log schemes with trivial log structure, and write $X$ for  $(X,\CO^\times_X )$.

The category of log schemes has finite inverse limits; the underlying scheme functor $(Z,\CM ) \mapsto Z$ commutes with inverse limits (\cite{K1} 1.6). 

Let $f: (Z,\CM )\to (S,\CL )$ be a map of log schemes. A {\it chart} for it consists
of charts $M\to \Gamma (Z,\CM)$, $L\to\Gamma (S,\CL )$ and a map of monoids $L\to M$ such that the evident diagram commutes. Such a chart  is {\it fine} if both $M$, $L$ are fine.

\proclaim{Proposition}  If $(S,\CL )$ is fine, $(Z,\CM )$ is integral quasi-coherent, then $f$ admits a chart \'etale locally. If $(Z,\CM )$ is also fine, then $f$ admits a fine chart \'etale locally.
\endproclaim 

\demo{Proof} The second claim is  \cite{K1} 2.10. We modify Kato's proof to 
cover the first claim. Let us construct an \'etale neighborhood  with a chart. Let  $\{ \ell_i \}_{i\in I}$ be a finite subset of $ \CL_s$ that generates $\CL_s /\CO^\times_{S, s}$. Let 
$K\subset L\subset \Bbb Z^I$ be the preimages of $\CO^\times_{S, s}$ and $\CL_s $ by the map $\Bbb Z^I \to \CL_s^{\text{gr}}$, $e_i\mapsto\ell_i$. 
Let $\{ k_j \}_{j \in J}$ be a base of $K$, so $k_j = \Sigma n_{j i} e_i - \Sigma n'_{j i}e_i$, where $n_{j i} ,n'_{j i}\in\Bbb N$. Then $ e_i $, $k_j$, and $-k_j$ generate $L$ as a monoid, with $k_j +\Sigma n'_{j i}e_i =\Sigma n_{j i} e_i$, $k_j  + (-k_j )=0$ being the full set of  relations as of {\it integral} monoid. As in \cite{K1} 2.10, after localizing $S$, the map $L\to \CL_s$ lifts to a chart $c_L : L\to \Gamma (S,\CL )$. 

Pick any $z\in Z$. Let us construct its \'etale neighborhood with a chart for $f$. Localizing $Z$, we can assume that $(Z,\CM )$ has  chart $M''$ and  $f^* (\ell_i )=m_i g_i$, where $m_i \in M''$, $g_i \in\Gamma (Z,\CO^\times_Z )$. Set
$h_j := f^* (c_L (k_j ))\in \Gamma (Z,\CO^\times_Z )$. 
Notice that $M' := M''\oplus \Gamma (Z,\CO^\times_Z )$, equipped with an evident chart map $c_{M'}$,
is another chart for $\CM$. The elements $(0, h_j ) + \Sigma n'_{j i}(m_i , g_i )$ and  $\Sigma n_{j i}(m_i , g_i )$ of $M'$  are identified by $c_{M'}$.
 Let $M$ be the {\it integral} monoid quotient  of $M'$ modulo the relations $(0, h_j ) + \Sigma n'_{j i}(m_i , g_i )= \Sigma n_{j i}(m_i , g_i )$. Then $c_{M'}$ factors through
 $c_M : M\to \Gamma (Z,\CM )$, which is again a chart, and we have
a map of monoids $L\to M$ which assigns to $e_i$, $k_j$ the images of $(m_i ,g_i )$, $(0, h_j )$ in $M$. This is the promised chart for $f$. 
\qed
\enddemo

 \remark{Exercises} (i)  Every quasi-coherent $(Z,\CM )$ admits the maximal closed integral log subscheme $(Z,\CM )^{\text{int}}$. The log scheme $(Z,\CM )^{\text{int}}$ is quasi-coherent, and the functor $(Z,\CM )\mapsto (Z,\CM )^{\text{int}}$ is right adjoint to the embedding of the category of quasi-coherent integral log schemes into that of quasi-coherent log schemes.
 \newline
 (ii) If $ (Z_i ,\CM_i )$, $i=1,2$,  are integral quasi-coherent log schemes over a fine log scheme
  $(S,\CL )$, then  $(Z_1 ,\CM_1 )\times_{(S,\CL )}  (Z_2 ,\CM_2 )$ is quasi-coherent.\footnote{Hint: Use the fact that \'etale locally both projections $(Z_i ,\CM_i )\to (S,\CL )$ admit charts with the same $L$ (follows from the proof of the proposition).} 
  \newline
  (iii) Suppose $(T,\CN )$ is integral and
   $(Z,\CM )\hra (T,\CN )$ is an exact closed embedding defined by a nil ideal $\CJ \subset \CO_T$. If $Z$ is affine, then for any $m\in \Gamma (Z,\CM )$ the set of its liftings to $\Gamma (T,\CN )$ is a $\Gamma (T, (1+\CJ )^\times )$-torsor.\footnote{Hint: $T$ is affine and $H^{>0}(T_{\et}, (1+\CJ )^\times )=0$. }  $(T,\CN )$  is quasi-coherent  if (and only if) such is  $(Z,\CM )$.\footnote{Hint: If $Z$ is affine and $M\to\Gamma(Z,\CM )$ is a chart for $\CM$, then $M\times_{\Gamma(Z,\CM )}\Gamma (T,\CN )\to \Gamma (T,\CN )$ is a chart for $\CN$.}
 \endremark

\remark{Remarks} (i) Let $f : (Z,\CM )\to (S,\CL )$ be a map of fine log schemes; suppose $S$ is affine, $Z/S$ is finitely presented. Then there is $f' : (Z',\CM' )\to (S',\CL' )$ having the same properties  with $S'$  affine of finite type over $\Bbb Z$, and a map  $(S,\CL )\to (S',\CL')$ such that $\CL$ is the pullback of $\CL'$ and $f$ isomorphic to the pullback of $f'$; for $f$ log smooth, one can find log smooth $f'$. 
\newline  (ii) Suppose $f$ as in (i) is log smooth, $Z$ is affine, and we have a closed exact embedding $(S,\CL )\hra (T,\CN )$ defined by a nil ideal $\CJ \subset \CO_T$; then $f$ is the pullback of some log smooth $f_T : (Z_T ,\CM_T )\to (T,\CN )$ defined uniquely up to an isomorphism.\footnote{By (i), we are reduced to the situation when $\CJ$ is nilpotent; now the assertion is \cite{K1} 3.14.}
 \endremark

\enspace 

 We denote by $\Bbb A^1_{(S,\CL )}$ and $\Bbb A^{(1)}_{(S,\CL )}$ the affine line and the{ \it logarithmic affine line}  over a log scheme $(S,\CL )$. Both  equal $\Spec\, \CO_S [t]$ as  schemes, the log structures are generated by, respectively, $\CL$ and $\CL \oplus \Bbb N$, the latter log structure map sends $n\in \Bbb N$ to $t^n$. For a log $(S,\CL )$-scheme $(Z,\CM )$, a map from it to $\Bbb A^1_{(S,\CL )}$ is (the same as) a section of $\CO_Z$, a map to $\Bbb A^{(1)}_{(S,\CL )}$ is a section of $\CM$. The group scheme $ \Bbb G_a$ acts on $\Bbb A^1_{(S,\CL )}$ by  translations,   $\Bbb G_m$ acts on $\Bbb A^{(1)}_{(S,\CL )}$ by homotheties.
For a  (possibly infinite) 
set $I$ we denote by $\Bbb A^{I}_{(S,\CL )}$, $\Bbb A^{(I)}_{(S,\CL )}$ the fiber products of $I$ copies of $\Bbb A^1_{(S,\CL )}$, $\Bbb A^{(1)}_{(S,\CL )}$.

\subhead{\rm 1.2}\endsubhead  {\it Log pd-schemes.} For us, a {\it log pd-scheme} is a log scheme $(T,\CM_T )$ equipped with a pd-ideal $(\CJ_T ,\delta )$ (so $\CJ_T$ is a quasi-coherent ideal in $\CO_T$, $\delta$ a pd structure on it). A {\it pd-thickening} of a
 log scheme $(Z,\CM )$ is an exact closed embedding  (\cite{K1} 3.1) of $(Z,\CM )$ into a log pd-scheme
as above such that $\CO_Z \iso \CO_T /\CJ_T$.   Log pd-schemes, i.e., pd-thickenings, form naturally a category. 
We often abbreviate $(T,\CM_T,\CJ_T ,\delta )$ to  $(T,\CM_T,\CJ_T )$; a pd-thickening as above is denoted by $(Z,T,\CM_T )$ or simply  $(Z,T )$.

Let $S^\sharp =(S,\CL,\CI ,\gamma )$ be a  log pd-scheme with $p\in\CO_S$  nilpotent, $p^n =0$. 
Below a {\it log scheme over $S^\sharp$}, or a {\it log $S^\sharp$-scheme}, means a log scheme $(Z,\CM )$ over $(S,\CL )$ such that $\gamma$ extends to $Z$ (i.e., to $\CI \CO_Z$). A {\it log pd-$S^\sharp$-scheme} is a log $S^\sharp$-scheme $(T,\CM_T )$ equipped with a pd-ideal $(  \CJ_T , \delta )$ such that $\gamma$ and $\delta$ extend to a pd structure on $\CJ_T +\CI \CO_T $. We can consider {\it pd-$S^\sharp$-thickenings} of  $(Z,\CM )$, etc. Notice that  $\CJ_T$ is a nil ideal
(since for $a\in \CJ_T$ one has $a^p\in p\CJ_T$ and $p^n =0$). Thus 
pd-$S^\sharp$-thickenings of $(Z,\CM_Z )$ have $Z_{\et}$-local nature.

The categories of log pd-schemes, log $S^\sharp$-schemes, and log pd-$S^\sharp$-schemes have finite inverse limits. For a group scheme $G$ we denote by $G^\sharp$ its pd-envelope at $1\in G$; this is a group pd-scheme. One has $ \Bbb G_m^{\sharp}((Z,T))=\Gamma (T,(1+\CJ_T )^\times )$.

\subhead{\rm 1.3}\endsubhead {\it Log pd-envelopes.} Suppose $S^\sharp$ from 1.2 is such that $(S,\CL )$ is 
quasi-coherent. Let $\CC_{S^\sharp}$ be the category whose objects are locally closed embeddings\footnote{I.e., $(Z,\CM )$ is a closed log subscheme of an open part of $(Y,\CN )$.}
 $i_Y : (Z,\CM )\hra (Y,\CN )$ of log $S^\sharp$-schemes such that $(Z,\CM )$ is integral quasi-coherent, $(Y,\CN )$ is quasi-coherent (the morphisms are maps of $(Y,\CN)$'s that preserve $(Z,\CM )$'s). Let $\CT_{S^\sharp}$ be the category of log pd-$S^\sharp$-thickenings $(Z,T )$ as in 1.2 such that $(T,\CM_T )$, hence $(Z,\CM)$, is integral quasi-coherent.

\proclaim{Theorem} The  evident functor $\CT_{S^\sharp} \to \CC_{S^\sharp}$  (forgetting of the pd structure on $\CJ_T$) admits a right adjoint.
\endproclaim

The theorem says that every
 $i_Y $ in $\CC_{S^\sharp}$ admits a  factorization $$(Z,\CM )\buildrel{i_T}\over\hra (T,\CM_T )\buildrel{r}\over\to (Y,\CN ) \tag 1.3.1$$ with $(Z,T)\in\CT_{S^\sharp}$ which is universal. Such $(Z,T)$ is    called  {\it the pd-$S^\sharp$-envelope} of $i_Y$.
  
\demo{Proof} Let us construct $i_T$. When $(Z,\CM )$ is fine, $(Y,\CN )$ is coherent, this was done by Kato  \cite{K1} 5.4. The idea of  the
 general  argument:  after \'etale localization and adding extra polynomial variables to $Y$, one can realize $i_Y$ as a filtered inverse limit of embeddings of fine log schemes,  which brings $i_T$ by loc.~cit.; the extra polynomial variables are then  factored off. Here are the details:

(a) Suppose we have some factorization $i_T$ as in (1.3.1) and $U\in Y_{\et}$. Consider the $U$-localizations
$i_U : (Z_U, \CM_{Z_U}):= (Z,\CM)\times_Y U \hra (U,\CN_U )$, $(T_{U}, \CM_{T_{U}}):= (T,\CM_T )\times_Y U$. If $(Z,T)$ is the pd-$S^\sharp$-envelope of $i_Y$, then $(T_{U}, \CM_{T_{U}})$ is the pd-$S^\sharp$-envelope of $i_U$. The converse is true if $U$ is a covering of $Y$.

(b) {\it The assertion of the theorem has \'etale local nature:} If $U$ is an \'etale covering of $Y$ and $i_{U}: (Z_{U},\CM_{Z_{U}})\hra (U,\CN_{U})$ admits the
pd-$S^\sharp$-envelope, then $i$ admits the pd-$S^\sharp$-envelope. Proof: Let us construct the the pd-$S^\sharp$-envelope $(Z,T)$.
We can assume that $Z$, $Y$, $U$ are affine schemes,  $Z$ is closed in $Y$.
By (a), if $U' \in Y_{\et}$ admits a morphism to $ U$, then
$i_{U'}$ admits the pd-$S^\sharp$-envelope. Take $U':= U\times_Y U$; let  $(T_{U},\CM_{T_U})$, $(T_{U'},\CM_{T_{U'}})$ be the pd-envelopes. Then $T_{U'}$ is an equivalence relation on $T_U$ by (a), and
 the schemes $T_{U}$, $T_{U'}$ are affine (since such are $Z_{U}$, $Z_{U'}$). Define $T$ as the affine scheme quotient of $T_U$ modulo the equivalence relation $T_{U'}$, i.e., $T$ is an affine scheme with $\Gamma (T,\CO_T )$ equal to the equalizer of $\Gamma (T_{U},\CO_{T_{U}})\rightrightarrows\Gamma (T_{U'},\CO_{T_{U'}})$. 
It is a pd-thickening of $Z$ over $Y$ relative to $(S, \CI, \gamma )$. The map $T_U \to T$ is an \'etale covering and $T_{U'}=T_U \times_T T_U$ (indeed, since our claim has \'etale local nature with respect to $Y$, to check it we can assume that $i$ admits pd-$S^\sharp$-envelope, and then it follows by (a)).
 We define the log structure $\CM_T$ by the \'etale descent from $\CM_{T_U}$.
 Then $(Z,T)$ is the pd-$S^\sharp$-envelope of $i_Y$ by (a).

(c) The pd-$S^\sharp$-envelope of $i_Y$ equals that of $(Z,\CM )\hra (Y,\CN )^{\text{int}}$, so we can assume  that $(Y,\CN )$ is integral. The pd-$S^\sharp$-envelope of $i_Y$ equals its  pd-$Y^\sharp$-envelope, where $Y^\sharp = (Y,\CN, \CI \CO_Y )$, so we can assume that  $(Y,\CN )= (S,\CL )$. By (b),
 we can assume that  $Z$, $Y$ are affine, $Z$ closed in $Y$, and one has integral monoids $M,N$ and charts  $ M\to \Gamma (Z,\CM )$, $N\to \Gamma (Y,\CN )$.

 Pick a set of generators $\{ m_i \}_{i\in I}$ of $M$. Let $i_{V} : (Z,\CM )\hra (V,\CK ):=  \Bbb A^{(I)}_{(Y,\CN )}$  be the lifting of $i_Y$  that corresponds to $\{ m_i \}$ (see 1.1).
 Let us show that
{\it $i_{V}$ admits a pd-$Y^\sharp$-envelope $(\tilde{T},\CM_{\tilde{T}})$.} 

For a subset $J $ of $I$ let $M(J)$ be the submonoid of $M$ generated by $m_{j}$, $j\in J$, and $\CM(J) \to \CO_Z$ be the log structure generated by  prelog one $M(J)_{Z}\to\CO_Z$; there is an evident morphism $\CM (J)\to\CM$. Let
$\CS$ be the set  of triples $ s=(N^s  ,I^s ,\gamma^s )$,  $N^s $ is a finitely generated submonoid of $N$, $ I^s$ is a finite subset of $I$,  $\gamma^s : N^s \to \Gamma (Z,\CM (I^s ))$ is a map of monoids which lifts the composition $N^s \subset N\to\Gamma (Y,\CN)\to \Gamma (Z, \CM )$. The natural order on $\CS$ makes it a directed set, and for every finitely generated $N'\subset N$ and finite $I'\subset I$ one can find $s\in\CS$ with $N'\subset N^s$, $I'\subset I^s$.

For $s\in\CS$  let $\CN^s \to \CO_Y$, $\CK^s \to \CO_{V}$ be the log structures generated by the prelog ones $N^s_Y \to\CO_Y$,  $N_{V}^s \oplus \Bbb N  [I^s ]_V \to \CO_{V}$;   set $\CM^s :=\CM (I^s )$. The map of monoids $N^s \oplus \Bbb N [I^s ]\to \Gamma (Z,\CM^s )$, which is $\gamma^s$ on the first component and the evident map on the second one, extends $Z\hra V$ to a closed embedding of fine log schemes $i^s_{V} :(Z, \CM^s )\hra (V,\CK^s )$. Let $(Z, \CM^s )\hra (\tilde{T}^s ,\CM_{\tilde{T}^s}) \to (V,\CK^s )$ be its pd-$Y^{s\sharp}$-envelope, where $Y^{s\sharp}=(Y,\CN^s ,\CI\CO_{Y})$. By \cite{K1} 5.4, it exists and $\tilde{T}^s$ is affine. When $s\in \CS$ varies, the pd-envelopes form an $\CS$-projective system. Its projective limit is the promised pd-$Y^\sharp$-envelope $(\tilde{T},\CM_{\tilde{T}})$  of $i_{V}$. Indeed, $(\tilde{T},\CM_{\tilde{T}})$  is evidently integral, hence it is quasi-coherent by Exercise (iii) in 1.1, so $(Z,\tilde{T})\in \CT_{Y^\sharp}$,  and the universality property is evident.

(d) As in 1.1,  $(V,\CK )$ carries the  $\Bbb G_m^I$-action.\footnote{For $g=(g_i)\in \Bbb G_m^I$ one has  $g^* (t_i )=g^{-1}_i t_i$.} By the universality property, the group pd-scheme $\Bbb G^{\sharp I}_m$ (see 1.2) acts on  $(\tilde{T},\CM_{\tilde{T}})$. 
 Let $(T,\CM_T )$ be the quotient log pd-$Y^\sharp$-scheme, so $T=\Spec\, (\CO (\tilde{T})^{\Bbb G^{\sharp I}_m})$ and  $\CM_T$ is the sheaf of $ \Bbb G^{\sharp I}_m$-invariant sections of $\CM_{\tilde{T}}$. Let us show that the evident map $i_T \! : (S,\CM )\to (T,\CM_T )$
 is the pd-$Y^\sharp$-envelope of $i_Y$.
 
By construction, $(T,\CM_T )$ is a log pd-$Y^\sharp$-scheme with $T$ affine and $\CM_T$ integral.
 $i_T $ is a closed embedding of log schemes since $(S,\CM )$ is a closed log subscheme of $(Y,\CN )$. It is exact since $i_V$ is exact. By Exercise (iii) in 1.1, $(T,\CM_T )$ is quasi-coherent. Thus $(Z,T)\in \CT_{Y^\sharp}$. 
  It remains to prove  the universality property. 
 
 Suppose 
 $(Z',T')\in\CT_{Y^\sharp}$ is such that the structure map $(T',\CM'_{T'})\to (Y,\CN )$  sends
 $(Z',\CM' )$ to $(Z,\CM )$; we want to show that there is a unique
 $\CT_{Y^\sharp}$-morphism  $(Z',T')\to (Z,T)$. We can assume that $Z'$ is affine, so the set of $\CT_{Y^\sharp}$-morphisms  $(Z',T')\to (Z, \tilde{T})$ is 
 a $\Bbb G^{\sharp I}_m ((Z',T'))$-torsor (see Exercise (iii) in 1.1). So there is a unique  $(Z',T')\to (Z,T)$ that can be lifted to $(Z, \tilde{T})$. To finish the proof, 
 it suffices to check that every  $s: (Z,T)\to  (Z, \tilde{T})$
 is a section of the  projection $p: (Z, \tilde{T})\to (Z,T)$. Now $ps= \id_{(Z,T)}$ amounts to $psp=p$, which is clear since $sp =g\cdot \id_{ (Z, \tilde{T})   }$ for some $g \in \Bbb G^{\sharp I}_m ( (Z, \tilde{T}))$.
   \qed\enddemo

 \subhead{\rm 1.4}\endsubhead {\it Log pd-smooth thickenings.} An object $(Z,T)$ of $\CT_{S^\sharp}$  (see 1.3) is said to be {\it pd-$S^\sharp$-smooth}  if  the next condition is satisfied: Suppose we have $(Z',T')\in \CT_{S^\sharp}$ 
 such that $Z'$ is affine (as a plain scheme); then
 any map of log $(S,\CL )$-schemes $(Z',\CM')\to (Z,\CM )$ can be extended to a morphism $(Z',T')\to (Z,T)$ in $\CT_{S^\sharp}$.

\remark{Remarks} (i) Suppose $(Z,\CM )$ is an integral quasi-coherent log $S^\sharp$-scheme, and we have its closed embedding $i$ into  $\Bbb A^{I}_{(S,\CL )}\times_{(S,\CL )}\Bbb A^{(J)}_{(S,\CL )}$, where $I$, $J$ are two sets 
 (see 1.1).  Then the pd-$S^\sharp$-envelope
 $(P,\CM_P  )$ of $i$ is  pd-$S^\sharp$-smooth. We call such pd-$S^\sharp$-thickenings {\it coordinate} ones.
   If $(Z,\CM )$ is log affine (see 1.1), then a coordinate thickening always exists, and an
 arbitrary $(Z,T)\in \CT_{S^\sharp}$  is pd-$S^\sharp$-smooth if and only if it is  retract of a coordinate one.\newline
(ii) Suppose  $(S,\CL )$ and $(Z,\CM )$ are fine, and
  $i_Y : (Z,\CM )\hra (Y,\CN )$ is a closed embedding of log $S^\sharp$-schemes with
 $(Y,\CM_Y )$ log smooth over $(S,\CL )$ (see \cite{K1} 3.3). Then its pd-$S^\sharp$-envelope
  is  pd-$S^\sharp$-smooth. 
 \newline (iii) For any integral quasi-coherent log affine $(Z,\CM )$ over $S^\sharp$  there is a universal $i$ as in (i): Take  $I=\Gamma (Z,\CO_Z )$,  $J=\Gamma (Z,\CM)$,  the embedding $i$ is the evident one. The corresponding $(P^{\text{univ}}, \CM_{P^{\text{univ}}}) $  depends on $(Z,\CM )$ in a functorial way.
  \endremark

\remark{Question} Is it true that property of being pd-$S^\sharp$-smooth is \'etale local?
\endremark

\subhead{\rm 1.5}\endsubhead {\it Log crystalline site.} For $S^\sharp$  as in 1.3 and  
a log $S^\sharp$-scheme $(Z,\CM )$ which is integral and quasi-coherent, objects of
the {\it  log crystalline site} $  ((Z,\CM)/S^\sharp )_{\crys} =(Z/S)^{\log}_{\crys}$ are pairs $(U,T)$ where $U$ is an \'etale $Z$-scheme, $(U,T)\in \CT_{S^\sharp}$ is
a pd-$S^\sharp$-thickening of $(U,\CM_U )$. The coverings are \'etale ones, i.e., collections of morphisms such that
the  maps of $T$'s form an \'etale covering.
The  structure sheaf $\CO_{Z/S}$ of $(Z/S)^{\log}_{\crys}$ is $\CO_{Z/S} (U,T ):=\Gamma (T,\CO_T )$; we have its pd-ideal $\CJ_{Z/S}$, $\CJ_{Z/S} (U,T ):=\Gamma (T,\CJ_T )$. One has an evident map of ringed sites $\iota : Z_{\et}\to (Z/S)^{\log}_{\crys} $, $\iota^{-1}((U,T)):=U$,
 $\CO_{Z/S}/\CJ_{Z/S} \iso \iota_* \CO_{Z_{\et}}$. The ringed site $(Z/S)^{\log}_{\crys}$ carries a log structure with the monoid sheaf  $\CM_{Z/S}$, $\CM_{Z/S}(U,T):=\Gamma (T,\CM_T )$, and $\iota$ is naturally a map of log ringed sites.
 There is a canonical morphism of topoi   $u^{\log}_{Z/S}: (Z/S)^{\log}_{\crys}\tilde{}\to Z\tilde{_{\et}}$, $u^{\log}_{Z/S *} (\CF )(U):= \Gamma ((U/S)^{\log}_{\crys},\CF )$, $u^{\log *}_{Z/S} (\CG )(U,T)= \CG (U)$.

\remark{Remark} 
For  $(Z,T)\in (Z/S)^{\log}_{\crys}$, any sheaf $\CF$ on $(Z/S)^{\log}_{\crys}$ yield naturally a sheaf $\CF_{(Z,T)}$  on $Z_{\et}=T_{\et}$ with $\Gamma (Z,\CF_{(Z,T)})=\CF (Z,T)$. The functor $\CF \mapsto \CF_{(Z,T)}$ is exact, so the  natural map $u^{\log}_{Z/S *} (\CF ) \to    \CF_{(Z,T)}$ yields
one $Ru^{\log}_{Z/S *} (\CF ) \to    \CF_{(Z,T)}$. \endremark

 \proclaim{Proposition} The category $(Z/S)^{\log}_{\crys}$ has non-empty finite inverse limits.
 \endproclaim
 
 \demo{Proof} It suffices to check that it has non-empty finite products and fiber products. 
  
 (a) For $(U_i ,T_i )\in (Z/S)^{\log}_{\crys}$, $i=1,2$, let us construct their product $(U,T)$. One has $U=U_1 \times_Z U_2$. To define $T$, consider the diagonal embedding $i_Y : (U,\CM_U )\hra (Y,\CN ):=(T_1 ,\CM_{T_1} )\times (T_2,\CM_{T_2} )$. The ideal of $\CO_Y$ generated by the pullbacks of  $\CJ_{T_i}+\CI \CO_{T_i}$ is a pd-ideal, so we have the log pd-scheme $Y^\sharp$. Let $(U,T')$  be  the pd-$Y^\sharp$-envelope of $i_Y$ (it is well defined since $(Y,\CN )$ is integral quasi-coherent).  Our $(T,\CM_T )$ is an exact closed log pd-subscheme of $(T',\CM_{T'})$ 
 defined as follows.  Let $p_i $ be the compositions $ (T',\CM_{T'})\to (T_i , \CM_{T_i})\to (S,\CL )$. Let $\CJ^0_{T'}$ be the ideal in $\CO_{T'}$ generated by local sections 
 $p_1^* (f)-p_2^* (f)$, $f\in \CO_S $, and $ (p_1^* ( \ell )/p^*_2 (\ell )) - 1$, $\ell \in \CL$. Here $p_i^* (\ell )$ are sections of $\CM_{T'}$ with the same image in $\CM_U$,
so their ratio is a section of $\CO^\times_{T'}$ that equals 1 on $U$.\footnote{Since $(U,\CM_U )\hra (T',\CM_{T'})$ is an exact embedding with nil ideal $\CJ_{T'}$,  $\CM_U$ is the quotient of $\CM_{T'}$ modulo the action of $1+\CJ_{T'} \subset \CO^\times_{T'}$, and this action is free since $\CM_{T'}$ is integral.} The ideal $\CJ^0_{T'}$ is 
 quasi-coherent since $(S,\CL )$ is quasi-coherent, and $\CJ^0_{T'}\subset \CJ_{T'}$. Now the ideal of $T$ in $\CO_{T'}$ is the pd-ideal generated by $\CJ^0_{T'}$. Since $p_i$ coincide on $(T,\CM_T )$, our $(U,T )$ is an object of $ (Z/S)^{\log}_{\crys}$. We leave it to the reader to check that it is the product of $(U_i ,T_i )$.

 (b) For morphisms $(U_i ,T_i )\to (V,Q)$, $i=1,2$, in $(Z/S)^{\log}_{\crys}$, let us construct their fiber product $(W, P)$. One has $W=U_1 \times_V U_2$. This is an open subset of $U_1 \times_Z U_2$; let  $(W, P')$ be the restriction to $W$ of the product of $(U_i ,T_i )$ (see (a)). Our $(W,P)$ is an exact closed log pd-subscheme of $(W, P')$ whose pd-ideal $\CJ'_{P'} \subset \CJ_{P'}$ is defined as follows. Let $q_i $ be the compositions $(P' ,\CM_{P'}) \to (T_i ,\CM_{T_i}) \to (Q,\CM_Q  )$.
 We can 
 work \'etale locally on $V$, so let us assume that $V$ is affine and $\CM_Q$ has  a chart $M\to \Gamma (Q,\CM_Q  )$. Now $\CJ'_{P'}$ is generated by sections $q_1^* (g)-q_2^* (g)$, $g\in \Gamma (Q,\CO_Q )$, and $ (q_1^* ( m )/q^*_2 (m )) - 1$, $m\in M$.
  \qed
 \enddemo

Let $f: (Z',\CM' )/S^{\prime\sharp}\to (Z,\CM )/S^\sharp$ be a map of the above data. A presheaf 
 $\CF'$ on $(Z'/S')^{\log}_{\crys}$ yields a presheaf  $f_{\crys\, *}(\CF' )$ on $(Z/S)^{\log}_{\crys}$ with
   $f_{\crys\, *}(\CF' )(U,T  ):=\Gamma (((U_{Z'},\CM'_{U_{Z'}} )/T^\sharp )_{\crys},\CF' )$, where $U_{Z'}:=U\times_Z Z'$, $T^\sharp :=(T,\CM_T, \CJ_T +\CI \CO_T )$. If $\CF'$ is a sheaf, then  $f_{\crys\, *}(\CF' )$ is a sheaf. 
   
\proclaim{Corollary} $f_{\crys\, *}$  defines  a morphism of  topoi $f_{\crys}\!: (Z'/S')^{\log}_{\crys}\tilde{}\to (Z/S)^{\log}_{\crys}\tilde{}$. 
\endproclaim

 \demo{Proof} For $(U',T')\in  (Z'/S')^{\log}_{\crys} $ let $C_f = C_f  (U',T')$ be the category of pairs $((U,T), g)$ where $(U,T)\in (Z/S)^{\log}_{\crys}$, $g: (U',T')\to (U,T)$ is a map of pd-$S^\sharp$-thickenings compatible with $f$. 
 The proposition implies that $C_f$ has finite inverse limits, so $C_f^\circ$ is directed.
  For a presheaf $\CF$ on $(Z/S)^{\log}_{\crys}$, its pullback $f^{-1}_{\crys} (\CF )$ assigns to  
  $(U',T')$ the colimit of the functor $((U,T), g)\mapsto \CF (U,T)$ on $C_f^\circ$.\footnote{ $C_f^\circ$ contains a small  cofinal subcategory, so the colimit is well defined.} If $\CF$ is a sheaf, then
 $f_{\crys}^* (\CF )$ is the sheaf associated with $f^{-1}_{\crys} (\CF )$. 
 Since  $C_f^\circ$ is directed,  $f^{-1}_{\crys} $ commutes with finite inverse limits, so same is true for $f_{\crys}^*$. We are done. \qed
 \enddemo
 
The evident maps $\CO_{Z/S}\to f_{\crys\, *}(\CO_{Z'/S'})$, $\CM_{Z/S}\to f_{\crys\, *}(\CM_{Z'/S'})$ make $f_{\crys}$ a morphism of log ringed topoi.

  \subhead{\rm 1.6}\endsubhead  {\it Log crystalline cohomology.}  For $ (U,T )\in (Z/S)^{\log}_{\crys}$, let $(U,T_* )$ be the restriction to $U \subset U\times_Z \ldots \times_Z U$ of  the standard simplicial object  of $ (Z/S)^{\log}_{\crys}$      with terms $(U,T_a ):=(U,T)^{[0,a]}$ (the product is computed in   $ (Z/S)^{\log}_{\crys}$, see Proposition  in 1.5).
The construction is natural and compatible with \'etale localization, hence
 any sheaf $\CF$ on $(Z/S)^{\log}_{\crys}$ yields a cosimplicial sheaf $(U,T) \mapsto \CF (U, T_*)$. If $\CF$ is a sheaf of abelian groups, then let
 $\CC^\cdot \CF$ be the normalized complex of this cosimplicial sheaf. The functor $\CC^\cdot$ is  exact. There is an evident projection $\alpha : \CC^\cdot \CF \to \CF$.
 
The embedding $H^0 \CC^\cdot \CF \hra \CF$ yields an isomorphism $
u^{\log}_{Z/S *} ( H^0 \CC^\cdot\CF ) \iso u^{\log}_{Z/S *} ( \CF )$, so we have a natural map $\beta : u^{\log}_{Z/S *} ( \CF )\to u^{\log}_{Z/S *} ( \CC^\cdot\CF )$ right inverse to $u^{\log}_{Z/S *}(\alpha )$.

Suppose we have  $(Z,T)\in (Z/S)^{\log}_{\crys}$. Due to exactness of $\CC^\cdot$,  
 the natural map $u^{\log}_{Z/S *} ( \CC^\cdot\CF )\to \CC^\cdot\CF_{(Z,T)}$  
yields one $R(u^{\log}_{Z/S *} \CC^\cdot )( \CF )\to \CC^\cdot \CF_{(Z,T)}$. 
 
\proclaim{Proposition} (i) One has $R(u^{\log}_{Z/S *} \CC^\cdot )(\CF )\iso Ru^{\log}_{Z/S *} (\CF )$. \newline
(ii) If  $(Z,P)  \in (Z/S)^{\log}_{\crys} $  is pd-$S^\sharp$-smooth, then  $R(u^{\log}_{Z/S *} \CC^\cdot )(\CF )\iso \CC^\cdot \CF_{(Z,P)}$. Thus $$Ru^{\log}_{Z/S *} (\CF )\iso \CC^\cdot \CF_{(Z,P)}.\tag 1.6.1$$ 
 \endproclaim
 
 \demo{Proof} Deriving $\alpha$ and $\beta$, we get $R(u^{\log}_{Z/S *} \CC^\cdot )(\CF )\leftrightarrows Ru^{\log}_{Z/S *} (\CF )$ whose composition in one direction is identity; we want to check they are mutually inverse. The problem is local, so
 we can assume  that there is $P$ as in (ii) and $Z$ is affine. Since $(Z,P)$ is pd-$S^\sharp$-smooth, one has $u^{\log}_{Z/S *} (\CF )\iso H^0 \CC^\cdot\CF_{(Z,P)}$. To prove (i), (ii), it suffices then to
 find for every $\CF$  an embedding $\CF \hra \tilde{\CF}$ such that $H^{\neq 0}\CC^\cdot\tilde{\CF}_{(Z,P)}=0$. 
 
For $(Z,T)\in (Z/S)^{\log}_{\crys}$ consider the simplicial object $(Z,P_* )\times (Z,T)$ augmented over $(Z,T)$. The augmentation admits an inverse up to homotopy: indeed, 
any map $(Z,T)\to (Z,P)$ yields one in the usual way, and such a map exists  since $(Z,P)$ is pd-$S^\sharp$-smooth. So the sheaf $(U,T')\mapsto \tilde{\CF}(U,T'):=\CF ((U,T')\times (Z,T))$ satisfies $H^{\neq 0}\CC^\cdot\tilde{\CF}_{(Z,P)}=0$. If
$(Z,T)$ is pd-$S^\sharp$-smooth (say, $(Z,T)=(Z,P)$), then
the evident map $\nu : \CF \to \tilde{\CF}$ is injective.\footnote{Indeed, for $U$ affine one can find a map $(U,T')\to (Z,T)$, which yields a left inverse to $\nu$.} We are done. \qed
 \enddemo
 
\remark{Remark} (i) Suppose $(Z,P)$ is a coordinate thickening as in  Remark (i) in 1.4.
The $\Bbb G_a^I \times \Bbb G_m^J$-action on
$\Bbb A^{I}_{(S,\CL )}\times_{(S,\CL )}\Bbb A^{(J)}_{(S,\CL )}$ (see 1.1)\footnote{An element $((h_i ),(g_j))_{i\in I, j\in J} \in \Bbb G_a^I \times \Bbb G_m^J$ acts as $ t_i \mapsto t_i - h_i$,  $ t_j \mapsto g^{-1}_i t_j$.} yields, by universality, a  $\Bbb G_a^{\sharp I} \times \Bbb G_m^{\sharp J}$-action on $(Z,P)$ (see 1.2). Then $(Z,P)$ is a  $\Bbb G_a^{\sharp I} \times \Bbb G_m^{\sharp J}$-torsor on $(Z/S)^{\log}_{\crys}$. Thus $(Z,P_a )= (\Bbb G_a^{\sharp I} \times \Bbb G_m^{\sharp J})^a \times (Z,P )$ and $(Z,P_* )$ is the  ``universal  simplicial $\Bbb G_a^{\sharp I} \times \Bbb G_m^{\sharp J}$-quotient"  for the $\Bbb G_a^{\sharp I} \times \Bbb G_m^{\sharp J}$-action on $(Z,P)$. \newline
(ii) If $(Z,P)$ is pd-$S^\sharp$-smooth, then it is a covering of $(Z/S)^{\log}_{\crys}$. \newline
(iii)  If  no pd-$S^\sharp$-smooth  $(Z,P)$  is available, one can compute   $Ru^{\log}_{Z/S *} (\CF )$ by either of the next procedures: \newline -  By Remark (i) in 1.4 and Proposition in 1.5, any \'etale hypercovering $p :U_\cdot \to Z$ with log affine $(U_0 ,\CM_{U_0})$  has a pd-$S^\sharp$-smooth thickening $(U_\cdot ,P_\cdot )$, so $Ru^{\log}_{Z/S *} (\CF )\iso Rp_{ *} \CC^\cdot \CF_{(U_\cdot ,P_\cdot )}$. \newline
 - Replace $Z_{\et}$ by an equivalent topology $Z_{\et'}$ formed by those   \'etale $U/Z$ that $(U,\CM_U )$ is log affine.  Then $U\mapsto \Gamma (U,\CC^\cdot \CF_{(U,P^{\text{univ}})})$ (see Remark (iii) in 1.4) is a presheaf on $Z_{\et'}$; its sheafification equals $Ru^{\log}_{Z/S *} (\CF )$. \newline -  Consider 
for each $U\in X_{\et}$ the category $\CS (U)$ of its pd-$S^\sharp$-smooth thickenings. Then\footnote{To define holim, one should take  care of the usual set-theoretic difficulties. } the presheaf
$ U \mapsto \text{holim}_{(U,P)\in \CS (U)}R \Gamma (U,  \CC^\cdot \CF_{(U,P)})$  represents $Ru^{\log}_{Z/S *} (\CF )$.
\endremark

\subhead{\rm 1.7}\endsubhead {\it Log $\CO$-crystals and  connections.} For $(U,T)\in (Z/S)^{\log}_{\crys}$ its  {\it de Rham pd-algebra}
 $\Omega^\cdot_{(U,T)/S}$ is the quotient of the log de Rham dg algebra $\Omega^\cdot_{(T,\CM_T )/(S,\CL )}$ (see \cite{K1} 1.7, 1.9)  modulo the relations $d(u^{[n]})=u^{[n-1]}du$, $u\in \CJ_T $. This is a sheaf of commutative dg $\CO_S$-algebras on $T_{\et}=U_{\et}$ whose terms are quasi-coherent $\CO_T$-modules (since $(T, \CM_T )$ and $(S,\CL )$ are quasi-coherent).  It carries 
 the {\it Hodge pd-filtration} $F^m $, $F^m \Omega^{ a}_{(Z, T  ) /S}:= \CJ_T^{[m-a]} \Omega^{ a}_{(Z, T  ) /S}$; one has $F^m \cdot F^\ell \subset F^{m+\ell}$. We get 
a sheaf    $\Omega^\cdot_{Z/S}$  of filtered commutative dg algebras on $(Z/S)^{\log}_{\crys}$, $\Omega_{Z/S}^{\cdot}(U,T ):= \Gamma (T, \Omega^{ \cdot}_{(U,T) /S})$. 

\remark{Exercises} (i) If $(U,T)$ is the pd-envelope of $(U,\CM) \hra (Y,\CN )$ as in (1.3.1), then
the composition $ r^* \Omega^1_{(Y,\CN )/(S,\CL )} \buildrel{r^*}\over\to \Omega^1_{(T,\CM_T)/(S,\CL )}\twoheadrightarrow  \Omega^1_{(U,T)/S}$ is an isomorphism. \newline (ii) If $(U,T)$ is a pd-smooth thickening and $U$ is affine, then $\Omega^1_{Z/S}(U,T)$ is a projective $\Gamma (T,\CO_T )$-module.\footnote{Our assertion is local, so we can assume that $(Z,\CM )$ is log affine; by Remark (i) in 1.4, we can assume that $(Z,T)$ is a coordinate thickening; now the assertion follows from (i).}
\endremark

\enspace

Here is another description of $\Omega^\cdot_{Z/S}$. For $(U,T)\in (Z/S)^{\log}_{\crys}$ let $(U,T_* )$ be be the simplicial object defined in 1.6. Let 
$\CV_{T_1}\subset \CO_{T_1}$ be the ideal of the subscheme $T=T_0$ of $T_1$, $\CV_{T_1}^{[2]}  \subset \CO_{T_1}$ be its divided powers square. Let  $T^\flat_* $ be the closed exact simplicial log subscheme of $T_*$ whose cosimplicial ideal in $\CO_{T_*}$ is generated by $\CV_{T_1}^{[2]}$. 
The normalized complex N$^\cdot (\CO_{T^\flat_*})$ equipped with 
the Alexander-Whitney product is a sheaf of dg algebras on $U_{\et}$ that  depends functorially on $(U,T)$. Thus we have a sheaf $\Omega^{\flat\cdot}_{Z/S} $ of dg algebras on $(Z/S)^{\log}_{\crys}$, $\Omega_{Z/S}^{\flat\cdot}(U,T):=\Gamma(U, \text{N}^\cdot (\CO_{T^\flat_*}))$.

\proclaim{Proposition} One has a natural  identification of dg algebras 
 $\phi^\cdot : \Omega^\cdot_{Z/S}\iso \Omega_{Z/S}^{\flat\cdot}
$. \endproclaim

\demo{Proof} We will construct a natural isomorphism $\varphi^\cdot : \Omega^\cdot_{(U,T)/S}\iso \text{N}^\cdot (\CO_{T^\flat_*})$
of sheaves of dg algebras on $U_{\et}$. It  satisfies
 the  properties  $\phi^0 =\id_{\CO_T}$ and  $\phi^1 (d\log (m))=(p^*_1 (m)/p^*_0 (m))-1\in \CV_{T_1}/\CV_{T_1}^{[2]}= $ N$^{1} (\CO_{T^\flat_*})$ for $m\in\CM_T$ (here $p_i : T^\flat_1 \to T$ are the projections) that determine $\varphi^\cdot$ uniquely.\footnote{Since $\Omega^\cdot_{(U,T)/S}$ is generated, as a dg algebra, by $\Omega^0_{(U,T)/S}$ and $d\log (\CM_T )\subset \Omega^1_{(U,T)/S}$.} 

(a) Recall that, by Dold-Puppe, the normalization functor N$^\cdot$ is an equivalence between the category of cosimplicial abelian groups and the category of complexes vanishing in negative degrees; let K$^*$ be the inverse equivalence. Both categories carry the usual symmetric tensor products $\otimes$. The functors  N$^\cdot$ and K$^*$ transform algebras to algebras (using the Alexander-Whitney and shuffle products respectively); both send associative algebras to associative ones, and  K$^*$ transforms commutative algebras to commutative ones. If $C^\cdot$ is a dg algebra, then
the identification of complexes N$^\cdot$K$^* (C^\cdot )\iso C^\cdot$ is compatible with the products.

Some explicit formulas: Below $\partial_i : [0,n-1]\to [0,n]$, $\sigma_i : [0,n]\to [0,n-1]
$ are the standard face and degeneration maps, $\sigma_i \partial_i = \sigma_i \partial_{i+1}=\id_{[0,n-1]}$.
Let $A^* $ be a cosimplicial abelian group and $C^\cdot $ be a complex that correspond one to another  by N$^\cdot$ and $K^*$. For a monotone map
$e : [0,m]\to [0,n]$ we write the cosimplicial structure map $e=e_A : A^m \to A^n$ also as $a\mapsto {}^e a$. One has  $C^n =\cap\, \Ker \sigma_i \subset A^n$, and $d=\Sigma (-1)^i \partial_i |_C$.
For $m\in [0, n]$ let $E(m,n)$ be the set of increasing injections $e : [0,m]\hra [0,n]$ such that $e (0)=0$.\footnote{The map $e\mapsto e([1,m])$ identifies $E(m,n)$ with the set of all order $m$ subsets of $[1,n]$.} Let ${}^e C^m \subset A^n$ be 
the image of $C^m \subset A^m$ by the 
(injective) map $e_A$.
One has a Dold-Puppe  direct sum decomposition $$A^n =\oplus_{m\le n}\oplus_{e\in E(m,n)} {}^e C^m , \tag 1.7.1$$ 
For a monotone  $g: [0,n]\to [0,\ell ]$ the   components of the map $g_A : A^n \to A^\ell$ are as follows: For
 ${}^e c \in {}^e C^m$ one has $g_A ({}^e c )=0$ if $ge: [0,m]\to [0,\ell ]$ is not injective, and $g_A ({}^e c )={}^{ge} c$ if $ge \in E(m,\ell)$.  Otherwise we have $[ge]\in E(m+1,\ell )$ such that $[ge]\partial_0= ge$, and $g_A ({}^e c)={}^{[ge]}(d c )- \Sigma_{1\le i\le m+1}(-1)^i \,\,\, {}^{[ge]\partial_i } c$. Here $\partial_i : [0,m]\to [0,m+1]$ are the usual face maps.

If $C^\cdot $ is a  dg algebra, then the corresponding shuffle product on $A^*$ looks as follows. For $e_1 \in E(m_1 ,n)$, $e_2 \in E(m_2 ,n)$, $c_1 \in C^{m_1}$, $c_2 \in C^{m_2}$ 
the shuffle product ${}^{e_1} c_1  {}^{e_2}c_2$ vanishes if $e([1,m_1 ])e_2 ([1,m_2 ])\neq\emptyset$; otherwise it equals $ \pm {}^{e} (c_1 c_2 )$ where $e \in E(m_1 +m_2 ,n)$ has image $\{0\} \cup e_1 ([1,m_1 ])\cup e_2 ([1,m_2 ])$ and $\pm$ is the sign of
the permutation $\sigma$ of $[1,m+m']$ such that $\sigma (a)$ equals $e^{ -1}e_1 (a)$ if $a\le m_1 $, and  $e^{ -1}e_2 (a-m)$ if $a> m_1$. Notice that the product on each ${}^e C^m$ vanishes if $m\neq 0$.

Suppose $C^\cdot$ is a strictly\footnote{Which means that $c^2 =0$ for $c$ of odd degree.} commutative dg algebra, so the shuffle product on $A^*$ is commutative. Let
 $V^n $ be the kernel of the degeneration map $A^n \twoheadrightarrow A^0$. Then $V^*$ is a cosimplicial ideal in $A^*$. It carries a unique pd structure such that  the cosimplicial structure maps are pd-morphisms and  $ c^{[k]}=0$ for  $c\in C^m \subset A^m$, $k>1$, $m>0$.\footnote{One has $V^n =\oplus_{m>0}\, {}^e C^m$ as in (1.7.1). For $v\in V^n$ let $v_{m,e}\in {}^e C^m$ be the components of $v$. If we have our pd-structure, then $v_{m,e}^{[>1]}=0$ implies that
 $v^{[k]}=\mathop\Sigma\limits_P \, \mathop\Pi\limits_{(m,e)\in P}  v_{(m,e)}$, where $P$ runs the set of all $k$ element subsets of $\sqcup_{m>0}\, E(m,n)$. Now take this formula as the definition of
 divided powers on $V^n$. The axioms of pd-structure are clear. To show that $g_A : A^n \to A^\ell$ are pd-morphisms, it suffices
to check that  $(g_A (v_{m,e} ))^{[k]}=0$ for $k>1$. This is immediate for $m>1$ and follows from strict commutativity for $m=1$. } More generally, suppose we have a pd-ideal $J\subset C^0 =A^0 $. Then the pd structures on $J$ and $V^*$ are compatible (i.e., extend to the pd structure on the cosimplicial ideal $J_A^*$ they generate that is compatible with the cosimplicial structure maps) if and only if for every $f\in J$, $n>0$ one has $d(f^{[n]})=f^{[n-1]}d(f)$.\footnote{Notice that
$J_A^n = {}^{\delta_0}\! J \oplus V^n$, where ${}^{\delta_0}\! J\subset {}^{\delta_0} C^0$,  $\{ \delta_0 \}=E(0,n)$. So we set $({}^{\delta_0}\! f +v)^{[k]}:=\Sigma_a {}^{\delta_0}\! f^{[a]} v^{[k-a]}$. One easily checks the pd structure axioms. The compatibility with cosimplicial structure is enough to check for the map $\partial_{0 A} : J \to J_A^1$. Then $\partial_{0 A}(f)^{[n]}= ({}^{\delta_0}\!f+ df )^{[n]}= {}^{\delta_0}\! f^{[n]}+ f^{[n-1]}df$ and $\partial_{0 A}(f^{[n]})={}^{\delta_0}\! f^{[n]}+d(f^{[n]})$, and we are done.}

\enspace

(b) Set $\Omega^\cdot := \Omega^\cdot_{(U,T)/S}$, and consider the cosimplicial commutative algebra $A^* 
 :=$ K$^* (\Omega^\cdot )$. By (a), the pd structure
 on $\CJ_T \subset \CO_T =A^0$ yields a pd structure on
 the corresponding cosimplicial ideal $J^*_A$; notice that $A^n /J^n_A = \CO_T /\CJ_T =\CO_U$.  The log structure $\CM_T$  extends naturally to a log structure on the simplicial scheme $\Spec\, A^*$: We define
 $\CM_{\Spec\, A^n}$ as the pullback of $\CM_T$ by either of
 the $n+1$ simplicial structure projections $p_i =\Spec (\delta_{i A}): \Spec\, A^n \to \Spec\, A^0 =T$, where  $\delta_i : [0,0]\to [0,n]$ is $\delta_i (0)=i$. The pullbacks are  identified as $p_i^* (m)= (1- \delta_{i,j A}(d\log (m)))p_j^* (m)$ for $i<j$, $m\in\CM_T$, where $\delta_{i,j}:[0,1]\to [0,n]$ is $\delta_{i,j}(0)=i$, $\delta_{i,j}(1)=j$, and $d\log (m)\in \Omega^1 \subset A^1$. The log and pd structures  make $(U, \Spec\, A^* )$  a simplicial object of  $(Z/S)^{\log}_{\crys}$. By the definition of $T_*$, there is a unique map $(U,\Spec\, A^* )\to (U,T_* )$ in  $(Z/S)^{\log}_{\crys}$ that equals $\id_T$ in degree 0. It evidently takes image in $ (U,T^\flat_* )$.

 We have defined a morphism of cosimplicial algebras $\CO_{T^\flat_*} \to A^*$, hence a morphism of dg algebras  $\psi^\cdot :$ N$^\cdot (\CO_{T^\flat_*})\to $ N$^\cdot (A^* )=\Omega^\cdot$. It remains to show that $\psi^\cdot$
is an isomorphism: the promised identification  $\varphi^\cdot$ is its inverse.

\enspace 

(c) Below N$^\cdot :=$ N$^\cdot (\CO_{T^\flat_*})$. Let us first check
 that $\psi^1 :$ N$^1 \to \Omega^1$ is an isomorphism. We have $T\hra T^\flat_1$ and the two retractions $p_0 ,p_1 : T^\flat_1 \to T$. Our
 N$^1 $ is the ideal of $T$ in $ \CO_{T^\flat_1}$; it  has square zero; $p_0$ yields a splitting $\CO_T\, \oplus $
N$^1 \iso\CO_{T^\flat_1}$ and an identification of log structures $p_0^* (\CM_T )\iso \CM_{T^\flat_1}$. The other retraction $p_1$ amounts then to an $\CO_T$-linear map $\phi^1 : \Omega^1_{(T,\CM_T)/(S,\CL )}\to$ N$^1 $ such that $\phi^1 (d(f))=p_1^* (f)-p_0^* (f)$, $\phi^1 (d\log (m))= p_1^* (m)/p_0^* (m) -1$. Since $p_i$ are pd-morphisms and divided powers of degree $>1$ vanish on N$^1 $,  $\phi^1$ factors through $\Omega^1$ (cf.~the computation in the last footnote in (a)). Since 
the images of $ \CO_T$ and $\CM_T$ by $p^*_0$ and $p_1^*$ generate $\CO_{T^\flat_1}$,
$\phi^1$ is surjective. Since $ \psi^1 \phi^1 =\id_{\Omega^1}$ by
 the construction of $\psi^1$,  $\phi^1$ is inverse to $\psi^1$.

\enspace

(d) For $n\ge 2$, let $\sigma_i : \CO_{T^\flat_n} \twoheadrightarrow  \CO_{T^\flat_{n-1}}$, $i\in [0,n-1]$, be the standard degeneration maps, $\sigma_i \partial_i =\sigma_i \partial_{i+1}= \id_{[0,n-1]}$. Then the kernel $\CV_{T_n^\flat}$ of the degeneration map $\CO_{T^\flat_n}\twoheadrightarrow \CO_T$ equals $\Sigma_i \, \Ker (\sigma_i )$, N$^n =\cap_i \Ker (\sigma_i )$, and
 $\Ker (\sigma_i )$ equals the ideal in  $\CO_{T^\flat_n}$ generated by $\delta_{i,i+1}($
  N$^1 )$. Thus $ \Ker (\sigma_i )^{[2]}=\CV_{T_n^\flat}$N$^n = ($N$^n )^{[2]}=  0$.
 
For $c\in $ N$^1$ one has $c\cup c=0$ (here $\cup$ is the Alexander-Whitney product).\footnote{
Indeed, $c\cup c:=\partial_2 (c)\partial_0 (c)=(\partial_2 (c)+\partial_0 (c))^{[2]}= (\partial_1 (c)+d(c))^{[2]}=0$, since $\partial_i (c)^{[2]} =\partial_1 (c) d(c)=d(c)^{[2]}=0$.} Thus the subalgebra of N$^\cdot $ generated by N$^{\le 1} $ is strictly commutative. This subalgebra is closed under differential,  since N$^{1} $ is generated, as an $\CO_T$-module, by cycles.\footnote{Namely, $d(f)$ and  $\phi^1 (d\log (m))$, where  $f\in \CO_T$, $m\in\CM_T$, see (c).} By the universality property of 
$\Omega^\cdot$ and (c), there is a unique map of dg algebras $\phi^\cdot : \Omega^\cdot \to$ N$^\cdot $ such that $\psi^\cdot \phi^\cdot = \id_{\Omega^\cdot}$. 

It remains to prove that $\phi^\cdot$ is surjective, or, equivalently, that the image $A^*$ of $K^* (\phi^\cdot )$ equals $\CO_{T^\flat_*}$. Since $\CO_{T^\flat_*}$ is generated, as cosimplicial algebra, by $\CO_{T^\flat_{ 1}}$, it suffices to show that $A^*$ is a subalgebra of  $\CO_{T^\flat_*}$. Notice that $A^n$, $n>0$, is generated, as an abelian group, by elements $\delta_0 (f)$, $f\in \CO_T$, and $e_A (\nu_1 \cup \ldots \cup \nu_m )$, $\nu_i \in \Omega^1$, $e\in E(m,n)$, $0<m\le n$.
One has $\delta_0 (f)e_A (\nu_1 \cup\ldots\cup \nu_m )=e_A (f\nu_1 \cup\ldots\cup \nu_m )$.
 Since $\delta_{0,j}(\nu)\delta_{0,j}(\nu')\in \delta_{0,j}((\text{N}^1)^2)=0$, one has $\nu_1 \cup \ldots \cup \nu_m := \delta_{0,1}(\nu_1)\delta_{1,2}(\nu_2 )\ldots \delta_{m-1 ,m}(\nu_m )= \delta_{0,1}(\nu_1)\delta_{0,2}(\nu_2 )\ldots \delta_{0 ,m}(\nu_m )$.\footnote{Indeed, by induction by $m$, one has     $\nu_1 \cup \ldots \cup \nu_m = \delta_{0,1}(\nu_1)\ldots \delta_{0 ,m-1}(\nu_{m-1} ) \delta_{m-1 ,m}(\nu_m )$. Since  $\delta_{m-1 ,m}(\nu_m )-\delta_{0,m}(\nu_m )= \delta_{0,m-1}(\nu_m )+ e_A (d(\nu_m )) $, where $e\in E(2,n)$, $e(1)=m-1$, $e(2) =m$, and
 $e_A (d(\nu_m ))=\Sigma\, \delta_{0,m-1}(\nu'_i )\delta_{m-1,m}(\nu''_i)$, we see that  $\delta_{0 ,m-1}(\nu_{m-1} ) \delta_{m-1 ,m}(\nu_m )=  \delta_{0 ,m-1}(\nu_{m-1} ) \delta_{0 ,m}(\nu_m )$, hence the assertion.} Thus  $e_A (\nu_1 \cup \ldots \cup \nu_m ) e_{A'} (\nu'_1 \cup \ldots \cup \nu'_{m'} )=\delta_{0,e(1)}(\nu_1)\ldots \delta_{0 ,e(m)}(\nu_m ) \delta_{0,e'(1)}(\nu'_1)\ldots \delta_{0 ,e'(m')}(\nu'_{m'} )$ vanishes if $e([1,m])\cap e'([1,m'])\neq \emptyset$ and otherwise  equals $\pm e''_A (\nu_1 \cup \ldots \cup \nu_m \cup \nu'_1 \cup \ldots \cup \nu'_{m'} )$, where $e'' \in E(m+m', n)$, $e''([1,m+m'])=e([1,m])\cup e'[1,m'])$. We are done.
  \qed \enddemo

A sheaf $\CF$ of $\CO_{Z/S}$-modules on  $(Z/S)^{\log}_{\crys}$ is said to be {\it  $\CO_{Z/S}$-crystal} if for every morphism $\phi : (U',T' )\to (U,T)$ in 
$(Z/S)^{\log}_{\crys}$ the pullback map $\phi^{-1}\CF_{(U,T)}\to \CF_{(U',T')}$
yields an isomorphism of $\CO_{T'}$-modules $\phi^{*}\CF_{(U,T)}\iso \CF_{(U',T')}$. For such an $\CF$, let $\Omega^\cdot_{Z/S} \CF$ be the normalization of the cosimplicial sheaf $(U,T)\mapsto \CF (U,T^\flat_* )$. By Proposition, 
$\Omega^\cdot_{Z/S} \CF$ is a dg 
$\Omega^\cdot_{Z/S}$-module and
 $\Omega^i_{Z/S} \CF = \Omega^i_{Z/S}\otimes_{\CO_{Z/S}} \CF$. So $\Omega^\cdot_{Z/S} \CF$ is the de Rham complex for the flat connection $\nabla :=d^0 : \CF \to \Omega^1_{Z/S}\otimes_{\CO_{Z/S}} \CF$.
 
 \proclaim{Theorem} Suppose $(Z,P )$ is a pd-$S^\sharp$-smooth thickening. Then the connection $\nabla_{(Z,P)}: \CF_{(Z,P)} \to \Omega^1_{(Z,P)/S}\otimes_{\CO_P}\CF_{(Z,P)}$ is quasi-nilpotent (see \cite{K1} 6.2). The functor $\CF \mapsto (\CF_{(Z,P)},\nabla_{(Z,P)})$ is an equivalence between the category of $\CO_{Z/S}$-crystals and that of $\CO_P$-modules equipped with an integrable quasi-nilpotent connection.
 \endproclaim

\demo{Proof} This is theorem 6.2 from \cite{K1} (it is stated in loc.~cit.~under the assumption that $\CM$ is fine; the proof works 
in our setting as well). \qed
\enddemo

\subhead{\rm 1.8}\endsubhead {\it Comparison with the de Rham cohomology.} Below we call  $Ru^{\log}_{Z/S *} (\CO_{Z/S} )$ and $R\Gamma ((Z/S)^{\log}_{\crys},\CO_{Z/S} )= R\Gamma (Z_{\et}, Ru^{\log}_{Z/S *} (\CO_{Z/S} ))$ simply {\it the log crystalline complexes}. These are E$_\infty$ algebras.\footnote{I.e.,  dg algebras which are commutative up to a coherent system of higher homotopies.  See \cite{HS} for a nice, if old-fashioned, initial exposition (which explains, in particular, the E$_\infty$ algebra structure on 
the log crystalline complexes defined explicitly using the Godement resolution).}
  Suppose $\CF$ is an $\CO_{Z/S}$-crystal, so
  $Ru^{\log}_{Z/S *} (\CF )$ is a $Ru^{\log}_{Z/S *} (\CO_{Z/S} )$-module, $R\Gamma ((Z/S)^{\log}_{\crys},\CF)$ is an
$R\Gamma ((Z/S)^{\log}_{\crys},\CO_{Z/S})$-module.

\proclaim{Theorem} (i) The evident projection $\Omega^\cdot_{Z/S} \to \CO_{Z/S}$ yields a quasi-isomorphism $Ru^{\log}_{Z/S *} (\Omega^\cdot_{Z/S}\CF )\iso Ru^{\log}_{Z/S *} (\CF )$.
\newline (ii) For any pd-$S^\sharp$-smooth
 $(Z,P)\in   (Z/S)^{\log}_{\crys} $ the natural map  $Ru^{\log}_{Z/S *} (\Omega^\cdot_{Z/S}\CF )\to (\Omega^\cdot_{Z/S}\CF )_{(Z,P)} 
 $ is a quasi-isomorphism. Therefore (cf.~\cite{K1} 6.4)
 $$Ru^{\log}_{Z/S *} (\CF )\iso (\Omega^{\cdot}_{Z/S}\CF )_{(Z,P)} . \tag 1.8.1$$
\endproclaim

\demo{Proof} The assertions are local, so we can assume that we have $(Z,P)$ as in (ii) and $(Z,\CM )$ is log affine (see 1.1).  
Pick a pd-$S^\sharp$-smooth $(Z,T)$ (say, a copy of $(Z,P)$).
Let $(Z,Q_{*})$ be the product of $(Z,T_* )$ and $(Z,P )$; this is a simplicial object of $(Z/S)^{\log}_{\crys}$ augmented over $(Z,P)$. Consider the total complex of the cosimplicial complex  $(\Omega^\cdot_{Z/T_*}\CF)_{(Z,Q_* )}
$. Let us show that the pullback maps are quasi-isomorphisms: 
$$(\CC^\cdot \CF)_{(Z,T)}= \CF_{(Z,T_* )} \iso  (\Omega^\cdot_{Z/T_*}\CF)_{(Z,Q_{* })} \buildrel{\sim}\over\leftarrow (\Omega^\cdot_{Z/S}\CF )_{(Z,P)}. \tag 1.8.2$$

First arrow: Let us check that  the maps $\CF_{(Z,T_i )} \to  (\Omega^\cdot_{Z/T_i }\CF)_{(Z,Q_i )}$ are quasi-isomorphisms. By Remark (i) in 1.4, $(Z,P)$ is a retract of a coordinate pd-$S^\sharp$-thickening $(Z,P')$. Let $(Z,Q'_i)$ be the product of $(Z,T_i )$ and  $(Z,P')$. Then $\Omega^\cdot_{(Z,Q_i )/T_i}$ is a retract of $\Omega^\cdot_{(Z,Q'_i)/T_i}$, so the assertion  for $\Omega^\cdot_{(Z,Q'_i )/T_i}$ implies that for $\Omega^\cdot_{(Z,Q_i )/T_i}$ (since   the map $\CF_{(Z,T_i )} \to H^0 (\Omega^\cdot_{Z/T_i }\CF)_{(Z,Q_i )}$ is injective). By Remark (i) in 1.6, $(Z,Q'_i)$ is a $\Bbb G_m^{\sharp I}\times\Bbb G_a^{\sharp J}$-torsor over $(Z,T_i)$, and we are done by the evident computation of the de Rham pd-complexes of $\Bbb G_m^\sharp $ and $ \Bbb G_a^\sharp$.

 Second arrow: Let us check that the maps  $(\Omega^a_{Z/S}\CF)_{(Z,P)}\to (\Omega^a_{Z/T_*}\CF)_{(Z,Q_{*})}$ are quasi-isomorphisms. For case $a=0$, see the proof of the proposition in 1.6; the general case follows by base change since $\Omega^a_{(Z,P)/S}$ is $\CO_{P}$-flat by Exercise (ii) in 1.7.

Since diagrams (1.8.2) are compatible with maps between $P$'s, we see that
the simplicial structure maps $(\Omega^\cdot_{Z/S}\CF)_{(Z,P_i)}\to(\Omega^\cdot_{Z/S}\CF)_{(Z,P_j)}$ are quasi-isomorphisms. Now both (ii) and (i) of the theorem follow from  the proposition in 1.6, q.e.d. \qed
\enddemo

\remark{Remarks} 
(i) If no global $(Z,P)$ as in above is available, then one can compute  $Ru^{\log}_{Z/S *} (\CF )$ using (1.8.1) and Remark (iii) in 1.6. \newline (ii) If $(S,\CL )$ is fine and $(Z,\CM )$ is log smooth over $(S,\CL )$, then (1.8.1) implies that $$Ru^{\log}_{Z/S *} (\CF )
\iso (\Omega^{\cdot}_{Z/S}\CF)_{(Z,\CM)}. \tag 1.8.3$$ (iii) The map $\Omega^\cdot_{Z/S}\to \CO_{Z/S}$ sends $F^m 
\Omega^\cdot_{Z/S}\CF$ to $\CJ_{Z/S}^{[m]}\CF$ (see 1.7 for the notation). The quasi-isomorphisms in the statement of the proposition are, in fact, filtered quasi-isomorphisms for these filtrations. The proof is the same - just replace ``quasi-isomorphism" in it by ``filtered quasi-isomorphism".
\newline (iv) Quasi-isomorphism (1.8.1) coincides with the composition of maps in (1.8.2) and, via (1.6.1) and 1.7, with (the normalization of) the restriction map of cosimplicial sheaves $\CF_{P_*}\to \CF_{P_*^\flat}$.  
\endremark

\subhead{\rm 1.9}\endsubhead  {\it Comparison with derived de Rham cohomology.}   We discuss log version of  Illusie's comparison theorem \cite{Ill2} Ch.~VIII, 2.2.8 due to Olsson \cite{Ol1} 6.10. 

Let   $ L\Omega^\cdot_{(Z,\CM ) /(S,\CL)}$ be the  derived log de Rham complex (we use Gabber's construction,  \cite{Ol1} \S 8 or \cite{Bh2} \S 6, to be recalled in a moment; see \cite{B} 3.1 for a short review and the notation used below). This is a commutative dg $\CO_S$-algebra on $Z_{\et}$ equipped with the Hodge filtration $F^m$.  Let us define a natural morphism of filtered commutative dg $\CO_S$-algebras $$  L\Omega^\cdot_{(Z,\CM ) /(S,\CL)}  \to Ru^{\log}_{Z/S *} ( \CO_{Z/S} ). \tag 1.9.1$$ 

Assume for simplicity that $S=\Spec \, A$, $\CL$ comes from a prelog structure $L\to A$.\footnote{Otherwise do the construction \'etale locally on $S$.}    For a log affine $U\in Z_{\et}$  set $B:= \Gamma (U,\CO_U )$, $M:=\Gamma (U,\CM_U )$, so we have 
  $(B,M)/(A,L)\in\CC_{(A,L)}$. Let $P(U)_\cdot =P_{(A,L)}(B,M)_\cdot$ be its canonical simplicial resolution (see \cite{Ol1} 8.3), and $\Omega^\cdot_{P(U)_\cdot /(A,L)}$ be the relative log de Rham complex, which is a simplicial dg algebra. Let $L\Omega^\cdot_{(B,M)/(A,L)}$ be the total complex, 
$L\Omega^n_{(B,M)/(A,L)}=\oplus_i  \, \Omega^{i}_{P(U)_{i-n} /(A,L)}  $, filtered 
by the Hodge filtration $F^m  :=\oplus_{i \ge m}\,  \Omega^{i}_{P(U)_{i-n} /(A,L)}$. We have the filtered complex of presheaves
    $U\mapsto L\Omega^\cdot_{(B,M)/(A,L)}
    $  on  $Z_{\et'}$; the associated filtered complex of sheaves is quasi-isomorphic to
 $ L\Omega^\cdot_{(Z,\CM ) /(S,\CL)}$.

 Let $(U,T^{a})$ be the pd-$S^\sharp$-completion of the embedding $(U,\CM_U ) \hra \Spec\, P(U)_a$;\footnote{Notice that $(U,T^0 )$ is $(U,P^{\text{univ}})$ from Remark (iii) in 1.4.} this is a cosimplicial object in $(Z/S)^{\log}_{\crys}$. Each 
$(U,T^{ a})$ is pd-$S^\sharp$-smooth  (see Remark (i) in 1.4), so one has the filtered quasi-isomorphism $Ru^{\log}_{Z/S *}(\CO_{Z/S})|_{U_{\et}}\iso \Omega^\cdot_{(U,T^{ a}) /S}$ (see 1.8). Let $\Omega^{\natural\cdot }_{Z/S} (U)$ be the total complex of the simplicial dg algebra $\Omega^\cdot_{Z/S}(U,T^\cdot )$, $\Omega^{\natural n }_{Z/S} (U)=\oplus_i  \, \Omega^{i}_{Z/S}(U,T^{ i-n})$, filtered by the Hodge-pd filtration. We see that $\Omega^{\cdot }_{Z/S} (U,T^{ 0} )\hra \Omega^{\natural \cdot }_{Z/S} (U )$ is a filtered quasi-isomorphism, so  $Ru^{\log}_{Z/S *}(\CO_{Z/S})$ is represented by the filtered complex of presheaves  $U\mapsto   \Omega^{\natural\cdot}_{Z/S} (U )$. 

Now the de Rham pullback  $r^* : \Omega^\cdot_{P(U)_\cdot /(A,L)} \to    \Omega^{ \cdot}_{Z/S} (U,T^\cdot )  $ for the map of log schemes $r: T^{\cdot}\to \Spec\, P(U)_\cdot$ (see  (1.3.1)) sends the Hodge filtration to the Hodge-pd one, and (1.9.1) is the  map between the total complexes.

\proclaim{Theorem} Suppose $(Z,\CM)$, $(S,\CL )$ are fine and $f: (Z,\CM )\to (S,\CL )$ is  an integral locally log complete intersection map. Then (1.9.1) yields quasi-isomorphisms  $$ L\Omega^\cdot_{(Z,\CM ) /(S,\CL)} /F^m \iso Ru^{\log}_{Z/S *} ( \CO_{Z/S} /\CJ_{Z/S}^{[m]}). \tag 1.9.2$$  \endproclaim

\demo{Proof}  (a) Let $f: (Z,\CM )\to (S,\CL )$ be a map of fine log schemes. Recall (see \cite{KS} 4.4.2, 4.4.4 or \cite{Ol1} 6.8) that $f$ is {\it locally log complete intersection} map if \'etale locally it can be factored as $(Z,\CM )\buildrel{i}\over\to (Y,M_Y )\buildrel{g}\over\to (S,\CL )$ with $g$ log smooth, $i$ an exact closed embedding, and $Z\hra Y$ a regular immersion. Then, by loc.~cit., for any other factorization $f=g'i'$ with the first two properties, the third one holds automatically.
We say (see \cite{K1} 4.1, 4.3) that $f$ is {\it integral at a closed point} $z\in Z$ if the map $\Bbb Z [ (\CL/\CO^\times_S)_{f(z)}]\to \Bbb Z [(\CM /\CO^\times_Z )_z ]$ induced by $f$ is flat, and $f$ is {\it integral} if it is integral at every $z\in Z$. If $f$ is integral at $z$, then it is integral on an \'etale neighborhood of $z$.\footnote{Proof: We can assume that $f$ admits a fine chart $c_M : M\to \Gamma (Z,\CM )$, etc., as in 1.1, and for any $m\in M$ the function $\alpha c_M (m)\in \Gamma (Z,\CO_Z )$ is invertible if it is invertible at $z$.  Then $f$ is integral. This follows from \cite{K1} 4.1 and the next observation: Let $h: Q\to P$ be a map of integral monoids that satisfies condition (iv) of  \cite{K1} 4.1(i). Then for any submonoid $P' \subset P$ the map of the integral monoid quotients $Q/h^{-1}(P') \to P/P'$  satisfies the same condition (iv).} 
So for 
$f$ as in Theorem  the local factorization $f=gi$  can be chosen with $g$ integral (so $Y$ is flat over $S$ by \cite{K1} 4.5).

(b)  We return to the theorem. Its assertion is local, so, combining \cite{K1} 3.5, 4.1, we can assume that there is a fine chart $L\to M$ for $f$ (see 1.1) such that $\Bbb Z [ M]$ is $\Bbb Z [L]$-flat, the map $L^{\text{gr}}\to M^{\text{gr}}$ is injective with the cokernel having prime to $p$ torison, and
the map $i: Z\to Y:= \Spec\, (\Bbb Z [M]/p^n )\times_{\Spec \, (\Bbb Z [L]/p^n )} S$ is a regular immersion with parameters $t_1,..\,, t_m \in \CO_Y$. 
Let $\CM_Y$ be  the log structure on $Y$ defined by the chart $M$, and $\CM'$, $\CL'$ be the log structures on 
$\Spec\, (\Bbb Z [M]/p^n )$, $ \Spec \, (\Bbb Z [L]/p^n )$
 defined by  $M$, $L$. 
The  log cotangent complex L$_{(\Spec\, (\Bbb Z[M]/p^n ),\CM' )/(
 \Spec \, (\Bbb Z [L]/p^n ), \CL' )}$ equals L$_{( \Bbb Z [M]/p^n ,M)/(\Bbb Z [L]/p^n ,L)}\iso (\Bbb Z [M]/p^n )\otimes (M^{\text{gr}}/L^{\text{gr}})$. By base change (we use the flatness),  L$_{(Y,\CM_Y )/(S,\CL )}\iso  \Omega^1_{(Y,\CM_Y )/(S,\CL )}$ $ \iso \CO_Y \otimes (M^{\text{gr}}/L^{\text{gr}})$, where $d\log (c_M (m))$ is identified with  $\alpha c_M (m)\otimes m$ (see 1.1 for the notation). One has (see \cite{Ol1} 8.22) L$_{(Z,\CM )/(Y,\CM_Y )}=$ L$_{Z/Y}\iso \CO_Z^m [1]$, the generators $e_1,\ldots, e_m$ of $
\CO_Z^m$ correspond to regular parameters $t_1,..\,, t_m \in \CO_Y$. By the transitivity, one has L$_{(Z,\CM )/(S/\CL )}=\CC one(\delta : \CO_Z^m \to \Omega^1_{(Y,\CM)/(S,\CL )}|_Z )$, $\delta (e_i )= dt_i |_Z$, so $\gr^\cdot_F L\Omega^\cdot_{(Z,\CM )/(S,\CL )}=L\Lambda^\cdot \CC one(\delta )$. By \cite{Ill1} Ch.~I, 4.3.2.1(ii), this is $\CO_Z \langle e_1 ,\ldots ,e_m \rangle $ $\otimes_{\CO_Y}\Omega^\cdot_{(Y,\CM_Y )/(S,\CL )}$, where the first factor is the divided powers polynomial algebra.

 We compute $Ru^{\log}_{Z/S *} ( \CO_{Z/S})$ using (1.8.1) with $(Z,P)$ equal to the pd-envelope of $i: (Z,\CM)\hra (Y,\CM )$. Since $i$ is exact, $P$ is the pd-envelope of $Z\hra Y$,
so $\oplus\, \CJ^{[m]}_P /\CJ^{[m+1]}_P \iso \CO_Z \langle e_1 ,\ldots ,e_m \rangle$. So, by Remark (iii) in 1.8,  $\gr_F^\cdot Ru^{\log}_{Z/S *} ( \CO_{Z/S})\iso \CO_Z \langle e_1 ,\ldots ,e_m \rangle \otimes_{\CO_Y}\Omega^\cdot_{(Y,\CM )/(S,\CL )}$. The above two identifications provide an isomorphism of dg algebras $ \gr^\cdot_F L\Omega^\cdot_{(Z,\CM )/(S,\CL )} \iso \gr_F^\cdot Ru^{\log}_{Z/S *} ( \CO_{Z/S})$, which clearly coincides with the associated graded isomorphism to (1.9.1), q.e.d.
 \qed
\enddemo

\remark{Remarks} (i) The assertion remains true, by the direct limit argument, if $(Z,\CM )$ is projective limit of log $(S,\CL )$-schemes $(Z_\alpha ,\CM_\alpha )$ as in the theorem with respect to a directed family of affine transition maps. \newline
(ii) In the original theorem \cite{Ol1} 6.10, Olsson uses his version of log cotangent complex; there $f$ need not be integral.  For integral $f$ Olsson's version of log cotangent complex coincides with Gabber's one (that we use) by \cite{Ol1} 8.34.
\newline
(iii) By Bhatt \cite{Bh2} 7.22, (1.9.1) is itself a quasi-isomorphism  if, in addition, 
$f\otimes\Bbb Z/p$ is of Cartier type and
$Z$, $S$ are $\Bbb Z/p^n$-flat. We will not use this result. 
\endremark

\subhead{\rm 1.10}\endsubhead {\it The Cartier isomorphism.} Suppose our  $S$ is an $\Bbb F_p$-scheme. Let $(Y,\CN )$ be any log $S^\sharp$-scheme. Then for any  $(V,T)\in (Y/S)^{\log}_{\crys}$ the Frobenius map $Fr_T$ kills $\CJ_T$ since it is a pd-ideal, i.e., $Fr_T$ factors as $T\to V\hra T$; denote the first arrow by $Fr'_T$. The datum of all maps $Fr'_T$ forms an extension of the canonical morphism of topoi $u^{\log}_{Y/S}: (Y/S)^{\log}_{\crys}\!\tilde{}\to Y\tilde{_{\et}}$ (see 1.5) to a morphism $Fr'_Y$ of the {\it ringed} topoi.

Let $f : (Z,\CM )\to (Y,\CN)$ be a log smooth map of Cartier type between fine log $S^\sharp$-schemes; assume that $Z$, $Y$ are quasi-compact and quasi-separated. 

\proclaim{Theorem} The $\CO_{Y/S}$-complex $Rf_{\crys *} \CO_{Z/S}$ carries a natural finite filtration $\con_{\cdot}$, called the conjugate filtration, together with canonical Cartier quasi-isomorphisms $$ C=C^q_f :\, \gr_q^{\con} Rf_{\crys *} \CO_{Z/S}\iso Fr^{\prime *}_Y Rf_* \Omega^q_{(Z,\CM)/(Y,\CN )}[-q]. \tag 1.10.1$$\endproclaim

\remark{Remark}  The conditions imply that $\Omega^q_{(Z,\CM )/(Y,\CN )}$ are locally free $\CO_{Z}$-modules of finite rank, $f$ is flat,
and $Rf_* \Omega^q_{(Z,\CM)/(Y,\CN )}$ are quasi-coherent  of finite Tor-dimension.  \endremark

\demo{Proof} For  $(V,T)$ as above, one has $(Rf_{\crys *} \CO_{Z/S})(V,T)=R\Gamma ((Z_V /T)_{\crys}^{\log},\CO_{Z_V /T})$ $=
R\Gamma (Z_V, Ru^{\log}_{Z_V /T} \CO_{Z_V /T})$. The canonical filtration on $
Ru^{\log}_{Z_V /T} \CO_{Z_V /T}$ yields thus a  filtration $\con_{(V,T)\cdot}$ on $(Rf_{\crys *} \CO_{Z/S})(V,T)$. The filtrations $\con_{(V,T)\cdot}$ are compatible with morphisms of $(V,T)$. They form the promised conjugate filtration
 $\con_\cdot$. 
 
 Identifications (1.10.1)  have local nature: they come from isomorphisms $$ C :  H^q Ru^{\log}_{Z_V /T *} ( \CO_{Z_V/T})\iso
 \Omega^q_{(Z_V,\CM )/(V,\CN)}\otimes_{f_V^{-1} \CO_V} f_V^{-1} \CO_T ,
 \tag 1.10.2$$ where $\otimes$ in the r.h.s.~is taken for the map $Fr^{\prime *}_T : \CO_V \to \CO_T$, by applying $R\Gamma (Z_V, \cdot )$.
 
To define (1.10.2), let us find
  a  map   $C^{-1}: \Omega^q_{(Z_V,\CM )/(V,\CN)}\otimes_{f_V^{-1} \CO_V} f_V^{-1} \CO_T \to  H^q Ru^{\log}_{Z_V /T *} ( \CO_{Z_V/T})$  of graded $f_V^{-1}\CO_T$-algebras such that for $f\in \CO_{Z_V}$, $m\in \CM$ one has $C^{-1}(f)=  f^* (t^p )$, $C^{-1}(d\log m)= m^* (d\log t )$. Here 
 $f$, $m$ at the r.h.s.~are viewed as the $\Bbb A^1_{(T,\CN_T )}$- and  $\Bbb A^{(1)}_{(T,\CN_T )}$-valued maps, and $t^p$, $d\log t$ are the usual cohomology
 classes in $ H^0_{\!\dR} (\Bbb A^{1}_{(T,\CN_T  )}/(T,\CN_T  ) )    =H^0 R\Gamma ((\Bbb A^1_{( T,\CN_T)}/(T,\CN_T))_{\crys},\CO_{\Bbb A^1_{(T,\CN_T)}/T} ) $ and $  H^1_{\!\dR} (\Bbb A^{(1)}_{(T,\CN_T )}/ (T,\CN_T))$ $= H^1 R\Gamma ((\Bbb A^{(1)}_{(T,\CN_T )}/(T,\CN_T ))_{\crys},\CO_{\Bbb A^1_{( T,\CN_T)}/T} ) $ (see 1.1, (1.8.3)). The above properties determine $C^{-1}$ uniquely, so it is enough to find one such $C^{-1}$ locally.  By Remark (ii) in 1.1, we can assume that $(Z_V ,\CM )$ extends in a log smooth way over $(T,\CN_T )$. Now we are in the setting of  \cite{K1} 4.12(1), which provides  (via (1.8.3)) $C^{-1}$. By loc.~cit., $C^{-1}$
 is an isomorphism, and 
  $C$ of (1.10.2) is its inverse. \qed 
\enddemo
 
\subhead{\rm 1.11}\endsubhead {\it Perfect crystals and base change.} We are in general situation of 1.5, so $S^\sharp$ is as in 1.3. Let $(Y,\CN )$ be any integral quasi-coherent log $S^\sharp$-scheme. 

A bounded complex $\CF^\cdot$ of $\CO_{Y/S}$-modules on  $(Y/S)^{\log}_{\crys}$ is an
{\it  $\CO_{Y/S}$-crystal} (in derived sense)  if for every  $\phi   : (V',T' )\to (V,T)$ in 
$(Y/S)^{\log}_{\crys}$  the pullback  map $\phi^{-1}\CF^\cdot_{(V,T)}\to \CF^\cdot_{(V',T')}$ yields a quasi-isomorphism of $\CO_{T'}$-complexes $L\phi^{*}\CF^\cdot_{(V,T)}\iso \CF^\cdot_{(V',T')}$. We say that $\CF^\cdot$ is  {\it perfect} if for every $(V,T)\in (Y/S)^{\log}_{\crys}$ the complex $\CF^\cdot_{(V,T)}$ is $\CO_T$-perfect.
Perfect crystals form  a full triangulated subcategory
  $D^{\pcr}(Y/S)$
of the derived category  $D^b ((Y/S)^{\log}_{\crys},\CO_{Y/S})$
 of $\CO_{Y/S}$-modules.

Let now $\theta : (Y^{\nu},\CN^{\nu} )/S^{\nu\sharp} \to (Y,\CN )/S^\sharp$ be a map of  data as above.

 \remark{Exercise} If $\CF^\cdot$ is a (perfect) $\CO_{Y/S}$-crystal, then $L\theta^*_{\crys} (\CF^\cdot )$ is a (perfect) $\CO_{Y^{\nu}/S^{\nu}}$-crystal. For any $(V,T)\in (Y/S)^{\log}_{\crys}$, $(V^{\nu},T^{\nu})\in (Y^{\nu}/S^{\nu})^{\log}_{\crys}$, and a map of log pd-thickenings $\theta_{T^{\nu}/T} : (V^{\nu},T^{\nu})\to (V,T)$ compatible with $\theta$, one has a canonical identification $L\theta^*_{\crys} (\CF^\cdot )_{(V^{\nu},T^{\nu})}\iso L\theta^*_{T^{\nu}/T} (\CF^\cdot_{(V,T)})$.
 \endremark

\enspace

The next version of base change theorem  \cite{K1} 6.10  is sufficient for our purposes. 
Let $f: (Z,\CM )\to (Y,\CN )$ be a log smooth integral  map  
of fine log schemes over  $S^\sharp_1 := S^\sharp  \otimes\Bbb F_p$; assume that  $Z/Y$ is  quasi-compact and separated.
For  $\theta$ as above, let $f^{\nu} : (Z^{\nu},\CM^{\nu} )\to (Y^{\nu},\CN^{\nu} )$ be the $\theta$-pullback of $f$.
Then $Z$ is flat over $Y$ (by \cite{K1} 4.5), $Z^\nu $ is flat over $Y^\nu\!$, and
 $(Z^{\nu},\CM^{\nu} )$ is integral  quasi-coherent (by \cite{K1} 4.3.1 and Exercise (ii) in 1.1).

 \proclaim{Theorem} (i) The complex $Rf_{\crys *}(\CO_{Z/S})$ is an $\CO_{Y/S}$-crystal, and the pullback map yields a
  canonical identification $$L\theta^*_{\crys} Rf_{\crys *}( \CO_{Y/S})\iso Rf^{\nu}_{\crys *}( \CO_{Z^{\nu}/S^{\nu}}). \tag 1.11.1$$
 (ii) If $Z$ is proper over $Y$ and $f$ is of Cartier type, then $Rf_{\crys *}(\CO_{Z/S})$ is perfect.
 \endproclaim

\demo{Proof} (i) We can assume that $Y$ is affine. 
Computing $Rf_{\crys *}(\CO_{Z/S})$ using  a finite covering of $Z$ by affine opens, we get a finite filtration on it with $\gr_\cdot  Rf_{\crys *}(\CO_{Z/S})$ equal to a finite direct sum of complexes of type $R(f|_U)_{\crys *}(\CO_{U/S})$, $U\subset Z$ is an affine open. Thus we can assume that $Z$ is affine.
Using Remark (ii) in 1.1 and (1.8.3), we see that each   $Rf_{\crys *}(\CO_{Z/S})_{(V,T)}$ can be realized as a finite complex of flat $\CO_T$-modules, for every  $\phi: (V',T')\to (V,T)$ in  $(Y/S)^{\log}_{\crys}$ the pullback map yields a quasi-isomorphism $L\phi^* Rf_{\crys *}(\CO_{Z/S})_{(V,T)}\iso Rf_{\crys *}(\CO_{Z/S})_{(V',T')}$, and the pullback map for $\theta$ yields  quasi-isomorphism (1.11.1).
  
(ii) Let us check that the complexes $Rf_{\crys *}(\CO_{Z/S})_{(V,T)}$ are $\CO_T$-perfect. By (i), $Rf_{\crys *}(\CO_{Z/S})_{(V,T)}\otimes^L_{\CO_T} \CO_{T_1}= Rf_{\crys *}(\CO_{Z/S_1})_{(V,T_1 )}$, where $T_1 
 :=T\otimes\Bbb F_p$. Since $p $ is nilpotent in $\CO_T$, it suffices to check
 that $Rf_{\crys *}(\CO_{Z/S_1})_{(V,T_1 )}$ is $\CO_{T_1}$-perfect, which follows from (1.10.1). \qed 
\enddemo

\remark{Remarks} (a) If  $(V,T)\in (Y/S)^{\log}_{\crys}$ is such that $\CJ_T^{[n]}=0$, then $Rf_{\crys *}(\CO_{Z/S})_{(V,T)}\iso Rf_{\crys *}(\CO_{Z/S}/ \CJ_{Z/S}^{[m]})_{(V,T)}$ for $m\ge n+\dim Z/Y$.\footnote{Proof: As in the proof of (i), we can assume that $Z$, $Y$, $V$ are affine. Then $(Z,\CM )_V$ can be extended to a log smooth scheme over $(T,\CN_T )$ by \cite{K1} 3.14. Now use  Remark (iii) in 1.8.}
 \newline
(b) If we drop the Cartier type assumption in (ii),
 then $Rf_{\crys *}(\CO_{Z/S})_{(V,T)}$ is still $\CO_T$-perfect for those $(V,T)$ that the ideal $\CJ_T$ is nilpotent.\footnote{ Proof: It suffices, by  the nilpotency,  to check that $Rf_{\crys *}(\CO_{Z/S})_{(V,T)}\otimes^L \CO_V $ is $\CO_V$-perfect; this complex equals $Rf_{\crys *}(\CO_{Z/S})_{(V,V)}= (Rf_* \Omega^\cdot_{(Z,\CM )/(Y,\CN )})_V$ by (i) and 1.8, which is  perfect.} In particular, $Rf_{\crys *}(\CO_{Z/S})_{(Y,Y)}$ is $\CO_Y$-perfect. So $Rf_{\crys *}(\CO_{Z/S})$ is perfect
  if $Y=S$ (use (i)). 
\endremark

\subhead{\rm 1.12}\endsubhead {\it The p-adic setting; absolute crystalline cohomology.}  One can generalize slightly the setting for log crystalline cohomology taking for $S^\sharp =(S,\CL,\CI )$ a
 {\it formal p-adic log pd-scheme}, which is the same as sequence of exact closed embeddings of log pd-schemes $S^\sharp_1 \hra S^\sharp_2 \hra \ldots $ such that $\CO_{S_{n}}=\CO_{S_{n+1}}\otimes \Bbb Z /p^{n}$, $\CI_{{n}}=\CI_{n+1} \CO_{S_{n}}$. Assume that $(S,\CL )$ is quasi-coherent, i.e., all $(S_n ,\CL_n )$ are quasi-coherent. For an integral quasi-coherent  log $S^\sharp$-scheme $(Z,\CM )$, which is a  log $S^\sharp_n$-scheme for $n$ sufficiently large, one defines its log crystalline site $(Z/S)^{\log}_{\crys}$  as in 1.5. One has fully faithful embeddings $(Z/S_n)^{\log}_{\crys}\hra (Z/S_{n+1})^{\log}_{\crys}\hra \ldots$, and
 $(Z/S)^{\log}_{\crys}= \cup (Z/S_n )^{\log}_{\crys}$. The constructions and results of 1.5 remain true in the present setting. 
 
 For a sheaf $\CF$  on 
$(Z/S)^{\log}_{\crys}$  we denote by $\CF_{(n)}$ its restriction to $(Z/S_n)^{\log}_{\crys}$. Then $\CF_{(n)} = i_{n\,\crys}^* (\CF )$, where $i_n  : (Z,\CM ) /S^\sharp_n \to (Z,\CM )/S^\sharp$ is the evident map. The functor $ i_{n\,\crys}^*$ is exact, and it admits an evident exact left adjoint $ i_{n\,\crys  !}$, so $ i_{n\,\crys}^*$ sends injective sheaves to injective ones. Therefore the functors $Rf_{\crys *}$, $Lf^*_{\crys}$ commute with the passage $\CF \mapsto \CF_{(n)}$, and one has $Ru^{\log}_{Z/S *}(\CF )= \holim_n  Ru^{\log}_{Z/S_n *}(\CF_{(n)} )$,
$R\Gamma ((Z/S)^{\log}_{\crys},\CF )=\holim_n R\Gamma ((Z/S_n )^{\log}_{\crys},\CF_{(n)}  )$. 

\remark{Example} For a complex $\CF^\cdot$ of sheaves on $(S_1 ,S)^{\log}_{\crys}$, one has $R\Gamma ((S_1 ,S_n )^{\log}_{\crys} ,\CF^\cdot_{(n)} )$ $=\CF^\cdot (S_1 , S_n )$,  $R\Gamma ((S_1 ,S )^{\log}_{\crys} ,\CF^\cdot  )=\holim_n \,\CF^\cdot (S_1 , S_n )$. Unless $\CF^\cdot$ is a crystal, the base change map $R\Gamma ((S_1 ,S )^{\log}_{\crys} ,\CF^\cdot  )\otimes^L \Bbb Z/p^n \to R\Gamma ((S_1 ,S_n )^{\log}_{\crys} ,\CF^\cdot_{(n)} )$ need not be a quasi-isomorphism.
\endremark
\enspace

Consider now $S^\sharp_n =\Spec\, (\Bbb Z/p^n )$ equipped  with the trivial log structure and  $\CI_n =p(\Bbb Z/p^n )$ with the standard pd structure. Any 
 integral quasi-coherent log $\Bbb F_p$-scheme $(Z,\CM )$  is automatically a log $S^\sharp$-scheme. We refer to $(Z,\CM)_{\crys (n)}:=((Z,\CM )/S_n^\sharp )_{\crys}$, $ (Z,\CM)_{\crys}:= ((Z,\CM )/S^\sharp )_{\crys}$ as  the {\it  absolute crystalline sites}. One has the  {\it absolute log crystalline complexes} $R u^{\log}_{Z/\Bbb Z_p *} (\CO_{Z/\Bbb Z_p } )$,    
 $R\Gamma_{\!\crys}(Z,\CM )_{(n)}$ $:= R\Gamma ((Z,\CM)_{\crys (n)},\CO_{Z/(\Bbb Z/p^n ) })$, and
 $R\Gamma_{\!\crys}(Z,\CM ):= R\Gamma ((Z,\CM)_{\crys},\CO_{Z/\Bbb Z_p })=  R\Gamma (Z_{\et},R u^{\log}_{Z/\Bbb Z_p *} (\CO_{Z/\Bbb Z_p} ) ) =\holim_n \, R\Gamma_{\!\crys}(Z,\CM )_{(n )}$. We denote by
 $H^i_{\crys}(Z,\CM )$ the {\it absolute crystalline cohomology} $H^i R\Gamma_{\!\crys}(Z,\CM )$. 

 If $(Z,\CM )$ is an integral quasi-coherent log $\Bbb Z_p$-scheme, then for $n\ge 1$ we set $(Z,\CM )_n := (Z_n ,\CM_n )$, where $Z_n :=Z\otimes \Bbb Z/p^n $ and $\CM_n$ is the restriction of  $\CM$ to $Z_n \subset Z$. We set  $R\Gamma_{\!\crys}(Z,\CM )  := R\Gamma_{\!\crys}(Z_1 ,\CM_1 )$, etc.
 
\remark{Remarks} (i) If $Z$ is a flat $\Bbb Z_p$-scheme, then $(Z,\CM )_n$ is a log $S^\sharp$-scheme, and its log crystalline complex equals   
$R\Gamma_{\!\crys}(Z,\CM )$ by the  crystalline invariance property.
\newline (ii) For a perfect field $k$ of characteristic $p$, let $S^\sharp_{k\, n}$ be $\Spec\, W_n (k)$ equipped  with the trivial log structure and  $\CI=pW_n (k)$ with the standard pd structure. If  $(Z,\CM )$ is an  integral quasi-coherent log $k$-scheme $(Z,\CM )$, then every its pd-thickening is automatically a $W(k)$-scheme, i.e.,  
$((Z,\CM)/ S^\sharp_k )_{\crys}=(Z,\CM )_{\crys} $. Thus the absolute crystalline complexes of $(Z,\CM )$ are $W(k)$-algebras, etc. 
 \endremark

\subhead{\rm 1.13}\endsubhead {\it A digression on difference equations and Dwork's trick.} Let $R$ be a ring, $\varphi_R$ be an endomorphism of $R$. Denote by
 $R_{\varphi}$  the associative algebra generated by its subring $R$ and $\varphi$ with relation $\varphi r=\varphi_R (r)\varphi$, $r\in R$. Thus a (left) $R_{\varphi}$-module is the same as  a left $R$-module $F$ equipped with a $\varphi_R$-semilinear endomorphism $\varphi_F$ called {\it $\varphi_R$-action} or simply {\it $\varphi$-action}; we usually abbreviate $(F,\varphi_F )$ to $F$. 
 Denote by $R_{\varphi}$-mod the abelian category of $R_{\varphi}$-modules and by  $D_\varphi (R )$ its derived category. 
  
 For $F_1 ,F_2 \in R_\varphi$-mod  consider the map $\delta  : \Hom_R (F_1 ,F_2 )\to\Hom^{(\varphi )}_{R}(F_1 ,F_2 )$,  $\delta =\delta' -\delta''$, $\delta' (\xi ):= \xi \varphi_{F_1}$, $\delta''(\xi):= \varphi_{F_2}\xi $, where $\Hom_R^{(\varphi )}$ is the group of $\varphi$-semilinear maps, so
 $\Ker (\delta )=\Hom_{R_\varphi}(F_1 ,F_2 )$. Set $\Hom^\natural_{R_\varphi}(F_1 ,F_2 ):=\CC one (\delta )[-1]$.
 
 \remark{Exercises} (i)  The map $R\Hom_{R_\varphi}(F_1 ,F_2 )\to R\Hom^\natural_{R_\varphi}(F_1 ,F_2 )$ is a quasi-isomorphism, and $R\Hom^\natural_{R_\varphi}(F_1 ,F_2 )$ equals $\CC one (R\Hom_R (F_1 ,F_2 )\to R\Hom^{(\varphi )}_{R}(F_1 ,F_2 ))[-1]$ where $R\Hom_R$, $R\Hom^{(\varphi )}_{R}$ are computed in the derived category of $R$-modules.\footnote{Hint: $\delta$ is surjective when $F_1$ is a free $R_\varphi$-module, and $R_\varphi$ is free as a (left) $R$-module.} \newline (ii) One has $\Hom^\natural_{R_\varphi}(F_1 ,F_2 )=\Hom_{R_\varphi}(F_1^l ,F_2 )= \Hom_{R_\varphi}(F_1 ,F^r_2 )$ where $F_1^l$, $F_2^r$ are natural left and right two-term resolutions of $F_1$,  $F_2$. Explicitly, $F_1^{l0}=R_{\varphi}\otimes_R F_1=\oplus_{n\ge 0} \varphi^{n*}_R F_1$,  where $\varphi^{n*}_R F_1:=R\otimes_{\varphi^n}F_1$,  $F_{1}^{l\,-1}= (\varphi^*_R F_1)^{l 0}$ (so if $F_1$ is $R$-projective, then the resolution $F_1^l$ is $R_\varphi$-projective), and $F_2^{r0}=\Pi_{n\ge 0}\phi_{R*}^n F_2$, $F_2^{r\,1}= (\varphi_{R *}F_2 )^{r\, 0}$.
  \endremark
 
\enspace
An $R_{\varphi}$-complex is said to be {\it $R$-perfect} if it is perfect as a  complex of $R$-modules; such objects form a thick subcategory $D^{\text{perf}}_{\varphi}(R)$ of $D(R_\varphi )$. Let $R^{\text{prf}}_{\varphi}$-mod be the category of  $R_{\varphi}$-modules which are
 finitely generated and projective as $R$-modules.  An $R_{\varphi}$-module $M$ is said to be {\it nondegenerate} if the $R$-linear extension $\varphi^l_M:
 \varphi^*_R M \to M$ of $\varphi_M$ 
 is an isogeny, i.e., $\varphi^*_R M\otimes \Bbb Q \iso M\otimes \Bbb Q$; an $R_{\varphi}$-complex $C$ is {\it nondegenerate} if  the composition $L\varphi^*_R C \otimes\Bbb Q \to    \varphi^*_R C\otimes \Bbb Q \to C   \otimes \Bbb Q$ is a quasi-isomorphism. Denote by $?^{\text{nd}}$ the subcategory of ? formed by nondegenerate objects, e.g.~we have a triangulated category  $D^{\text{perf}}_{\varphi}(R)^{\text{nd}}$.

\proclaim{Lemma} (i)  $R^{\text{prf}}_{\varphi}$-mod generates
 $D^{\text{perf}}_{\varphi} (R)$ as a triangulated category.
 \newline (ii)
 If $R$ is a mixed characteristic dvr, then $R^{\text{prf}}_{\varphi}$-mod$\,{}^{\text{nd}}$ generates $D^{\text{perf}}_{\varphi} (R)^{\text{nd}}$. 
\endproclaim

\demo{Proof} (i)  Let $C$ be an $R$-perfect $R_{\varphi}$-complex. As an $R$-complex, it is quasi-isomorphic to a complex of finitely generated projective $R$-modules of finite amplitude $[a,b]$, $a\ge b$. We show that 
$C$ lies in  the triangulated subcategory $D^{\text{prf}}_{\varphi} (R)$ of 
 $D^{\text{perf}}_{\varphi} (R)$ generated by $R^{\text{prf}}_{\varphi}$-mod  using induction by $a-b$. 

  If $a=b$, then $C\iso (H^a C)[-a]$, $H^a C \in R^{\text{prf}}_{\varphi}$-mod, and we are done.
Suppose $n=a-b \ge 1$. Then $H^a C$ is a finitely generated $R$-module. Pick a surjection $P\twoheadrightarrow H^a C$ where $P$ is a finitely generated projective $R$-module. The $\varphi$-action on $H^a C$ can be lifted to $P$, i.e., there is a $\varphi$-action $\varphi_P$ on $P$ such that $P \twoheadrightarrow H^a C$ is a map of 
of $R_{\varphi}$-modules. Consider the projective $R_{\varphi}$-resolution $P^l$ of $ P=(P,\varphi_P)$. The above surjection lifts to a map of $R_{\varphi}$-complexes $g: P^l [-a]\to C$. Since $P^l $ is
quasi-isomorphic to $P$,   $C$ lies in $D^{\text{prf}}_{\varphi} (R)$ if (and only if) $\CC one (g)$ lies in $D^{\text{prf}}_{\varphi} (R)$, and the latter assertion is  true by the induction assumption, q.e.d. 
 
(ii) To make the proof of (i) work in the present situation, it suffices to check that for $C$ nondegenerate and any $P\twoheadrightarrow H^a C$ as in loc.~cit., one can choose  $\varphi_P$  so that $(P,\varphi_P )$ is nondegenerate.  Let $\varphi_P$ be any lifting of the $\varphi$-action on $H^a C$. Let $Q\subset P$ be the kernel of the map  $P\to  H^a C \otimes\Bbb Q$. Then $\varphi_P$ preserves $Q$, $P/Q $ is a projective $R$-module, and $(P/Q)\otimes\Bbb Q \iso H^a C \otimes\Bbb Q$, so $( P/Q ,\varphi_{P/Q})\in R^{\text{prf}}_{\varphi}$-mod$\,{}^{\text{nd}}$. We can modify  $\varphi_P$ by adding to it any map $\varphi^*_R P \to Q\subset P$ which is sufficiently small in $p$-adic topology (here $p$ is the residual characteristic of $R$). If this map is  sufficiently general, then the resulting $\varphi_Q$, hence $\varphi_P$, is an isogeny; we are done. 
 \qed
\enddemo

Suppose now  $R$ is a $p$-adically complete commutative algebra, $R\iso \limleft R/p^n R$, and  $I\subset R$ is a closed (for the $p$-adic topology) ideal preserved by $\varphi_R$  such that the induced endomorphism $\varphi_W$ of $W:=R/I$ is invertible and
$\varphi_R$  is  {\it topologically nilpotent} on $I$ (i.e., $\varphi_R$ acts  nilpotently on $I/pI$, hence on $I/p^n I$). Then $I$ consists of all $r\in R$ such that $\varphi_R^n (r)\to 0$. The projection $R\twoheadrightarrow W$ admits a unique section $ s: W\to R$ compatible with the action of $\varphi$ (which is automatically a ring homomorphism).\footnote{To see this, notice that for every (set-theoretic) section  $s' : W\to R$ the sequence of sections 
$\varphi_R^n s' \varphi_W^{-n}$ converges $p$-adically, and its limit $s$ does not depend on the choice of $s'$.} The corresponding
 base change functors 
$\pi^* : W_{\varphi}$-mod $\rightleftarrows R_{\varphi}$-mod $: i^*$, 
$ \pi^* M := R\otimes_{W}M$, $i^* N := W\otimes_R N =N/IN$, $i^* \pi^* M=M$,  preserve the subcategories $?^{\text{prf}}$-mod; the derived functors $L\pi^* ,L i^*$ 
preserve  the subcategories $D^{\text{perf}}_{\varphi}(?)$. They preserve the subcategories of non-degenerate objects as well.

We say that $\varphi$ is {\it strongly topologically nilpotent} on $I$ if for every $m>0$  one can find a  finite filtration of $I$ by closed $\varphi_R$-invariant ideals  such that $ \varphi_R (\gr^\cdot I )\subset p^m \gr^\cdot I $. 
 
 \proclaim{Proposition} If  $W$ is a mixed characteristic dvr and $\varphi$ is strongly topologically nilpotent on $I$, then the functors $L\pi^* ,Li^*$ yield mutually inverse equivalences $$D^{\text{perf}}_{\varphi}(W)^{\text{nd}}\otimes\Bbb Q\buildrel{\sim}\over\leftrightarrow D^{\text{perf}}_{\varphi}(R)^{\text{nd}}\otimes\Bbb Q . \tag 1.13.1$$
 \endproclaim
 
\demo{Proof} (a) Let us show that {\it for any $P\in W^{\text{prf}}$-mod$\,{}^{\text{nd}}$ and $Q\in R^{\text{prf}}$-mod, one has $R\Hom_{W_{\varphi}} (P,IQ )\otimes\Bbb Q =0$.}

For any $M\in W_{\varphi}$-mod consider the two-term complex $\CC^\cdot (M):=\Hom^\natural_{W_{\varphi}} (P,M)$, $\CC^0 (M)=\Hom_W (P,M)$, $\CC^1 (M)=\Hom_W (\varphi^*_W P,M)$. By Exercises,  $R\Hom_{W_{\varphi}} (P,M)$ $=\CC^\cdot (M)$. 
Let $\CC_{\star}^\cdot (M)$ be the complex with the same components as $\CC^\cdot (M)$ and the differential $\delta'(\xi )= \xi \varphi^l_P$. Since $P$ is nondegenerate,  $\CC_{\star}^\cdot (M)\otimes\Bbb Q$ is acyclic. 

 The functor $\CC^\cdot $ is exact. We prove that $\CC^\cdot (IQ)\otimes\Bbb Q$ is acyclic by defining a finite filtration $IQ^{\cdot}$ on $IQ$ such that for $F:=\gr^\cdot IQ$ one has $ \CC^\cdot (F )\simeq  \CC_{\star}^\cdot (F )$.

For $m$ sufficiently large, there is $\psi : P\to\varphi^*_W  P$ with $\varphi^l_P \psi =p^{m-1} \id_P$, $\psi \varphi^l_P =p^{m-1} \id_{\varphi^*_W P}$. By the condition of the proposition, one can find a finite filtration $I^{(\cdot )}$ on $I$ by closed $\varphi$-invariant ideals such that $\varphi_R (\gr^{(\cdot )}I)\subset p^m \gr^{(\cdot )}I$. Set  $IQ^{\cdot}:= I^{(\cdot )}Q$. 
 
 One has $\varphi_F (F)\subset p^m F$. Therefore $\delta''$ on 
 $ \CC^\cdot (F )$ is divisible by $p^m$. Set 
$\chi :=\psi^{\tau }(p^{-m}\delta'') \in \End ( \CC^0 ( F ))$; then $\delta'' = p\delta'\chi$, i.e., $\delta(1-p\chi)=\delta'$. Since $\CC^0 (F)$ is p-adically complete,  $1 - p\chi$ is invertible, so it yields $\CC_{\star}^\cdot (F)\iso \CC^\cdot (F)$, q.e.d. 

(b)   The lemma and (a) imply that for every $P\in D^{\text{perf}}_{\varphi}(W)^{\text{nd}}$ and  $Q\in D^{\text{perf}}_{\varphi}(R)$ one has $ R\Hom_{R_{\varphi}}(L\pi^* P, Q)\otimes\Bbb Q \iso  R\Hom_{W_{\varphi}}(P,Li^* Q)\otimes\Bbb Q$. Thus the two functors $L\pi^* : D^{\text{perf}}_{\varphi}(W)^{\text{nd}}\otimes\Bbb Q\leftrightarrows D^{\text{perf}}_{\varphi}(R)^{\text{nd}}\otimes\Bbb Q : Li^*$ are adjoint.
  
Since $Li^*$ is left inverse to $ L \pi^*$, it  remains to show that for $Q\in D^{\text{perf}}_{\varphi}(R)^{\text{nd}}$ the adjunction $L\pi^* Li^* (Q\otimes\Bbb Q )\to Q\otimes\Bbb Q$ is a quasi-isomorphism. We can replace $Q\otimes\Bbb Q$ by $\CC one (L\pi^* Li^* (Q\otimes\Bbb Q) \to Q\otimes\Bbb Q )$, so it suffices to check that $Li^* (Q\otimes\Bbb Q )=0$ implies $Q\otimes\Bbb Q =0$. Since $Q$ is nondegenerate, the support $S$ of $Q\otimes\Bbb Q$, which is a closed subset of Spec$\, (R\otimes\Bbb Q )$, is $\varphi_R$-invariant. If $S$ is non-empty, then this implies that $\Spec\, (W\otimes\Bbb Q) $ lies in $S$.\footnote{Indeed, otherwise there is $r\in R\smallsetminus I$  that  vanishes on $S$. Suppose $r\,$mod$I \in p^n W^\times$. Since
$\varphi_R$ is topologically nilpotent on $I$ and invertible on $W$, for $m\gg 0$ one has $\varphi^{m*}_R (r)= p^n (a_m + pb_m )$, where $a_m \in s(W^\times )$ and $b_m \in I$. Hence  $\varphi^{m*}_R (r)\in p^n R^\times$   since $R$ is $p$-adically complete. Since $S$ is $\varphi_R$-invariant, $\varphi^{m*}_R (r)$ vanishes on $S$, i.e., $S=\emptyset$.}  Since
$Q\otimes\Bbb Q$ is perfect, this contradicts the assumption $Li^* (Q\otimes\Bbb Q) =0$, q.e.d.
\qed\enddemo

\subhead{\rm 1.14}\endsubhead {\it Frobenius crystals.} Suppose $S^\sharp$ as in the beginning of 1.12 is equipped with an endomorphism $\varphi_{S^\sharp}$ whose restriction to $(S_1 ,\CL_1 )$ is the Frobenius map $Fr_{(S_1 ,\CL_1 )}$. Then for any integral quasi-coherent $(Y ,\CN )$ over $S^\sharp_1$ the endomorphism
$\varphi = (Fr_{(Y,\CN )} ,\varphi_{S^\sharp} )$ of $(Y ,\CN )/S^\sharp$ acts on $(Y/S)^{\log}_{\crys}$. A {\it  Frobenius $\CO_{Y/S}$-module} is a  pair $(\CF ,\varphi_{\CF})$, where $\CF$ is an $\CO_{Y/S}$-module on $(Y/S)^{\log}_{\crys}$, $\varphi_{\CF}$ is a $\varphi_{\crys}$-action on $\CF$, i.e., a map $\varphi_{\CF}:\CF \to \varphi_{\crys *} (\CF )$; if $\CF$ is an 
$\CO_{Y/S}$-crystal, we call $(\CF ,\varphi_{\CF})$  an {\it $F$-crystal}. Frobenius $\CO_{Y/S}$-modules form an abelian category; 
let $D_{\varphi} ((Y/S)^{\log}_{\crys},\CO_{Y/S})$ be its derived category, and  $D^{\pcr}_{\varphi}(Y /S )= D^{\pcr}_{\varphi}((Y,\CN) /S^\sharp ) $ be the thick subcategory of {\it perfect F-crystals}, i.e., those $(\CF ,\varphi_{\CF})$ that $\CF$ is a perfect crystal (see 1.11). Such an  $(\CF ,\varphi_{\CF})$   is {\it nondegenerate} if the map $L\varphi_{\crys}^* (\CF )\to \CF$ that comes from $\varphi_{\CF}$ is an isogeny, i.e., it is a quasi-isomorphism in $D^{\pcr}(Y /S )\otimes\Bbb Q$; the corresponding category is denoted by $D^{\pcr}_{\varphi}(Y /S )^{\nd}.$ 

A morphism $\theta : (Y^\nu ,\CN^\nu)/S^{\nu\sharp} \to (Y,\CN )/S^\sharp$ compatible with $\varphi_{S^{\nu\sharp}}$, $\varphi_{S^\sharp}$ yields the pullback functor $\theta^*_{\crys}$ between the categories of Frobenius $\CO$-modules. The derived functor $L\theta^*_{\crys}$ preserves the subcategories of (nondegenerate) perfect $F$-crystals, and it is compatible with $L\theta^*_{\crys}$ from 1.5 via the forgetful functor
$(\CF,\varphi_{\CF})\mapsto \CF$.

\remark{Remarks} (i) The endofunctor $L\varphi^*_{\crys}$ of $ D^{\pcr}_{\varphi}(Y /S )^{\nd}\otimes \Bbb Q$ is canonically isomorphic to the identity functor. \newline
(ii) For a perfect F-crystal $(\CF,\varphi_{\CF})$, its nondegeneracy is a $Y_{\et}$-local property. Suppose $Y$ is affine and $P=(P,\CN_P )$ is its pd-$S^\sharp$-smooth formal thickening (i.e.,  $P_n $ are pd-$S_n^\sharp$-smooth thickenings of $Y$) equipped with a lifting $\varphi_P$ of $\varphi_{S^\sharp}$. Write $P=\text{Spf}\, R$, $R$ is a $p$-adically complete ring, so $P_n =\Spec\, R_n$ and we have endomorphism $\varphi_R =\varphi_P^*$ of $R$. Then $\CF (P):=\holim_n \, \CF (P_n )$ is a perfect $R$-complex equipped with a $\varphi_R$-action $\varphi_{\CF(P)}$.
Let $\varphi^l_{\CF(P)}: R\otimes^L_{R,\varphi_R}\CF (P)= \holim_n \, R_n \otimes^L_{R_n ,\varphi_{R_n}}\CF (P_n ) \to \CF (P)$ be its $R$-linear extension. Then
 $\CF$ is nondegenerate if and only if $\varphi^l_{\CF(P)}\otimes\Bbb Q$ is a quasi-isomorphism.
\endremark
\enspace
We usually abbreviate $(\CF,\varphi_{\CF})$ to $\CF$, and  denote by
$\Hom_{\varphi}(\CF_1 ,\CF_2  )$ the group  of Frobenius $\CO_{Y/S}$-module morphisms. One has an exact sequence $0\to \Hom_{\varphi}(\CF_1 ,\CF_2  )$ $\to \Hom (\CF_1 ,\CF_2  )\to \Hom (\CF_1 ,\varphi_{\crys *} (\CF_2  ))$, the last arrow is $\xi\mapsto \varphi_{\CF_2}\xi - \varphi_{\crys *}(\xi )\varphi_{\CF_1}$.

\proclaim{Lemma} For any $\CF_1 ,\CF_2 \in D^+_{\varphi} ((Y/S)^{\log}_{\crys},\CO_{Y/S})$ one has a canonical quasi-isomor- phism $R\Hom_{\varphi}(\CF_1 ,\CF_2 )\iso \CC one\, (R\Hom (\CF_1 ,\CF_2 )\to R\Hom (\CF_1 ,R\varphi_{\crys *} (\CF_2  )))[-1]=  \CC one\, (R\Hom (\CF_1 ,\CF_2 )\to R\Hom (L\varphi^*_{\crys }(\CF_1) ,\CF_2  ))[-1]  $.
\endproclaim

\demo{Proof} The above exact sequence yields a map of complexes $\alpha : R\Hom_{\varphi}(\CF_1 ,\CF_2  )\to \CC one\, (R\Hom (\CF_1 ,\CF_2  )\to R\Hom (\CF_1 ,R\varphi_{\crys *} (\CF_2  )))[-1]$. Let us check that $\alpha$ is a quasi-isomorphism. 

The forgetful functor $(\CF ,\varphi_{\CF})\mapsto \CF$ from Frobenius $\CO_{Y/S}$-modules to $\CO_{Y/S}$-modules, admits a right adjoint $\CF \mapsto \CF_{(\varphi )}$. Explicitly, $ \CF_{(\varphi )}=\Pi_{n\ge 0}\varphi_{\crys *}^n (\CF )$ and $\varphi_{\CF_{(\varphi )}}$ is the  projection  $\Pi_{n\ge 0}\varphi_{\crys *}^n (\CF )\twoheadrightarrow \Pi_{n\ge 1}\varphi_{\crys *}^n (\CF )$. The functor $\CF \mapsto \CF_{(\varphi )}$ is left exact and sends injective objects to injective ones. Thus every Frobenius $\CO_{Y/S}$-module admits an embedding into  $ \CG_{(\varphi )}$ where $\CG$ is some injective $\CO_{Y/S}$-module. Therefore it suffices to check that $\alpha$  is a quasi-isomorphism assuming that $\CF_1 =\CF $ is any Frobenius $\CO_{Y/S}$-module   and $\CF_2 =  \CG_{(\varphi )}$ with $\CG$ injective. Then $R\Hom_{\varphi}(\CF,\CG_{(\varphi )})\buildrel{\sim}\over\leftarrow\Hom(\CF,\CG )\iso 
\CC one\, (\Hom (\CF ,\CG_{(\varphi )}  )\to \Hom (\CF ,\varphi_{\crys *} (\CG_{(\varphi )} ))[-1] $ $\iso
\CC one\, (R\Hom (\CF ,\CG_{(\varphi )}  )\to R\Hom (\CF ,R\varphi_{\crys *} (\CG_{(\varphi )} )))[-1]$, q.e.d. \qed
\enddemo

Suppose now our $(Y, \CN )$ is a fine log scheme, $Y$ is affine, and there is $P$ as in Remark (ii) such that $R$ has no $p$-torsion.
Let $f: (Z ,\CM )\to (Y ,\CN )$ be a log smooth map of Cartier type with
$(Z ,\CM )$  fine and $Z$ is proper over $Y$. Consider $\CF := Rf_{\crys *}(\CO_{Z /S})$.
 By the theorem in 1.11, $\CF  \in D^{\pcr}_{\varphi}(Y /S )$. The next  result  is a log version \cite{HK} 2.24 of a theorem of Berthelot-Ogus   \cite{BO2} 1.3:  

\proclaim{Theorem} The perfect F-crystal $ \CF$
  is nondegenerate. 
\endproclaim

\demo{Proof} We use notation from Remark (ii). Set $C:= \CF (P)=R\Gamma ((Z/P )^{\log}_{\crys}, \CO_{Z /P})$, so $C_n :=C\otimes\Bbb Z/p^n = 
\CF (P_n )=R\Gamma ((Z/P_n)^{\log}_{\crys}, \CO_{Z /P_n})$, $C=\holim_n \, C_n$. 
Then $C$ carries the Frobenius endomorphism $\varphi_{C}=\varphi_{\CF (P)}$, and we want to prove that $ \varphi^l_{C} \otimes\Bbb Q$ is a quasi-isomorphism.

Set $\CG_{n} := Ru^{\log}_{Z /P_n *} (\CO_{Z /P_n})$,
$\CG := Ru^{\log}_{Z /P *} (\CO_{Z /P})=\holim_n \, \CG_{n} $. These are $R_n$- and $R$-complexes of sheaves on $Z_{\et}$, and (i) of the theorem in 1.11 implies that $\CG_{n} = \CG\otimes^L \Bbb Z/p^n $. They carry natural Frobenius $\varphi_R$-actions $\varphi_{\CG_n}$, $\varphi_{\CG}$. 
One has $C_n =R\Gamma (Z_{\et},\CG_n )$,  $C= R\Gamma (Z_{\et},\CG )$, and  $\varphi_{C}$ comes from  $\varphi_{\CG}$. Let $\varphi^l_{\CG_n}: R_n \otimes^L_{R_n ,\varphi_{R_n}} \CG_n \to  \CG_n$ 
be the $R_n$-linear extensions of $\varphi_{\CG_n}$; set
 $\varphi_{\CG}^l :=\holim_n \, \varphi^l_{\CG_{n}} :  R\hotimes^L_{R,\varphi_R}\CG:=  \holim_n \, (R_n \otimes^L_{R_n ,\varphi_{R_n}} \CG_{n} )  \to  \CG$. Since  
 $R\Gamma (Z_{\et},R\hotimes^L_{R,\varphi_R}\CG)=
  \holim_n (R_n \otimes^L_{R_n ,\varphi_{R_n}}\! C_n ) $ and  $R\Gamma (Z_{\et}, \cdot \otimes\Bbb Q)= R\Gamma (Z_{\et}, \cdot )\otimes\Bbb Q$  (for
 $Z$ is quasi-compact quasi-separated), it suffices to show that
$\varphi_{\CG}^l \otimes\Bbb Q :  (R\hotimes^L_{R,\varphi_R}\CG )\otimes\Bbb Q  \to  \CG\otimes\Bbb Q$ is a quasi-isomorphism.

Now our assertion is $Z$-local, so to check it we can assume that $Z$ is affine. Choose a log smooth (formal) lifting $(T ,\CM_T  )/(P,\CN_P )$ of $(Z ,\CM )/(Y,\CN )$ together with a lifting $\varphi_{T}$ of the Frobenius compatible with $\varphi_P$. By (1.8.1), one has $(\CG ,\varphi_{\CG} )= (\Omega^\cdot  ,\varphi^*_T )$ where $\Omega^\cdot =\limleft \Omega^\cdot_n$, $\Omega^\cdot_n :=  \Omega^\cdot_{(T ,\CM_T  )_n /(P,\CN_P )_n }$. Consider the $p$-adic filtration on $\Omega^\cdot$; let $\Phi^*$ be its shift, i.e., $\Phi^m\Omega^\cdot $ is the maximal subcomplex of $\Omega^\cdot$ such that $\Phi^m \Omega^n \subset p^m \Omega^n$. Since $\varphi_T^* (\Omega^i )\subset p^i \Omega^i$, one has $\varphi_T^* (\Omega^\cdot )\subset \Phi^0\Omega^\cdot $.
Since $(\Omega^\cdot / \Phi^0\Omega^\cdot )\otimes\Bbb Q =0$ and $\Omega^i$ are $R$-flat, to finish the proof it suffices to check that the $R$-linear map $\varphi_T^{*l} : R\hotimes_{R,\varphi_R}\, \Omega^\cdot \to \Phi^0 \Omega^\cdot$ is a quasi-isomorphism.

Our complexes of sheaves are $p$-adically complete, have no $p$-torsion, and $\Phi^m \Omega^\cdot = p^m \Phi^0 \Omega^\cdot$ for $m\ge 0$. Thus it suffices to show that $ \varphi_{T_1}^{*l} :  R_1 \otimes_{R_1 ,\varphi_{R_1}} \!\Omega^\cdot_1 \to \gr^0_{\Phi}\Omega^\cdot$ is a quasi-isomorphism. Let $
\bar{\gr}^0_{\Phi}\Omega^\cdot$ be the quotient of $ \gr^0_{\Phi}\Omega^\cdot$
modulo the subcomplex generated by the images of $p^{i+1}\Omega^{i}\subset \Phi^0\Omega^i$. This subcomplex is acyclic, and the multiplication by $p^{-i}$
yields an isomorphism $ \bar{\gr}^0_{\Phi}\Omega^i \iso H^i \Omega^\cdot_1$.
The composition $R_1 \otimes_{R_1 ,\varphi_{R_1}} \!\Omega^\cdot_1 \to \gr^0_{\Phi}\Omega^\cdot \to \bar{\gr}^0_{\Phi}\Omega^\cdot$ is an isomorphism of complexes: indeed, its components coincide, via the previous identification, 
 with the Cartier isomorphism $C^{-1}$ from (1.10.2). We are done.
\qed \enddemo

\remark{Remark}  In fact, according to \cite{BO1} 8.20,
$\varphi_T^{*l} : R\hotimes_{R,\varphi_R} \Omega^\cdot \iso \Phi^0\Omega^\cdot$ comes from a natural global quasi-isomorphism  $R\hotimes_{R,\varphi_R}\CG \iso \Phi^0 \CG$. \endremark

\remark{Exercise}  Using $r$-iterated Cartier isomorphism, identify
 the $r$-th differential of the spectral sequence for the $p$-adic filtration on $\Omega^\cdot$ with the de Rham differential. 
\endremark

\subhead{\rm 1.15}\endsubhead {\it   $(\varphi ,N)$-modules and the Fontaine-Hyodo-Kato torsor.} Let $k$ be a perfect field of characteristic $p$,  $K_0 :=\text{Frac}\, W(k)$, 
$\varphi$ be the Frobenius automorphism of $W(k)$ and $K_0$. As in \cite{F2} 4.2, a {\it $\varphi$-module} over $K_0$ is a pair $(V,\varphi )$, where $V$ is
a finite-dimensional $K_0$-vector space, $\varphi =\varphi_V$ is a $\varphi$-semilinear automorphism of $V$; a
{\it $(\varphi ,N)$-module} is a triple $(V,\varphi ,N)$, where $(V,\varphi )$ is
a  $\varphi$-module   and $N=N_V$ is a $K_0$-linear endomorphism of $V$ such that $N\varphi  =p\varphi N$ (then $N$ is automatically nilpotent). One says that $V$ is {\it effective} (or of nonnegative slope) if it contains a $\varphi_V$-invariant $W(k)$-lattice, hence one preserved by both $\varphi_V$ and $N_V$. The category $(\varphi ,N)$-mod of $(\varphi ,N)$-modules is naturally a Tannakian tensor $\Bbb Q_p$-category, and $(V,\varphi ,N)\mapsto V$ is a fiber functor over $K_0$.\footnote{The same is true for the category of $\varphi$-modules.} Let $(\varphi ,N)^{\eff}$-mod be its abelian tensor subcategory of effective modules, and $D_{\varphi ,N}(K_0 )^{\eff} \subset D_{\varphi ,N}(K_0 )$ be the corresponding bounded derived categories. 

We usually abbreviate $(V,\varphi ,N)$ to $V$. Let $\varphi^* V=K_0 \otimes_{\varphi}V$ be a copy of $V$ equipped with the $\varphi$-twisted $K_0$-action,  $\varphi_{\varphi^* V}:=\varphi^* (\varphi_V )$, and $N_{\varphi^* V}:=p\varphi^*(N_V )$. Notice that $\varphi_V$ viewed as a map $\varphi^* V \to V$ is a morphism of  $(\varphi ,N)$-modules, which we denote again as $\varphi_V$.

For $(\varphi ,N)$-modules
$V_1 , V_2$, we denote by
 $\Hom_{\varphi,N}(V_1 ,V_2 )$  the group of $(\varphi ,N)$-module morphisms, and by $\Hom (V_1 ,V_2 )$  that of $K_0$-linear maps. Let $\Hom^\natural_{\varphi,N}(V_1 ,V_2 )$ be the complex $\Hom (V_1 ,V_2 )\to\Hom (\varphi^* V_1 ,V_2 )\oplus \Hom (V_1 ,V_2 )  \to \Hom (\varphi^* V_1 ,V_2 )$ supported in degrees $[0,2]$ with the differential $d^0 (\xi )= (\varphi_2 \xi - \xi\varphi_1, N_2 \xi - \xi N_1 )$, $d^1 (\chi, \psi )=N_2 \chi - p\chi N_1 - p\varphi_2 \psi + \psi\varphi_1$. Thus $\Hom_{\varphi,N}(V_1 ,V_2 )=H^0 \Hom^\natural_{\varphi,N}(V_1 ,V_2 )$.

\remark{Exercise} Show that $R\Hom_{\varphi,N}(V_1 ,V_2 )\iso \Hom^\natural_{\varphi,N}(V_1 ,V_2 )$. \endremark
\remark{Remark} Complexes $\Hom^\natural_{\varphi,N}$ compose naturally, so
they provide a dg category structure on $(\varphi ,N)$-modules. By Exercise, its homotopy category equals $D_{\varphi ,N}(K_0 )$.
\endremark

\enspace

Let $ \CL^0 $ be {\it the canonical} log structure on $S^0 := \Spec\, W(k)$: it is generated by a prelog structure 
$W(k)\smallsetminus \{ 0\} \to W(k)$ or, equivalently, by  one $\Bbb Z_{\ge 0} \to W(k)$, $m\mapsto p^m$. As in 1.12, we have log schemes $(S^0 ,\CL^0 )_n$.
Let now  $ (Y, \CN )$ be an integral quasi-coherent log scheme; let
$\pi: (Y, \CN )\to  (S^0 ,\CL^0  )_1 = (\Spec\, k ,\CL^0_{ 1} )$ be a morphism of log schemes. 
Notice that for a fixed $k$-scheme structure on $Y$ such $\pi$ amounts to a section  $l= \pi^* (\bar{p})\in \Gamma (Y,\CN )$ such that $\alpha (l)\in\Gamma (Y,\CO_Y )$ equals 0; here $\bar{p}\in \CL^0_{ 1}$ is the image of $p\in \CL^0$. Sometimes we write $l= l_{\pi}$, $\pi =\pi_l$.

 Consider  the absolute crystalline topology $(Y,\CN )_{\crys}$ (see 1.12) and the  category of nondegenerate perfect $F$-crystals
$D^{\pcr}_{\varphi} (Y )^{\nd}$.

\proclaim{Theorem-construction} (i) There is a natural functor $$\epsilon_{\pi}=\epsilon_{l} : D_{\varphi ,N}(K_0 )^{\eff}\to D^{\pcr}_{\varphi} (Y )^{\nd}\otimes\Bbb Q.  \tag 1.15.1$$ (ii) $\epsilon_\pi$ is compatible with base change, i.e.,  for any  $\theta: (Y', \CN' )\to (Y, \CN )$  one has a canonical identification $\epsilon_{\pi\theta} \iso L\theta^*_{\crys}\epsilon_{\pi}$.  For any $a\in k^{\times}$, $m\in\Bbb Z_{>0}$ there is a canonical identification $\epsilon_{al^m}(V,\varphi ,N)\iso \epsilon_l (V,\varphi, mN)$.
 \newline (iii) Suppose that $Y$ is a local scheme  with residue field $k$ and nilpotent maximal ideal,  $\CN /\CO^\times_Y =\Bbb Z_{\ge 0}$, and  map $\pi^* :\CL^0_{ 1}/k^\times \to \CN /\CO^\times_Y $ is injective. Then (1.15.1) is an equivalence of categories. \endproclaim

So, by (iii), we have the basic equivalence of categories $\epsilon := \epsilon_{\bar{p}}: D_{\varphi ,N}(K_0 )^{\eff}\iso D^{\pcr}_{\varphi} (  (S^0 ,\CL^0 )_1  )^{\nd}\otimes\Bbb Q$ and, by (ii), a canonical identification $\epsilon_\pi = L\pi^*_{\crys} \epsilon$. 

\demo{Proof} (i) For $\nu \in \Gamma (Y,\CN )$ let $\lambda_{\nu}$ 
be the preimage of $\nu$ by the 
the surjection of the monoid sheaves $\CN_{Y/\Bbb Z_p} \to \iota_* \CN$ on  $(Y,\CN )_{\crys}$ (see 1.5). This is a $(1+\CJ_{Y/\Bbb Z_p})^\times$-torsor by Exercise (iii) in 1.1. For any  $m\in\Bbb Z_{\ge 0}$ the $m$ power endomorphism of $\CN$ yields an identification $\lambda_{\nu}^m \iso \lambda_{\nu^m}$ of $(1+\CJ_{Y/\Bbb Z_p})^\times$-torsors.

Consider the $(1+\CJ_{Y/\Bbb Z_p})^\times$-torsor 
$\lambda_{\pi}:=\lambda_l$. Since $\phi^* (l)=l^p$, the action of $\varphi_{\crys}$ on $\CN_{Y/\Bbb Z_p}$ (see 1.5) yields a canonical identification of 
$(1+\CJ_{Y/\Bbb Z_p})^\times$-torsors  $$\varphi_{\crys}^* (\lambda_l )\iso \lambda_l^p  . \tag 1.15.2$$

Let $V=(V,\varphi ,N)$ be a $(\varphi ,N)$-module. Pick a $W(k)$-lattice $V^{(0)}\subset V$ preserved by $\varphi$ and $N$. Consider the ``constant" crystal $V^{(0) }_{Y/\Bbb Z_p}:=V^{(0)}\otimes_{W(k)} \CO_{Y/\Bbb Z_p}$.
The  group sheaf $(1+\CJ_{Y/\Bbb Z_p})^\times$ acts naturally on $V^{(0)}_{Y/\Bbb Z_p}$ by $\CO_{Y/\Bbb Z_p}$-linear automorphisms:
namely, $f\in (1+\CJ_{Y/\Bbb Z_p})^\times$ acts as $f^N =\exp (N\log f):=\Sigma N^i (\log f)^{[i]} $; the sum is finite since $N$ is nilpotent. Let $\epsilon_{\pi} (V^{(0)} )$ be the $\lambda_l$-twist of 
$V^{(0)}_{Y/\Bbb Z_p}$, i.e., the sheaf  of $(1+\CJ_{Y/\Bbb Z_p})^\times$-invariant maps 
$\lambda_l \to V^{(0)}_{Y/\Bbb Z_p}$. 

This is naturally a $F$-crystal: Indeed, by (1.15.2) one has a canonical identification $\varphi^*_{\crys}\epsilon_{\pi} (V^{(0)} ) \iso \epsilon_{\pi} (\varphi^* V^{(0)} )$, where $
\varphi^* V^{(0)} $ is a copy of $V^{(0)}$ viewed as a $W(k)$-lattice in the $(\varphi ,N)$-module $\varphi^* V$,\footnote{Indeed, $\varphi^*_{\crys}(V^{(0)}_{Y/\Bbb Z_p})= (\varphi^* V^{(0)})_{Y/\Bbb Z_p}$, so, by (1.15.2), $\varphi^*_{\crys}\epsilon_{\pi} (V^{(0)} )$ is the 
$\lambda^p_l$-twist of $(\varphi^* V^{(0)})_{Y/\Bbb Z_p}$ for  $(1+\CJ_{Y/\Bbb Z_p})^\times$-action $f\mapsto f^{\varphi^* N_V}$, which is the same as $\lambda_l$-twist for the action $f\mapsto f^{p\varphi^* N_V}=f^{N_{\varphi^* V}}$, q.e.d. } so  the morphism $\varphi_V : \varphi^* V^{(0)} \to V^{(0)}$ yields the Frobenius action $\varphi^*_{\crys}\epsilon_{\pi} (V^{(0)} )\to \epsilon_{\pi} (V^{(0)} )$.

Our $F$-crystal is clearly nondegenerate. 
Set $\epsilon_{\pi} (V):= \epsilon_{\pi} (V^{(0)} )\otimes\Bbb Q$; the construction does not depend on the choice of $ V^{(0)} $.  
 Extending the construction to complexes in the evident way, we get (1.15.1).

\enspace

(ii) The identification $\epsilon_{\pi\theta}(V^{(0)}) = \theta_{\crys}^* \epsilon_{\pi\theta}(V^{(0)})=L\theta_{\crys}^* \epsilon_{\pi\theta}(V^{(0)})$ comes from the 
isomorphism of $(1+\CJ_{Y'/\Bbb Z_p})^\times$-torsors
$\theta_{\crys}^* \lambda_{\pi} \iso \lambda_{\pi\theta}$ obtained from the monoids
morphism   $ \theta_{\crys}^* : \theta_{\crys}^{-1} \CN_{Y/\Bbb Z_p} \to  \CN'_{Y'/\Bbb Z_p}$. The identification $\epsilon_{al^m}(V^{(0)}, \varphi ,N)\iso \epsilon_{l}(V^{(0)}, \varphi ,mN)$ comes from the isomorphism of  $(1+\CJ_{Y/\Bbb Z_p})^\times$-torsors $\lambda_l^m \iso \lambda_{al^m}$ obtained from the map $\CN_{Y/\Bbb Z_p} \to
\CN_{Y/\Bbb Z_p}$, $\delta \mapsto [a]\delta^m$; here $[a]\in W(k)^\times$ is the Teichm\"uller lifting of $a$ (we use Remark (ii) in 1.12).

\enspace
(iii) (a) {\it Let us check the assertion in case when} $\pi =\id_{(S^0 ,\CL^0 )_1}$:

Let $R$ be $p$-adic completion of the divided powers algebra $W(k)\langle t\rangle$, $\CL_E$ 
be  the log structure on $E:=\Spec\, R$ generated by $t$. One has an exact embedding $i:  (S^0 ,\CL^0 )_1  \hra (E,\CL_E )$, $i^* (t )= l=\bar{p}$, and
a  $\varphi$-action $\varphi_R =\varphi$ on $(E,\CL_E )$, $\varphi^*_R (t )=t^p$,  that extends the Frobenius endomorphism of $ (S^0 ,\CL^0 )_1 $. Then $(E ,\CL_{E})_n$ (see 1.12) is a coordinate pd-$W_n (k)$-thickening of $ (S^0 ,\CL^0 )_1$. The trivialization $t$ of $ \lambda_{\bar{p}}$ yields an isomorphism $\epsilon  (V^{(0)}) (E,\CL_E ) \iso V^{(0)}_R := V^{(0)}\otimes_{W(k)} R$, $\epsilon :=\epsilon_{\bar{p}}$, that identifies the connection $\nabla$ on $ 
\epsilon  (V^{(0)}) (E,\CL_E )$ (see the end of 1.7) with the connection with potential  $Nd\log t$ on $V^{(0)}_R$, i.e., $\nabla (v\otimes r)= N(v)\otimes r d\log t + v\otimes dr$.

{\it $\epsilon$ is fully faithful:} It suffices to show that for two $(\varphi ,N)$-modules $V_1$, $V_2$ the map $\epsilon: R\Hom_{\varphi ,N} (V_1 ,V_2 )\to R\Hom_{\varphi}(\epsilon  (V_1 ), \epsilon  (V_2 ))$ is a quasi-isomorphism. Since $\epsilon  (V_i^{(0)})$ are locally free $\CO_{S^0_1/\Bbb Z_p}$-modules, 
  Lemma in 1.14 combined with (1.8.1) shows that
$R\Hom_{\varphi}(\epsilon  (V_1 ), \epsilon  (V_2 ))$ is log de Rham complex of $\CC one(\CH om (\epsilon  (V_1^{(0)}),\epsilon  ( V^{(0)}_2 ))\to \CH om (\varphi_{\crys}^*\epsilon (V_1^{(0)}),\epsilon  ( V^{(0)}_2 ))_{(Y,E)}[-1]$ tensored by $\Bbb Q$, i.e.,  the cone of endomorphism $\nabla_{t\partial_t}$ of the complex
$\CC one ((V_1^* \otimes V_2 )\otimes R \to (\varphi^* V_1^* \otimes V_2 ) \otimes R)\otimes\Bbb Q [-1]$. By Exercise above, $R\Hom_{\varphi ,N} (V_1 ,V_2 )$ is the cone of endomorphism $N$ of the complex $\CC one (V_1^* \otimes V_2 \to \varphi^* V_1^* \otimes V_2 )[-1]$. The map $\epsilon$ is the evident embedding of 
the latter complex into  the former one. The 
$\varphi$-action on $R$ satisfies the conditions of Proposition in 1.13.\footnote{To check that the $\varphi$-action on $I=t R$ is strongly topologically nilpotent, use the finite filtration $I\supset t I \supset t^2 I \supset \ldots \supset t^n I$, where $n$ is sufficiently large.  } Thus  our embedding is a quasi-isomorphism by (1.13.1), and we are done.
 
 {\it $\epsilon$ is essentially surjective}: By the previous paragraph, it suffices to check that for every $(\CF ,\varphi_{\CF})\in D^{\pcr}_{\varphi} (S^0_1 )^{\nd}$ the cohomology $H^i \CF$ belong to $D^{\pcr}_{\varphi} (S^0_1 )^{\nd}\otimes \Bbb Q$ and lies in the essential image of $\epsilon $.
Consider $\CF (E,\CL_E ):=\holim_n \,\CF ((E,\CL_E )_n )$. This is a perfect $R$-complex equipped with a nondegenerate $\varphi$-action, so  $\CF (E,\CL_E )\in D^{\perf}_{\varphi}(R)^{\nd}$.  By
(1.13.1), $H^i \CF (E,\CL_E )\otimes \Bbb Q$, viewed as an $R$-module equipped with $\varphi$-action, can be written in a canonical way as $V\otimes_{W(k)}R$, where $V$ is an effective $\varphi$-module. Then $\nabla_{t\partial_{t}}\in\End_{K_0}H^i \CF (E,\CL_E )\otimes \Bbb Q $ preserves $V$ (as follows, say, from part (a) of the proof of loc.~cit.); set $N_V :=  \nabla_{t\partial_{t}}|_V$. Then $(V,\varphi_V , N_V )$ is an effective $(\varphi ,N)$-module and $ H^i \CF =\epsilon_\pi (V)$, q.e.d.

(b) {\it The general case:}
The conditions on $Y$ assure that high enough power of Frobenius $Fr_{(Y,\CN )}^m$ factors through the embedding $i: (Y_{\text{red}}, \CN_{Y_{\text{red}}})\hra (Y,\CN )$, i.e., we have $\psi : (Y,\CN ) \to  (Y_{\text{red}},\CN_{Y_{\text{red}}})$ such that $i\psi = Fr_{(Y,\CN )}^m$, $\psi i= Fr_{ (Y_{\text{red}}, \CN_{Y_{\text{red}}})}^m$.     This implies  that $Li^*_{\crys}: D^{\pcr}_{\varphi} (Y )^{\nd}\otimes\Bbb Q \iso
D^{\pcr}_{\varphi} (Y_{\text{red}} )^{\nd}\otimes\Bbb Q$ by Remark (i) in in 1.14. Thus, by (ii), it is enough to prove assertion (iii) for $(Y_{\text{red}}, \CN_{Y_{\text{red}}})$. Equivalently, it is enough to consider the case when $Y$ is reduced. Then
 there is an isomorphism $\pi_0 : (Y,\CN )\iso (S^0 ,\CL^0 )_1$ of log $k$-scheme that corresponds to a generator $l_0$ of $\CN$. One has $l=al_0^m$ for some $a\in k^\times$, $m\in \Bbb Z_{>0}$. Since $\epsilon_{l_0}$ is an equivalence by (a), (ii) shows that $\epsilon_l$ is an equivalence as well, q.e.d.
\qed\enddemo

Let us describe functor (1.15.1) on the level of sections. Suppose our 
 $Y$ is affine, and let $(T,\CN_T )$ be a $p$-adic pd-thickening of $(Y,\CN )$, i.e., a sequence of objects  $(Y,T_n )$ in $(Y,\CN )_{\crys}$ with $T_{n}=T_{n+1}\otimes\Bbb Z/p^{n}$. Thus we have a $p$-adic algebra $A=\limleft A_n$, $A_n :=  A/p^n$, equipped with a pd-ideal $J$ such that $Y=\Spec\, A/J$, $T_n =\Spec\, A_n$. A complex $\CF$ of $\CO_{Y/\Bbb Z_p}$-modules on $(Y,\CN )_{\crys}$ yields a complex of $A$-modules $\CF (Y,T):= \holim_n \CF (Y,T_n )$.
 
 The $(1+\CJ_{T_n})^\times$-torsors $\lambda_{\pi (Y,T_n)}$, see part (i) of the proof above, are trivial (since $T_n$ is affine), so $\lambda_{ A} := \limleft \Gamma (T_n , \lambda_{\pi (Y,T_n)})$ is a $(1+J)^\times$-torsor. Let $\tau_{ A_{\Bbb Q}}$ be the  $A_{\Bbb Q}$-torsor obtained from 
 $\lambda_{A} $ by the pushout by  $(1+J)^\times \buildrel{\log}\over\lra J \to A_{\Bbb Q}$, $A_{\Bbb Q}:=A\otimes\Bbb Q$. We call $\tau_{ A_{\Bbb Q}}$, as well as the $\Bbb G_a$-torsor $\Spec\, A_{\Bbb Q}^\tau$ over $\Spec\, A_{\Bbb Q}$ with sections $\tau_{ A_{\Bbb Q}}$, the {\it Fontaine-Hyodo-Kato} torsor. Let $N_{\tau}$ be the $ A_{\Bbb Q}$-derivation of 
 $ A^\tau_{\Bbb Q}$ which is the action of the generator of Lie$_{\Bbb G_a}$.

 For a $(\varphi ,N)$-module $V$ the group $\Bbb G_a$ acts on the vector space  $V$ so that $N_V$ is the action of the generator of Lie$_{\Bbb G_a}$. Let $V^\tau_{A_{\Bbb Q}}$ be the $\tau_{ A_{\Bbb Q}}$-twist of $V_{A_{\Bbb Q}}:=V\otimes_{K_0} A_{\Bbb Q}$. Thus $V^\tau_{A_{\Bbb Q}}$ is the $A_{\Bbb Q}$-module of maps $v: \tau_{ A_{\Bbb Q}} \to V_{A_{\Bbb Q}}$ such that $v(\tau + f)= \exp (f N_V )(v(\tau))$ for any $\tau \in \tau_{ A_{\Bbb Q}}$, $f\in A_{\Bbb Q}$. Equivalently,  $V^\tau_{A_{\Bbb Q}}=(V\otimes_{K_0} A_{\Bbb Q}^\tau )^{\Bbb G_a}= (V\otimes_{K_0} A_{\Bbb Q}^\tau )^{N=0}$; here $N$ acts on the tensor product by the Leibniz rule $N:= N_V \otimes 1 +1\otimes N_\tau$. 
  It follows directly from the construction that $$\epsilon_\pi (V)(Y,T)= V^\tau_{A_{\Bbb Q}}. \tag 1.15.3$$   
 
\remark{Remarks} (i) One has a natural map $a: \tau_{A_{\Bbb Q}}\to 
A_{\Bbb Q}^\tau$ of $A_{\Bbb Q}$-torsors  which assigns to
$\tau \in \tau_{A_{\Bbb Q}}$ a function $a_\tau \in A_{\Bbb Q}^\tau$  whose value on any $\tau' \in\tau_{ A_{\Bbb Q}}$, viewed as a section $\Spec\, A_{\Bbb Q} \to \Spec\, A_{\Bbb Q}^\tau$, is
 $\tau -\tau' \in A_{\Bbb Q}$. The algebra $A_{\Bbb Q}^\tau$ is freely generated over $A_{\Bbb Q}$ by $a (\tau )$ for any $\tau$, and the derivation $N_{\tau}$ is characterized by property $N_{\tau}(a_\tau )=1$. 
\newline
(ii)  By (1.15.2), for every lifting $\varphi_T$ of Frobenius endomorphism of $(Y,\CN )$ to $(T,\CN_T )$ the corresponding endomorphism $\varphi_{A_{\Bbb Q}}$ of $A_{\Bbb Q}$ extends canonically to an endomorphism $\varphi_\tau$ of $A^\tau_{\Bbb Q}$ such that $N_\tau \varphi_\tau =p\varphi_\tau N_\tau$. Hence we get a Frobenius action on  $V^\tau_{A_{\Bbb Q}}$, and (1.15.3) is compatible with the Frobenius actions.   \newline
(iii) We can view $\lambda_A$ as the preimage of $l$ by the map $\Gamma (T, \CN_T^{\text{gr}})\to\Gamma (Y,\CN^{\text{gr}})$ or by $\Gamma (T, \CN_T^{\text{gr}}/k^\times)\to\Gamma (Y,\CN^{\text{gr}}/k^\times )$,\footnote{Here $k^\times$ acts on $\CN_T$ via the Teichm\"uller section $k^\times \to W(k)^\times$ and Remark (ii) in 1.12.  } 
so $\tau_{A_{\Bbb Q}}$ is the preimage of $l$ by $\Gamma (T, \CN_T^{\text{gr}}/k^\times)\otimes\Bbb Q \to\Gamma (Y,\CN^{\text{gr}}/k^\times)\otimes\Bbb Q$. Thus every $t \in  \Gamma (T, \CN_T )$ whose image in $\Gamma (Y,\CN^{\text{gr}}/k^\times)\otimes\Bbb Q$ equals $l^{v(t)}$, $v(t)\in \Bbb Q^\times$, yields a trivialization
$t^{v(t)^{-1}}$ of  $\tau_{A_{\Bbb Q}}$.
\newline
(iv) Let $(Y,\CN )$ be as in (iii) of the theorem. Let 
$i$ be the embedding $ (\Spec\, k, \CL' )\hra (Y,\CN )$, $\CL' :=\CN |_{\Spec\, k}$, and $\CL\,\tilde{}$ be the log structure on $\Spec\, W(k)$ 
defined by the prelog one $\CL' \to k \to W(k)$, the right arrow is the Teichm\"uller section, so $(\Spec\, W(k) , \CL\,\tilde{}\, )$ is a $p$-adic thickening of $  (\Spec\, k, \CL' )$.
Then for $\CF\in D^{\pcr}_{\varphi} (Y)^{\nd}$, the complex $\epsilon_\pi^{-1} (\CF )$  as an object of $D_{\varphi}(W(k))\otimes\Bbb Q$ is equal to $Li^{*}_{\crys} (\CF )(\text{Spec}\, W(k),\CL\,\tilde{}\,)$. Indeed, by (ii), (iii) of the theorem it is enough to check the assertion for $Y=\Spec\, k$ and $\CF = \epsilon_\pi (V)$. Then our
identification $V \iso \epsilon_\pi (V)(\Spec\, W,\CL\,\tilde{}\, )$  comes from the trivialization of
 the torsor $\lambda_l (\Spec\, W , \CL\,\tilde{}\, )$  by the image of $l$ by  $\CN =\CL' \to \CL\,\tilde{}$. \endremark

\subhead{\rm 1.16}\endsubhead {\it Hyodo-Kato theory.} Let $K$ be a $p$-adic field, i.e., a complete discretely valued field of characteristic 0 with perfect residue field $k=O_K/\fm_K$. Set $K_0 =\text{Frac}\, W(k)\subset K$. Let $v$ be the valuation on $K$ normalized so that $v(p)=1$.
Let  $\CL =\CL_{K} $ be  {\it the canonical } log structure on $S=S_{K}:=\Spec\, O_{K} $ generated by  prelog one $ O_{K} \smallsetminus \{ 0\} \to O_{K}$. Set $S^0 :=\Spec\, W(k)$ and let $\CL^0$ be the canonical log structure on $S^0$; we have the evident map $\pi : (S,\CL )\to (S^0 ,\CL^0 )$.  
As in 1.12, we have log $W_n (k)$-schemes $(S ,\CL )_n /(S^0 ,\CL^0 )_n$, etc. 
Since  $\pi_1 : (S,\CL )_1 \to (S^0 ,\CL^0 )_1$ satisfies the condition of (iii) in Theorem in 1.15, we have  the equivalences of categories  
$\epsilon : D_{\varphi ,N}(K_0 )^{\eff}\iso D^{\pcr}_{\varphi} ((S^0 ,\CL^0 )_1  )^{\nd}\otimes\Bbb Q$ and  
$\epsilon_{\pi_1}=L\pi_{1\,\crys}^* \,\epsilon : D_{\varphi ,N}(K_0 )^{\eff}\iso D^{\pcr}_{\varphi} ((S,\CL)_1 )^{\nd}\otimes\Bbb Q$.

Let $f: (Z_1 ,\CM_1 )\to (S_1 ,\CL _1)$ be a log smooth map of Cartier type with
$(Z_1 ,\CM_1 )$  fine and
$Z_1$ is proper over $S_1$.  
By the theorem in 1.14,\footnote{applied to $Z=Z_1$, $Y=S_1$, and $S=$Spf$\, (\Bbb Z_p )$.} $Rf_{\crys *}(\CO_{Z_1 /\Bbb Z_p})$ is a nondegenerate perfect F-crystal on $(S ,\CL )_{1\,\crys}$. Set $$R\Gamma_{\!\text{H}\!\text{K}}(Z_1,\CM_1 ):= \epsilon_{\pi_1}^{-1} Rf_{\crys *}(\CO_{Z_1 /\Bbb Z_p})\otimes\Bbb Q     \in D_{\varphi ,N}(K_0 ). \tag 1.16.1$$ This is the {\it Hyodo-Kato complex} of $(Z_1,\CM_1 )$.

\remark{Remarks} (i) Let $\CL_{01}  $ be the restriction of $\CL_1$ to $S_{ 01}  :=\Spec \, k \subset S_1$, so we have $i: (S_{ 01}  ,\CL_{ 01}  )\hra (S_1 ,\CL_1 )$ and $\pi_{ 01}  :=i\pi_1$. Let $f_0 : (Z_{ 01}  ,\CM_{ 01}  )\to (S_{ 01}  ,\CL_{ 01} )$ be the $i$-pullback of $f$. Then $Rf_{0\crys *} \CO_{Z_{ 01}  /\Bbb Z_p}) = Li^*_{\crys} Rf_{\crys *} \CO_{Z_1 /\Bbb Z_p}) $ by base change (1.11.1). Thus, by (ii) of
 the theorem in 1.15, one has  $R\Gamma_{\!\text{H}\!\text{K}}(Z_1,\CM_1 )= \epsilon_{\pi_{ 01} }^{-1} Rf_{0\crys\, *}(\CO_{Z_{ 01}  /\Bbb Z_p})\otimes\Bbb Q $, so the Hyodo-Kato complex depends only on $f_0$. \newline
 (ii) Remark (iv) in 1.15 implies that the above definition amounts to the original one from \cite{HK}  \S 3.
\endremark

\enspace

The Hyodo-Kato complex controls, up to isogeny, the relative log crystalline cohomology for all
base changes of $f$. Namely, suppose we have an  integral quasi-coherent log scheme $ (Y ,\CN )$ with $Y$ affine, a
 map $\theta:  (Y ,\CN )\to (S ,\CL )_1$, and a $p$-adic pd-thickening $(T,\CN_T )$ of $(Y,\CN )$ as in the end of 1.15, so $T=\Spec\, A$, $Y=\Spec\, A/J$. 
 Let $f_Y : (Z_{1Y},\CM_{1Y})\to  (Y ,\CN )$ be the $\theta$-pullback of $f$.

\proclaim{Theorem} (i) The $A$-complex $R\Gamma_{\!\crys}(( Z_{1Y} /T)^{\log}_{\crys} ,\CO_{Z_{1Y} /T})$ is perfect, and one has  $R\Gamma_{\!\crys}(( Z_{1Y} /T_n )^{\log}_{\crys} ,\CO_{Z_{1Y} /T_n })=R\Gamma_{\!\crys}(( Z_{1Y}/T)^{\log}_{\crys} ,\CO_{Z_{1Y}/T})\otimes^L \Bbb Z/p^n$.
\newline
(ii) There is a canonical Hyodo-Kato quasi-isomorphism of $A_{\Bbb Q}$-complexes
$$\iota : R\Gamma_{\!\text{H}\!\text{K}}(Z_1,\CM_1 )^\tau_{A_{\Bbb Q}}\iso R\Gamma_{\!\crys}(( Z_{1Y} /T )^{\log}_{\crys} ,\CO_{Z_{1Y}}/T)\otimes\Bbb Q.\tag 1.16.2$$  If we have a Frobenius lifting   $\varphi_T$, then $\iota$ commutes with its action.
\endproclaim

 \demo{Proof}  (i) By the theorem in 1.11 and Exercise in loc.~cit., $
R f_{Y\crys\, *}(\CO_{Z_{1Y} /\Bbb Z_p})$ is a perfect $\CO_{Y/\Bbb Z_p }$-crystal. Its value on $(Y , T_n )$ equals $R\Gamma_{\!\crys}(( Z_{1Y} /T_n )^{\log}_{\crys} ,\CO_{Z_{1Y} /T_n})$. So the latter complex is  $A_n$-perfect and $
R\Gamma_{\!\crys}(( Z_{1Y} /T_{n+1} )^{\log}_{\crys} ,\CO_{Z_{1Y} /T_{n+1}})
\otimes^L_{A_{n+1}}A_n =R\Gamma_{\!\crys}(( Z_{1Y} /T_n )^{\log}_{\crys} ,\CO_{Z_{1Y} /T_n})$. Since one has $R\Gamma_{\!\crys}(( Z_{1Y} /T)^{\log}_{\crys} ,\CO_{Z_{1Y} /T})
=\holim_n R\Gamma_{\!\crys}(( Z_{1Y} /T_n )^{\log}_{\crys} ,\CO_{Z_{1Y} /T_n})
$, we are done.

(ii) One has $\epsilon_{\pi_1} R\Gamma_{\!\text{H}\!\text{K}}(Z_1,\CM_1 )= Rf_{\crys *}(\CO_{Z_1 /\Bbb Z_p})\otimes\Bbb Q$. So, by 
(ii) in Theorem in 1.15 and
base change (see 1.11), $\epsilon_{\pi_1 \theta} R\Gamma_{\!\text{H}\!\text{K}}(Z_1,\CM_1 )= Rf_{Y\crys *}(\CO_{Z_{1Y} /\Bbb Z_p})\otimes\Bbb Q$. We get (1.16.2) evaluating this crystal on $(Y,T)$ and using (1.15.3). For the Frobenius action, see Remark (ii) in 1.15. \qed
 \enddemo

\remark{Example} Let $(Z,\CM )$ be a fine log scheme log smooth over $(S,\CL)$ such  that $Z$ is proper over $S$ and $(Z,\CM )_1 =(Z_1 ,\CM_1 )$ is of Cartier type over $(S,\CL )_1$. Consider $\theta =\id_{(S ,\CL )_1}$ and $(T ,\CN_T )=(S,\CL )$, so $A = O_K$. Then
  $R\Gamma ((Z_1 /S_n )^{\log}_{\crys}, \CO_{Z_1 /S_n})=R\Gamma (Z_n, \Omega^\cdot_{(Z ,\CM )_n /(S,\CL)_n} )$ by (1.8.3) and the
invariance property of crystalline topology. Thus $R\Gamma ((Z_1 /S )^{\log}_{\crys}, \CO_{Z_1 /S})=R\Gamma (Z, \Omega^\cdot_{(Z ,\CM ) /(S,\CL)} )$, and (1.16.2) becomes
 $$\iota_{\dR } :  R\Gamma_{\!\text{H}\!\text{K}}(Z_1,\CM_1 )_K^\tau  \iso R\Gamma (Z_{K}, \Omega^\cdot_{(Z_{K} ,\CM_{K} ) /K} ), \tag 1.16.3$$ where $Z_{K}:=Z\otimes_{O_{K}}\! K$. A choice of a non-zero element $t$ in the maximal ideal of $O_K$ yields, as in Remark (iii) in 1.15,  a trivialization of $\tau_K$, hence an identification $R\Gamma_{\!\text{H}\!\text{K}}(Z_1,\CM_1 )_K \iso R\Gamma_{\!\text{H}\!\text{K}}(Z_1,\CM_1 )_K^\tau$; its composition with (1.16.3) is
 the classical Hyodo-Kato quasi-isomorphism from \cite{HK} 5.1.\endremark

\subhead{\rm 1.17}\endsubhead {\it Absolute crystalline cohomology of $O_{\bar{K}}/p$.}   Let $(Y,\CN )$ be an integral quasi-coherent log scheme over $\Bbb Z/p$ such that the monoid $\CN /\CO_Y^\times$ is uniquely $p$-divisible. The next lemma was pointed out by the referee:

\proclaim{Lemma} One has $Y_{\crys}=(Y,\CN )_{\crys}$, i.e.,  every pd-thickening $P$ of $Y$ over $\Bbb Z/p^n$  carries a unique log structure $\CN_P$. \endproclaim 

\demo{Proof}  Let us construct $\CN_P$; its uniqueness follows from Exercise below.  Let $\CN^{(n)}$ be a copy of $\CN$, and 
$\alpha^{(n)}: \CN^{(n)}\to \CO_P$ be the map $a\mapsto \tilde{\alpha} (a)^{ p^n}$, where $\tilde{\alpha} (a)$ is any lifting of $\alpha (a)\in \CO_Y$ to $\CO_P$. This map is well defined (indeed,  for  $b\in  \CO_Y$ and $c\in \CJ_{P}$ one has $(b+c)^{p^n} =b^{p^n}$ since $c^p \in p\CJ_P$). It is a prelog structure on $P$; our $\CN_P$ is the corresponding log structure. The map $\CN^{(n)} \to \CN$, $a \mapsto a^{p^n}$, lifts the embedding $Y\hra P$ to an embedding of log schemes $(Y,\CN )\hra (P,\CN_P )$, which is exact due to the condition of the lemma. 
\qed \enddemo

\remark{Exercise} Show that $\CN_P$ satisfies the next universal property: For any $(Z,T, \CM_T)\in \CT_{\Bbb Z/p^n}$ (see 1.3) every map of log schemes $ (Z,\CM )\to (Y,\CN )$ and a map of pd-thickenings $T\to P$ that restrict to the same map $Z\to Y$ lift in a unique manner to a morphism log pd-thickenings $(Z,T,\CM_T )\to (Y, P,\CN_P )$.
\endremark

\enspace 

As above, set $K_0 =\text{Frac}\,   W(k)$, $k$ is a perfect field. Let $\bar{K}$ be an algebraic closure of $K_0 $, $O_{\bar{K}}$ be its ring of integers, $\bar{k}= O_{\bar{K}}/\fm_{\bar{K}}$  the residue field. Let $\bar{\CL}$ be the canonical log structure on $\bar{S}:=\Spec\,O_{\bar{K}}$ generated by the prelog structure $O_{\bar{K}}\smallsetminus \{ 0\} \to O_{\bar{K}}$. Let $v$ be the normalized valuation on $\bar{K}$, $v(p)=1$, so we have $v: \bar{\CL} / O_{\bar{K}}^\times =  \bar{\CL}_1 / (O_{\bar{K}}/p)^\times \iso \Bbb Q_{\ge 0}$.

Consider Fontaine's ring $\Acrys$  from \cite{F1} 2.2, 2.3. This is a p-adically complete ring such that  A$_{\crys\, n}= \text{A}_{\crys} /p^n$ is a universal pd-thickening of $O_{\bar{K}}/p$ over $\Bbb Z/p^n $. Let $J_{\crys\, n} $ be the pd-ideals, A$_{\crys\, n}/ J_{\crys\, n} =O_{\bar{K}}/p$. Set $E_{\crys\, n}:=\Spec\, \text{A}_{\crys\, n}$.

  The log structure $\bar{\CL}_{1}$ on $\Spec\, ( O_{\bar{K}}/p) = \bar{S}_1$ satisfies the condition of the  lemma. So it extends in a unique manner to an integral log structure $\CL_{\crys\, n} $ on $E_{\crys\, n}$. By Exercise, the pd-thickening $(\bar{S},\bar{\CL})_1 \hra (E_{\crys\,n}, \CL_{\crys\, n})$ is universal, i.e.,  for $(Z,T)\in \CT_{\Bbb Z/p^n}$ every map $h : (Z,\CM )\to (\bar{S},\bar{\CL})_1$ of log  schemes extends in a unique way to a $\CT_{\Bbb Z/p^n}$-map $h_T : (Z,T)\to (S_1 ,E_{\crys\,n})$.

The Frobenius map $\varphi$ lifts  to $(E_{\crys\,n} ,  \CL_{\crys\, n} )$  by universality. 

The log structures $\CL_{\crys\, n} $ are mutually compatible. Thus we have the log structure $\CL_{\crys}=\limleft  \CL_{\crys\, n}$ on $E_{\crys} =\Spec\, \Acrys$ such that $ \CL_{\crys}|_{E_{\crys\,n}} =\CL_{\crys\, n} $. Explicitly,  the identifications $\CL_{\crys\, n}|_{E_{\crys\,m}}$ $\iso \CL_{\crys\, m}$
for $n\ge m$ come from the maps $\bar{\CL}^{(n)}_{1}\to \bar{\CL}^{(m)}_{1}$, $l^{(n)}\mapsto \varphi^{n-m} (l^{(n)})$, and the log structure  $\CL_{\crys}$ comes from the prelog one $\alpha_\varphi = \limleft\alpha^{(n)} : \CL_\varphi := \limleft  \bar{\CL}^{(n)}_{1}  \to \Acrys$.

\remark{Exercise} One has $\CL_\varphi \iso \{ \nu\in \CL_{\crys}:\varphi (\nu )=\nu^p \}$.\footnote{Hint: $\dim_{\Bbb Q_p} \{ b\in \Bcrys^+ : \varphi (b)=pb,  b\in F^1\}=1$.}
\endremark

\enspace

Our $(\bar{S}_1 ,E_{\crys\, n})$ is a final object of the absolute crystalline site $(\bar{S}_1,\bar{\CL}_1 )_{\crys (n)}$ (see 1.12) by Lemma, so for any sheaf $\CF$ its global sections  are equal to $\CF (\bar{S}_1 ,E_{\crys \,n} )$. The \'etale topology of $E_{\crys\, n}$ is trivial, so the higher cohomology vanish and $$R\Gamma ((\bar{S}_1 ,\bar{\CL}_1  )_{\crys (n)}, \CF )=\CF (\bar{S}_1 ,E_{\crys\, n} ). \tag 1.17.1$$ Thus for a  sheaf $\CF$ on $(\bar{S}_1 ,\bar{\CL}_1  )_{\crys}$ one has $$R\Gamma ((\bar{S}_1 ,\bar{\CL}_1  )_{\crys}, \CF )=\CF (\bar{S}_1 ,E_{\crys} ):=\holim_n \,\CF (\bar{S}_1 ,E_{\crys\, n} ). \tag 1.17.2$$
 In particular, $$R\Gamma_{\!\crys}(\bar{S}, \bar{
\CL})_{(n)} = \text{A}_{\crys \, n}, \quad R\Gamma_{\!\crys}(\bar{S}, 
\bar{\CL}) = \text{A}_{\crys }. \tag 1.17.3$$

Recall that $\Bcrys^+ := \text{A}_{\text{crys}\,\Bbb Q} $. The $\Bcrys^+$-algebra $\text{A}^\tau_{\text{crys}\,\Bbb Q}$, defined by the projection $(\bar{S},\bar{\CL})_1 \to (S^0 ,\CL^0 )_1$ (see 1.15), equals Fontaine's ring $\Bst^+$. To see this, consider the subset
$\CL_{\varphi \bar{p}}$ of $\CL_\varphi$ formed by  sequences $(l^{(n)})$ that start with $l^{(0)}=\bar{p}$, i.e., $\CL_{\varphi \bar{p}}$ is the intersection of $\CL_\varphi$ and $\lambda_{\Acrys}$ in $\CL_{\crys}$. Let $\chi: \CL_{\varphi \bar{p}}\to \text{A}^\tau_{\text{crys}\,\Bbb Q}$ be the composition $\CL_{\varphi \bar{p}}\hra \lambda_{\Acrys}\to \tau_{\text{A}_{\text{crys}\,\Bbb Q}}\to 
\text{A}^\tau_{\text{crys}\,\Bbb Q}$, the second arrow is the canonical map that defines $\tau_{\text{A}_{\text{crys}\,\Bbb Q}}$, the third one is the map $a$ from Remark (i) in 1.15. It follows from loc.cit.~that there is a unique morphism of $\Bcrys^+$-algebras $\text{A}^\tau_{\text{crys}\,\Bbb Q}\to \Bst^+$ that sends  $\chi (l^{(n)})$ to the element $\lambda (l^{(n)})$ from  \cite{F1} 3.1.4, and that morphism is an isomorphism. It evidently commutes with the $\varphi$- and $N$-actions (see \cite{F1} 3.2.1, 3.2.4).

\subhead{\rm 1.18}\endsubhead {\it Absolute crystalline cohomology of log schemes over $O_{\bar{K}}/p$.} 
Let $\bar{f}: (Z\bar{_1},\CM\bar{_1})$ $\to (\bar{S}_1 , \bar{\CL}_1  )$
be a map of log schemes, its source is integral quasi-coherent. Then for any sheaf $\CF$  on $(Z\bar{_1} ,\CM\bar{_1 } )_{\crys}$ (see 1.12)  one has a natural identification $$R\Gamma (
(Z\bar{_1},\CM\bar{_1})_{\crys},\CF )\iso R\Gamma (((Z\bar{_1},\CM\bar{_1})/(\bar{S}_1 ,E_{\crys} ))_{\crys},\CF ). \tag 1.18.1$$ Namely, consider  $\bar{f}_{\crys} :(Z\bar{_1},\CM\bar{_1} )_{\crys}\!\!\tilde{}\,\, \to (\bar{S}_1 ,\bar{\CL}_1  )_{\crys}\!\!\tilde{}\,\,$; then the l.h.s.~of (1.18.1) is $
R\Gamma  (  (\bar{S}_1 ,\bar{\CL}_1  )_{\crys}  ,  R\bar{f}_{\crys*}\CF )$, the same is true for the r.h.s.~by (1.17.2).

\enspace 

Let  $K\subset \bar{K}$ be a finite extension of $K_0$, and
$\theta =\theta_{K} : (\bar{S}, \bar{\CL} )\to (S,\CL )=(S_{K},\CL_{K})$ be the map defined by the embedding $K\hra \bar{K}$. 
Suppose our $\bar{f}$ is the base change of $f:  (Z_1 ,\CM_1 )\to (S,\CL )_1$ by $\theta_1$,  i.e., we have  $\theta_{Z1}  : (Z\bar{_1} ,\CM\bar{_1} )\to (Z_1 ,\CM_1)$ such that the square $(\bar{f},f,\theta_1 ,\theta_{Z1} )$ is Cartesian. Suppose, as in 1.16, that $f$ is log smooth of Cartier type, $(Z_1 ,\CM_1 )$ is fine, and $Z_1$ is proper over $S_1$. 
Applying the  theorem in 1.16 with $Y=\bar{S}_1$, $T=E_{\crys}$ and 
using (1.18.1) for $\CF = \CO_{Z\bar{_1}/\Bbb Z_p}$, we get (recall that $\Bcrys^+ :=\text{A}_{\text{crys}\,\Bbb Q} $; see 1.12 for the rest of the notation):

\proclaim{Theorem} (i) $R\Gamma_{\!\crys} (Z\bar{_1} ,\CM\bar{_1} )$ is a perfect $\Acrys$-complex, and
$R\Gamma_{\!\crys} ( Z\bar{_1} ,\CM\bar{_1})_{(n)} =   R\Gamma_{\!\crys} (
Z\bar{_1} ,\CM\bar{_1} )\otimes^L \Bbb Z/p^n = R\Gamma_{\!\crys} (Z\bar{_1} ,\CM\bar{_1}
 )\otimes_{\Acrys}\! \text{A}_{\crys\, n}$.
 \newline (ii) There is a canonical  quasi-isomorphism of $\Bcrys^+$-complexes $$\iota_{\crys} :
 R\Gamma_{\!\text{H}\!\text{K}} (Z_1,\CM_1 )_{\Bcrys^+}^{\tau} \iso R\Gamma_{\!\crys}(Z\bar{_1} ,\CM\bar{_1} )\otimes\Bbb Q  \tag 1.18.2 $$ compatible with the action of $\varphi$. \qed
\endproclaim

\remark{Remarks} (i) Suppose our $\bar{f}$ is such that the datum of $(K, (Z_1 ,\CM_1 ),\theta_{Z1} )$ as above exists, but we don't want to specify one. All such data form a category in an evident manner; since $\bar{S}_1$ is faithfully flat over $S_1$, this is an ordered set, which we denote by $\Xi_1$. In fact, $\Xi_1$ is directed. For a morphism $(K', (Z'_1 ,\CM'_1 ),\theta'_{Z1} )\to (K, (Z_1 ,\CM_1 ),\theta_{Z1} )$ in $\Xi_1$ one has a canonical base change identification of the Hyodo-Kato complexes  $ R\Gamma_{\!\text{H}\!\text{K}} (Z_1,\CM_1 )\otimes_{K_0}K'_0 \iso R\Gamma_{\!\text{H}\!\text{K}} (Z^{\prime }_1,\CM^{\prime }_1 )$, and isomorphisms (1.18.2) are compatible with it. Set $$R\Gamma_{\!\text{H}\!\text{K}}(Z\bar{_1} ,\CM\bar{_1} ):= \limright_{\,\Xi_1}\! R\Gamma_{\!\text{H}\!\text{K}} (Z_1,\CM_1 ); \tag 1.18.3$$ this is a complex of $(\varphi ,N)$-modules over $K_0^{\text{nr}}$ (the maximal unramified extension of $K_0$ in $\bar{K}$) functorial with respect to morphisms of $X$'s,  and (1.8.3) provides a canonical isomorphism  compatible with the action of $\varphi$ $$\iota_{\crys}: R\Gamma_{\!\text{H}\!\text{K}}(Z\bar{_1} ,\CM\bar{_1} )^\tau_{\Bcrys^+} \iso 
R\Gamma_{\!\crys}(Z\bar{_1} ,\CM\bar{_1} )\otimes\Bbb Q.  \tag 1.18.4 $$ 
Since $\Spec\,\Bst^+$ is a $\Bbb G_a$-torsor over $\Spec\, \Bcrys^+$, (1.18.4) amounts to a quasi-isomorphism of $\Bst^+$-complexes compatible with the action of $N$ and $\varphi$  $$\iota_{ \crys }: R\Gamma_{\!\text{H}\!\text{K}}(Z\bar{_1} ,\CM\bar{_1} )_{\Bst^+} \iso 
R\Gamma_{\!\crys}(Z\bar{_1} ,\CM\bar{_1} )\otimes^L_{\Acrys}\Bst^+ . \tag 1.18.5$$ 
Here $ R\Gamma_{\!\text{H}\!\text{K}}(Z\bar{_1} ,\CM\bar{_1} )_{\Bst^+ }:= R\Gamma_{\!\text{H}\!\text{K}}(Z\bar{_1} ,\CM\bar{_1} )\otimes_{K_0^{\text{nr}}} \Bst^+ $. 

(ii) Suppose $(Z\bar{_1},\CM\bar{_1})/(\bar{S},\bar{\CL})_1$ is reduction mod $p$ of a log scheme $(Z\,\bar{},\CM\,\bar{}\,)/(\bar{S},\bar{\CL})$. Assume that there exists a datum $(K,(Z,\CM ),\theta_Z )$, where $K\subset \bar{K}$ is a finite extension of $K_0$, $(Z,\CM )/(S,\CL )$ is a log scheme that {\it satisfies conditions of Example in 1.16}, and $\theta_Z : (Z\,\bar{},\CM\,\bar{}\,)\to (Z,\CM )$ is an identification of $(Z\,\bar{},\CM\,\bar{}\,)$ with the $\theta$-pullback of $(Z,\CM )$. Again, such data form a directed set $\Xi$, and the reduction mod $p$ map $\Xi \to \Xi_1$ is cofinal. Isomorphisms (1.16.3) are compatible with morphisms in $\Xi$, and their $\Xi$-colimit is a natural isomorphism (here $Z\bar{_{\bar{K}}}:= Z\,\,\bar{}\otimes_{O_{\bar{K}}}\bar{K}$)
$$ \iota_{\dR}: R\Gamma_{\!\text{H}\!\text{K}}(Z\bar{_1} ,\CM\bar{_1} )_{\bar{K}}^\tau \iso R\Gamma (Z\bar{_{\bar{K}}}, \Omega^\cdot_{(Z\bar{_{\bar{K}}},\CM\bar{_{\bar{K}}})/\bar{K}}). \tag 1.18.6$$
\endremark

\subhead{\rm 1.19}\endsubhead {\it Log de Rham complex in characteristic 0.} 
Let $(Y,\CN)$ be an integral fine  log scheme log smooth over a field $F$ of characteristic 0, and $Y^0$ be the open subset where the log structure is trivial. The embedding $j: Y^0\hra Y$ is affine. 
The next result is due to Ogus \cite{Og2} 1.3; the key idea of the proof is borrowed from \cite{D}.

\proclaim{Theorem} If the sheaf of groups $\CN^{\text{gr}}/\CO^\times_{Y}$ has trivial torsion (e.g.~if $\CN$ is saturated), then
 the natural map $r_{Y}:\Omega^\cdot_{(Y, \CN )/F} \to j_* \Omega^\cdot_{Y^0/F}$ is a quasi-isomorphism. 
 \endproclaim
\demo{Proof} (i) We want to prove that $\CC one (r_Y )$ is acyclic. This is a complex of
quasi-coherent   $\CO_{Y}$-modules whose differentials are differential operators.
Let $C_Y$ be corresponding complex of induced $\CD$-modules on $Y$, see \cite{S}; here $\CD$-modules on a singular variety are understood  in the usual way (using closed embeddings into a smooth variety). By loc.~cit., $\CC one (r_Y )$ is quasi-isomorphic to 
 the de Rham complex dR$(C_Y )$, so it suffices to show that $C_Y$ is acyclic. We proceed by induction by $\dim Y$.

 (ii) Our claim is \'etale local, so, by \cite{K1} 3.5, we can assume that $Y=\Spec\, F [N]$ for a fine monoid $N$, $\CN$ comes from the  $N$-chart. Since $\CN^{\text{gr}}/\CO^\times_{Y}$ has trivial torsion, we can assume that $N^{\text{gr}}$ has trivial torsion. Then $T:= \Spec\, F[N^{\text{gr}}]$ is a torus which acts on $Y$, and $Y^0 =T$ is the open orbit. 
 
  (iii) $C_{Y}$ is acyclic outside $Y^T$: For a closed $y\in Y$ not fixed by $T$, let us find an \'etale neighborhood $U$ of $y$ such that $C_U$ is acyclic. Pick  $n\in N\smallsetminus\{ 0\}$ with $\alpha (n) (y)\neq 0$ and then a 1-parameter subgroup $G\subset T$ such that $n|_G$ is nontrivial; set $Z:=\alpha (n)^{-1}(\{1\})\subset Y$. Our $U$ is $ G\times Z\to Y$, $(g,z)\mapsto g(z)$. Since $(U,\CN_U)=G\times (Z,\CN|_Z )$,   $C_U$ is quasi-isomorphic to the pullback of $C_Z$ by the projection $U\to Z$; since $C_Z$ is acyclic by the induction assumption, $C_U$ is acyclic.
 
  (iv) By (iii), $C_Y = i_* Ri^! (C_Y )$ where $i:Y^T \hra Y$. Since  $Y^T$ is a single point if nonempty, $Ri^! (C_Y )=
 R\Gamma_{\!\dR}(Y, C_{Y})$, which is the cone of $\Gamma (r_Y): \Gamma (Y,\Omega^\cdot_{(Y,\CN )/F})\to \Gamma (T,\Omega^\cdot_{T /F})$ since $Y$ is affine. It remains to check that $\Gamma (r_Y)$ is a quasi-isomorphism.
 The $T$-action on $Y$ yields an $N^{\text{gr}}$-grading on the complexes. Since $\Omega^i_{T /F}=\CO_{T}\otimes \Lambda^i N^{\text{gr}}$,   the $n$-component  $\Gamma (T ,\Omega^\cdot_{T /F})_n $,  $n\in N^{\text{gr}}$,  equals $ F\otimes \Lambda^\cdot  N^{\text{gr}}$ with differential $\ell \mapsto n\wedge \ell$. Since $\Omega^i_{(Y,\CN )/F}=\CO_{Y}\otimes \Lambda^i N^{\text{gr}}$, the map $\Gamma (r_Y )$ is injective and its image is the sum of components $\Gamma (T ,\Omega^\cdot_{T /F})_n $ for $n\in N$.  We are done since $\Gamma (T ,\Omega^\cdot_{T /F})_n $ is acyclic for $n\neq 0$, hence for $n\in N^{\text{gr}}\smallsetminus N$. \qed\enddemo

\head 2. The $h$-sheaf $\CA_{\crys}$ and the crystalline Poincar\'e lemma. \endhead

\subhead{\rm 2.1}\endsubhead The next general format will be of use. Let  $\CV ar_F$ be the category of algebraic varieties over a field $F$ (i.e., separated reduced $F$-schemes of finite type). Let  $\CV ar_{F\h}$ be the corresponding h-site (see \cite{SV} or \cite{B} 2.4); for $X\in \CV ar_F$ let $X_{\h}$ be the  h-site of $X$.  For a complex $P$ of abelian groups we denote by $P_{\CV ar_{F\h}}$ the corresponding complex of constant h-sheaves; as in \cite{B} (1.1.1), $P\hotimes \Bbb Z_p := \text{holim}_n P\otimes^L \Bbb Z/p^n$, $P\hotimes \Bbb Q_p :=( P\hotimes \Bbb Z_p )\otimes\Bbb Q$. 

Let $\CA$ be a complex of h-sheaves on $\CV ar_F$. Set $A:= \CA (\Spec\, F)$; one has an evident canonical morphism $A_{\CV ar_{F\h}}\to \CA$. We say that $\CA$ is {\it $p$-adically constant} if the map
$(A\otimes^L \Bbb Z/p)_{\CV ar_{F\h}}= A_{\CV ar_{F\h}}\otimes^L \Bbb Z/p \to \CA\otimes^L \Bbb Z/p$ is a quasi-isomorphism in the derived category of h-sheaves $D(\CV ar_{F \h})$.

\remark{Remarks} (i) For a $p$-adically constant $\CA$ the map $(A\otimes^L \Bbb Z/p^n )_{\CV ar_{F\h}}\to \CA\otimes^L \Bbb Z/p^n$ is automatically a quasi-isomorphism for every $n>0$. \newline
(ii) $p$-adically constant complexes form a thick subcategory of  $D(\CV ar_{F \h})$. 
\endremark

\proclaim{Proposition} (i) For a $p$-adically constant $\CA$ as above and  $X\in \CV ar_F$, one has canonical quasi-isomorphisms  $$R\Gamma (X_{\h}, \CA )\otimes \Bbb Z/p^n \iso R\Gamma (X_{\et},\Bbb Z/p^n  )\otimes_{\Bbb Z/p^n}( A\otimes \Bbb Z/p^n ). \tag 2.1.1$$ (ii) If, in addition, $F$ is algebraically closed and $p$ is prime to char $F$, then one has a canonical quasi-isomorphism  $$R\Gamma (X_{\h}, \CA )\hotimes \Bbb Z_p \iso R\Gamma_{\!\et}(X,\Bbb Z_p )\otimes_{\Bbb Z_p}( A\hotimes \Bbb Z_p ). \tag 2.1.2$$ If $\CA$ is an E$_{\infty}$ algebra, then these are quasi-isomorphisms of E$_{\infty}$ algebras. 
 \endproclaim

\demo{Proof} (i) (2.1.1) is the composition of canonical quasi-isomorphisms $R\Gamma (X_{\h}, \CA )\otimes^L \Bbb Z/p^n = R\Gamma (X_{\h}, \CA \otimes^L \Bbb Z/p^n ) \buildrel{\sim}\over\leftarrow
R\Gamma (X_{\h}, A \otimes^L \Bbb Z/p^n )\buildrel{\sim}\over\leftarrow R\Gamma (X_{\et}, A \otimes^L \Bbb Z/p^n ) \buildrel{\sim}\over\leftarrow R\Gamma (X_{\et}, \Bbb Z/p^n )\otimes^L_{\Bbb Z/p^n} (A \otimes^L \Bbb Z/p^n )$, the first $\buildrel{\sim}\over\leftarrow$ comes from Remark (i), the second one comes since, by  Deligne's cohomological descent, the \'etale and h-cohomology with torsion coefficients coincide (see Remark in \cite{B} 3.4),  the third one is \cite{G} 3.3. 

(ii) By the condition on $F$, the complex $R\Gamma_{\!\et}(X,\Bbb Z_p ):= \text{holim}_n R\Gamma (X_{\et}, \Bbb Z/p^n )$ is $\Bbb Z_p$-perfect and 
$ R\Gamma_{\!\et}(X,\Bbb Z_p )\otimes^L_{\Bbb Z_p} \Bbb Z/p^n \iso  R\Gamma (X_{\et}, \Bbb Z/p^n )$. Thus (2.1.1) can be rewritten as 
$R\Gamma (X_{\h}, \CA )\otimes^L \Bbb Z/p^n \iso R\Gamma_{\!\et}(X,\Bbb Z_p )\otimes^L_{\Bbb Z_p} (A \otimes^L \Bbb Z/p^n )$. Applying holim$_n$, we get (2.1.2). \qed
\enddemo

\subhead{\rm 2.2}\endsubhead From now on $K$ is a $p$-adic field as in 1.16. Let $(V,\bar{V})$ be an ss-pair over $\bar{K}$, see \cite{B} 2.2(c).
As in \cite{B} 3.2, we view it as a log $W(k)$-scheme with underlying scheme $\bar{V}$. The final object $\Spec (\bar{K},O_{\bar{K}})$ of  $\CV ar^{\ss}_{\bar{K}}$  is  $(\bar{S},\bar{\CL})$ (see 1.17), so  $(V,\bar{V})$ is a log scheme over $(\bar{S},\bar{\CL})$.

As in 1.12, one has the absolute log crystalline complexes $R\Gamma_{\!\crys} (V,\bar{V} )_{(n)}$ and $R\Gamma_{\!\crys} (V,\bar{V} )=\holim_n \, R\Gamma_{\!\crys} (V,\bar{V} )_{(n)}$. By (1.17.3), 
 $R\Gamma_{\!\crys}(\Spec (\bar{K},O_{\bar{K}})) =  \text{A}_{\crys} $.

\enspace

According to the lemma in  \cite{B} 4.1, there is a finite extension $K$ of $K_0$, $K\subset \bar{K}$, and a log smooth integral map $f: (Z,\CM )\to (S,\CL )=(\Spec\, O_{K}, \CL_{K})$ with $f_1$ of Cartier type, together with an identification of $(V,\bar{V})/ (\bar{S},\bar{\CL})$ with the pullback of $(Z,\CM )/ (S,\CL )$
 by $\theta : (\bar{S},\bar{\CL})\to (S,\CL )$ (see 1.18). 
 By (i) of Theorem in 1.18, one has:

\proclaim{Proposition}  $R\Gamma_{\!\crys} (V,\bar{V} )$ is a perfect $\Acrys$-complex and $R\Gamma_{\!\crys} (V,\bar{V} )\otimes^L \Bbb Z/p^n
 \iso R\Gamma_{\!\crys} (V,\bar{V} )_{(n)} $. \qed
\endproclaim

\remark{Remark} Enlarging $K$, we find that any finite diagram of $(V,\bar{V})$'s comes from a diagram of $(Z,\CM )$'s  over common $(S,\CL )$. \endremark

\subhead{\rm 2.3}\endsubhead
Let $\CA_{\crys}$ be  h-sheafification of the presheaf 
 $(V,\bar{V}) \mapsto R\Gamma_{\!\crys} (V,\bar{V} )$ on  $\CV ar^{\ss}_{\bar{K}}$  (see \cite{B} 2.6).\footnote{To see  $\CA_{\crys}$ explicitly, one computes the complex of presheaves $(V,\bar{V})\mapsto  R\Gamma_{\!\crys} (V,\bar{V} )$ using Godement's resolution,  sheafifies it  for the h-topology on $\CV ar^{\ss}_{\bar{K}}$, and views the result as a complex of h-sheaves on $\CV ar_{\bar{K}}$. } This is an h-sheaf of E$_{\infty}$  $\Acrys$-algebras on $\CV ar_{\bar{K}}$ equipped with the Frobenius action $\varphi$. Since h-sheafification is exact,  $\CA_{\crys\, n}:= \CA_{\crys}\otimes^L \Bbb Z/p^n$ equals the h-sheafification of 
 the presheaf 
 $(V,\bar{V}) \mapsto R\Gamma_{\!\crys} (V,\bar{V} )_{(n)}$ by Proposition in 2.2.

\proclaim  {Theorem {\rm (the crystalline  Poincar\'e lemma)}}  $ \CA_{\crys}
$ is $p$-adically constant. \endproclaim

\demo{Proof} Since $(\text{A}_{\crys\, 1 })_{\CV ar_{\bar{K} \h}}\iso H^0\CA_{\crys\, 1}$,  we need to show that $H^{>0} \CA_{\crys\, 1}=0$. It suffices to show that every $(V,\bar{V})\in \CV ar^{\ss}_{\bar{K}}$ admits an h-covering $(V',\bar{V}')\to (V,\bar{V})$ such that the  map $H^{>0}R\Gamma_{\!\crys}(V,\bar{V})_{(1)}  \to H^{>0}R\Gamma_{\!\crys}(V',\bar{V}')_{(1)}$ vanishes. By the next lemma, any composition of $\dim V +1$  p-negligible h-coverings  (see \cite{B} 4.3) does the job, so we are done by the theorem in \cite{B} 4.3. 

Let $(V^{m+1} , \bar{V}^{m+1} )\buildrel{\psi_m}\over\lra \ldots \buildrel{\psi_1}\over\lra (V^1 ,\bar{V}^1)$ be  p-negligible  maps in $\CV ar^{\ss}_{\bar{K}}$.

\proclaim{ Lemma} 
The composition 
$\tau_{>0}R\Gamma_{\!\crys}(V^1 ,\bar{V}^1)_{(1)}  \to \tau_{>0}R\Gamma_{\!\crys}(V^{m+1} ,\bar{V}^{m+1})_{(1)}$ vanishes   if $m> \dim V^1$.
\endproclaim

\demo{Proof of Lemma} Choose $K$ and $f_i : (Z^i ,\CM^i )\to (S,\CL )=(\Spec\, O_{K},\CL_{K})$ for $(V^i ,\bar{V}^i )$ as in 2.2 so that $\psi_i$ come from 
morphisms $\psi_i : (Z^{i+1} ,\CM^{i+1} )\to (Z^{i} ,\CM^{i} )$ over $(S,\CL )$. Since $ R\Gamma (Z^i_1 , \Omega^a_{(Z^i,\CM^i )_1 /(S,\CL )_1})= R\Gamma (Z^i  ,\Omega^a_{(Z^i,\CM^i ) /(S,\CL )}) )\otimes^L_{O_{K}} (O_{K}/p )$  and $\psi_i $ is  $p$-negligible, we know that 
 the $\psi_i^*$ morphisms
 $ R\Gamma (Z^i_1 , \Omega^a_{(Z^i,\CM^i )_1/(S,\CL )_1})\to R\Gamma (Z^{i+1}_1 \! , \Omega^a_{(Z^{i+1}\! ,\CM^{i+1} )_1 /(S,\CL )_1})$, $a>0$, and
  $\tau_{>0} R\Gamma (Z^i_1 ,\CO_{Z^i_1})\to \tau_{>0} R\Gamma (Z^{i+1}_1 \! ,\CO_{Z^{i+1}_1})$ 
vanish as maps in the derived category of $O_{K}/p\,$-modules.\footnote{For $O_{\bar{K}}/p$ is faithfully   $O_{K}/p$-flat  and $R\Gamma (Z^i , \Omega^a_{(Z^i,\CM^i )/(S,\CL )})$ are  $O_{K}/p$-perfect complexes.} We are done now by (1.10.1), since the span of the conjugate filtration on the source is $[0,\dim V^1 ]$. 
\qed
\enddemo
\enddemo

\subhead{\rm 2.4}\endsubhead
For $X\in\CV ar_{\bar{K}}$ set 
$ R\Gamma_{\!\crys} (X ):=R\Gamma (X_{\h},   \CA_{\crys})$.
 This is an  E$_\infty$ $\text{A}_{\crys}$-algebra 
 equipped with the Frobenius action $\varphi$. The Galois group Gal$(\bar{K}/K_0)$ acts on $\CV ar_{\bar{K}}$, and it acts on $X\mapsto R\Gamma_{\!\crys} (X )$ by transport of structure. In particular, if $X$ is defined over an extension $K\subset \bar{K}$ of $K_0$,  $X=X_{K}\otimes_{K} {\bar{K}}$, then Gal$(\bar{K}/K)$ acts naturally on $R\Gamma_{\!\crys} (X )$. We write $H^n_{\crys} (X ):=H^n R\Gamma_{\!\crys} (X )$.

By the theorem in 2.3, the proposition in 2.1 provides
 canonical quasi-isomorphisms of E$_\infty$ A$_{\crys\, n}$-algebras $$R\Gamma_{\!\crys} (X ) \otimes^L \Bbb Z/p^n \iso R\Gamma(X_{\et},\Bbb Z/p^n )\otimes^L_{\Bbb Z/p^n} \text{A}_{\crys\, n} \tag 2.4.1$$ and, since  $A_{\crys}\iso A_{\crys}\hotimes\Bbb Z_p$,  a canonical quasi-isomorphism of E$_\infty$ A$_{\crys}$-algebras
 $$ R\Gamma_{\!\crys} (X ) \hotimes \Bbb Z_p \iso R\Gamma_{\!\et} (X,\Bbb Z_p )\otimes_{\Bbb Z_p} \Acrys . \tag 2.4.2$$ 
 
 \subhead{\rm 2.5}\endsubhead To control $R\Gamma_{\!\crys} (X)\otimes \Bbb Q$, we use its connection with the de Rham cohomology provided by the Hyodo-Kato theory:  
  
 For an ss-pair $(V,\bar{V})$ over $\bar{K}$, set
 $R\Gamma_{\!\text{H}\!\text{K}}(V,\bar{V}):=R\Gamma_{\!\text{H}\!\text{K}}((V,\bar{V})_1 )$, see Remark (i) in 1.18.   By loc.~cit., one has natural isomorphisms $\iota_{\crys} : R\Gamma_{\!\text{H}\!\text{K}}(V,\bar{V})^\tau_{\Bcrys^+}\iso R\Gamma_{\!\crys} (V,\bar{V} )\otimes\Bbb Q$, $\iota_{\dR}: R\Gamma_{\!\text{H}\!\text{K}}(V,\bar{V})^\tau_{\bar{K}} \iso R\Gamma_{\!\dR}(V_{\bar{K}},\bar{V}_{\bar{K}} ).$ Let $\CA_{\!\text{H}\!\text{K}}$ be h-sheafification of the presheaf
$(V,\bar{V})\mapsto R\Gamma_{\!\text{H}\!\text{K}}(V,\bar{V})$ on $\CV ar_{\bar{K}}^{\ss}$; this is an h-sheaf of E$_\infty$ $K_0^{\text{nr}}$-algebras on $\CV ar_{\bar{K}}$ 
equipped with $\varphi$-action and  locally nilpotent derivation $N$ such that $N\varphi =p\varphi N$.  For $X\in\CV ar_{\bar{K}}$ set 
$ R\Gamma_{\!\text{H}\!\text{K}} (X ):=R\Gamma (X_{\h},   \CA_{\text{H}\!\text{K}})$,  $H^n_{\text{H}\!\text{K}} (X ):= R\Gamma_{\!\text{H}\!\text{K}} (X )$. We get canonical quasi-isomorphisms\footnote{Here  $R\Gamma_{\!\dR}$ is Deligne's version of the de Rham cohomology, 
see  \cite{B} 3.4.}
$$\iota_{\crys} : R\Gamma_{\!\text{H}\!\text{K}}(X)^\tau_{\Bcrys^+}\!\iso R\Gamma_{\!\crys} (X)\otimes\Bbb Q,\quad  \iota_{\dR}: R\Gamma_{\!\text{H}\!\text{K}}(X)^\tau_{\bar{K}} \iso R\Gamma_{\!\dR}(X ) \tag 2.5.1$$ compatible with the Gal$(\bar{K}/K_0 )$-action; here ${}^\tau_{\Bcrys^+}\!$, ${}^\tau_{\bar{K}}$ are the crystalline and de Rham Fontaine-Hyodo-Kato twists (they commute with the passage to h-sheafification and $R\Gamma$).

\proclaim{ Proposition} (i) For any $(V,\bar{V})\in\CV ar^{\ss}_{\bar{K}}$ the canonical maps $R\Gamma_{\! \crys}(V,\bar{V})\otimes\Bbb Q\to R\Gamma_{\! \crys}(V)\otimes\Bbb Q$,
$R\Gamma_{\!\text{H}\!\text{K}}(V,\bar{V})\to R\Gamma_{\!\text{H}\!\text{K}}(V)$
are quasi-isomorphisms. \newline(ii) For every $X\in \CV ar_{\bar{K}}$ the cohomology groups $H^n_{\crys} (X)\otimes \Bbb Q$, $H^n_{\text{H}\!\text{K}}(X)$
 are free $\Bcrys^+$-modules, resp.~$K_0^{\text{nr}}$-vector spaces, of rank equal to $\dim H^n_{\et} (X,\Bbb Q_p )$. The same is true for the relative cohomology groups for a map of varieties.
 \endproclaim

\demo{Proof} (i) The map $R\Gamma_{\dR}(V,\bar{V})\to R\Gamma_{\dR}(V)$ is a quasi-isomorphism by usual mixed Hodge theory (see (i) of the proposition  in \cite{B} 3.4). Using $\iota_{\dR}$, we see that $R\Gamma_{\!\text{H}\!\text{K}}(V,\bar{V})\iso R\Gamma_{\!\text{H}\!\text{K}}(V)$; applying $\iota_{\crys}$, we get $R\Gamma_{\! \crys}(V,\bar{V})\otimes\Bbb Q\iso R\Gamma_{\! \crys}(V)\otimes\Bbb Q$.

(ii) One has $\dim_{\bar{K}} H^n_{\dR}(X)= \dim H^n_{\et} (X,\Bbb Q_p )$ (see (ii) of the proposition  in \cite{B} 3.4). Now use $\iota_{\dR}$ and $\iota_{\crys}$ as in (i).\footnote{  $ H^n_{\crys} (X)\otimes \Bbb Q$ is a free $\Bcrys^+$-module being a twist of  $H^n_{\text{H}\!\text{K}}(X)\otimes_{K_0^{\text{nr}}} \Bcrys^+$  by a trivial torsor.}
Ditto for the relative cohomology.  \qed
 \enddemo

\remark{Exercise} Assertion (ii) remains valid for any finite diagram of varieties.
 \endremark

\subhead{\rm 2.6}\endsubhead Assertion (i) of the above proposition  can be generalized as follows. Let $(Z,\CM )/(S,\CL )$ be as in Example at the end of 1.16 with $K$ a subfield of $\bar{K}$. 
Denote by $(Z\,\bar{},\CM\,\bar{} \,)$ the pullback of $(Z,\CM )$ by $(\bar{S},\bar{\CL})\to (S,\CL )$ (see 1.17, 1.18 for the notation), by $(\bar{X},\CM_{\bar{X}})$  the generic fiber of $(Z\,\bar{},\CM\,\bar{} \,)$, i.e., the pullback of $(Z,\CM )$ to $\bar{K}$, and by $X\subset \bar{X}$ the open subset of triviality of the log structure.
 Let us define
 canonical maps $$R\Gamma_{\!\crys}(Z\,\bar{},\CM\,\bar{}\, ) \to R\Gamma_{\!\crys}(X), \quad R\Gamma_{\!\text{H}\!\text{K}} (Z_1,\CM_1 )\to R\Gamma_{\!\text{H}\!\text{K}}(X). \tag 2.6.1$$
 Consider the arithmetic pair $(X,Z\,\bar{}\,)$ over $\bar{K}$, and   pick any h-hypercovering $(V_\cdot ,\bar{V}_\cdot )$
 of $(X,Z\,\bar{}\,)$
  by ss-pairs (see \cite{B} \S 2).
One has an evident map of log schemes\footnote{We view pairs as log schemes as in \cite{B} 3.2.}  $(X,Z\,\bar{}\,)\to 
(Z\,\bar{},\CM\,\bar{} \,)$, so   $(V_\cdot ,\bar{V}_\cdot )$ is a simplicial log scheme over $(Z\,\bar{},\CM\,\bar{} \,)$. Our maps are compositions $ R\Gamma_{\!\crys}(  Z\,\bar{},\CM\,\bar{}\,   ) \to R\Gamma_{\!\crys}(V_\cdot ,\bar{V}_\cdot )\to   R\Gamma_{\!\crys}(V_\cdot ) \buildrel{\sim}\over\leftarrow R\Gamma_{\!\crys}(X)$ and $  R\Gamma_{\!\text{H}\!\text{K}} (Z_1,\CM_1 )\to R\Gamma_{\!\text{H}\!\text{K}}(Z\,\bar{},\CM\,\bar{} \,)\to
R\Gamma_{\!\text{H}\!\text{K}}(V_\cdot ,\bar{V}_\cdot )\to R\Gamma_{\!\text{H}\!\text{K}}(V_\cdot  )\buildrel{\sim}\over\leftarrow R\Gamma_{\!\text{H}\!\text{K}}(X)$.

Suppose the 
 sheaf of groups $\CM^{\text{gr}}_{\bar{X}}/\CO^\times_{\bar{X}}$   has trivial torsion (e.g.~$\CM_{\bar{X}}$ is saturated).

\proclaim{Proposition}
The maps of  (2.6.1)   yield quasi-isomorphisms 
$R\Gamma_{\!\crys}(  Z\,\bar{},\CM\,\bar{}\,   )\otimes\Bbb Q $ $ \iso R\Gamma_{\!\crys}(X)\otimes\Bbb Q$, $R\Gamma_{\!\text{H}\!\text{K}} (Z_1,\CM_1 )\otimes_{K_0}\!K_0^{\text{nr}}\iso R\Gamma_{\!\text{H}\!\text{K}}(X)$.
\endproclaim 

\demo{Proof} 
Our maps are compatible with identifications $\iota_{\crys}$ of (1.18.2) and (2.5.1), so it suffices to check the claim for $R\Gamma_{\!\text{H}\!\text{K}}$.  Due to
isomorphisms $\iota_{\dR}$ of (1.16.3) and (2.5.1), it is enough to show that the restriction  $ R\Gamma (\bar{X},\Omega^\cdot_{(\bar{X}, \CM_{\bar{X}} )/\bar{K}})\to R\Gamma (X,\Omega^\cdot_{X/\bar{K}})$ $=R\Gamma_{\!\dR}(X)$ is a quasi-isomorphism, which follows from Ogus' theorem in 1.19.\qed \enddemo

\head 3. The  Fontaine-Jannsen conjecture. \endhead

\subhead{\rm 3.1}\endsubhead  For $X\in\CV ar_{\bar{K}}$ we  define   {\it the crystalline period map}
$$\rho_{\crys} : R\Gamma_{\!\crys} (X ) \to R\Gamma_{\!\et}(X,\Bbb Z_p )\otimes_{\Bbb Z_p} \text{A}_{\crys} \tag 3.1.1$$  as the composition of the evident map $
R\Gamma_{\!\crys} (X ) \to R\Gamma_{\!\crys} (X )\hotimes\Bbb Z_p$ with identification (2.4.2). Composing $\rho_{\crys }\otimes\Bbb Q$ with isomorphism $\iota_{\crys}$ from (2.5.1), we get
 $$\rho_{\text{H}\!\text{K}  } : R\Gamma_{\!\text{H}\!\text{K}}(X)^\tau_{\Bcrys^+}\! \to R\Gamma_{\!\et}(X,\Bbb Q_p )\otimes_{\Bbb Q_p} \text{B}^+_{\crys}. \tag 3.1.2$$ These are morphisms of  E$_\infty$ $\text{A}_{\crys}$- and $\Bcrys^+$-algebras equipped  with the Frobenius action $\varphi$ (it acts on the target via the second factor). The Galois group Gal$(\bar{K}/K_0 )$ acts on $\CV ar_{\bar{K}}$ and  all the functors; $\rho_{\crys}$, $\rho_{\text{H}\!\text{K}}$ are compatible with this action. Thus if $X$ is defined over an extension  $K\subset \bar{K}$ of $K_0$, $X=X_{K}\otimes_{K}\bar{K}$, then Gal$(\bar{K}/K)$ acts on both terms of (3.1.1), (3.1.2) and $\rho_{\crys}$, $\rho_{\text{H}\!\text{K}}$ commute with the Galois action.

\subhead{\rm 3.2}\endsubhead  We use the notation from 1.17. Below 
$\Bbb Q_p /\Bbb Z_p (1)$ is the subgroup of $p^\infty$-roots of 1 in $ O_{\bar{K}}^\times$, so 
  $\Bbb Q_p (1)$ is projective limit of the system $\ldots \buildrel{p}\over\to \Bbb Q_p /\Bbb Z_p (1) \buildrel{p}\over\to \Bbb Q_p /\Bbb Z_p (1)$. Thus one has an embedding
 $\Bbb Q_p (1)\hra \CL_{\varphi}$,
 $(\varepsilon^{(n)}) \mapsto (\varepsilon^{(n)} \,\text{mod}\, pO_{\bar{K}})\in\CL_{\varphi}$, hence a canonical embedding $\Bbb Q_p (1)\hra 
\Acrys^\times$. The image of $\Bbb Z_p (1)$ lies in $(1+ J_{\crys})^\times$; applying $\log :(1+ J_{\crys})^\times  \to J_{\crys}$, we get  $l: \Bbb Z_p (1) \hra J_{\crys}\subset \Acrys$.  As  in  \cite{F1} 2.3.4, set $\Bcrys  :=\Acrys [l(t)^{-1}]$, where $t$ is a  generator of $\Bbb Z_p (1) $.  Inverting $l(t)$ implies inverting $p$,\footnote{For the pd structure on $J_{\crys}$ provides $p^{-1}l(t)^p \in J_{\crys}$.} so $\Bcrys \supset \Bcrys^+ $.

\proclaim{Theorem} The $\Bcrys$-linear extensions  of $\rho_{\crys}$ and $\rho_{\text{H}\!\text{K}}$ are quasi-isomorphisms: for any $X\in\CV ar_{\bar{K}}$ one has $\rho_{\crys} : R\Gamma_{\!\crys} (X )\otimes_{\Acrys}\! \Bcrys \iso R\Gamma_{\!\et}(X,\Bbb Q_p )\otimes \Bcrys $, $\rho_{\text{H}\!\text{K}  } : R\Gamma_{\!\text{H}\!\text{K}}(X)^\tau_{\Bcrys}\! \iso R\Gamma_{\!\et}(X,\Bbb Q_p )\otimes  \text{B}_{\crys}$.
\endproclaim

\demo{Proof} It is very similar to that of the $\rho_{\dR}$ counterpart in \cite{B} 3.6: there is a  calculation for circle (Lemma  below), the rest comes by a general functoriality argument.

(a) {\it The case of $X=\Bbb G_m=\Bbb P^1 \smallsetminus\{0,\infty\}$}: Let $t$ be the standard parameter on $\Bbb G_m$, and $\bar{\Bbb G}_m$ be $\Bbb P^1$ viewed as a $\Bbb G_m$-equivariant compactification of $\Bbb G_m$. So $(\Bbb G_{m \bar{K}}, \bar{\Bbb G}_{m \bar{S}})$ is an ss-pair over $\bar{K}$; denote the corresponding log $\bar{S}$-scheme by $Y_{\bar{S}}$. Its log structure is generated by $\bar{\CL}$ and $t$. Let $Y_{\crys}$ be 
$\bar{\Bbb G}_{m E_{\crys}}$ equipped with the log structure on $\bar{\Bbb G}_{m E_{\crys}}$  generated by $t$ and $\CL_{\crys}$. 

By (ii) of the proposition in 2.5, it suffices to consider the group $H^1$.  
Consider the canonical map $R\Gamma_{\!\crys}(Y_{ \bar{S}})\to R\Gamma_{\!\crys}(\Bbb G_{m \bar{K}})$. Since $Y_{\crys}$ is a pd-smooth object of $(Y_{\bar{S}_1} /(E_{\crys},\CL_{\crys}))_{\crys}$, the $R\Gamma$ of its de Rham complex equals $
R\Gamma_{\!\crys}(Y_{ \bar{S}})$ (see (1.8.1), (1.18.1)), so $H^1_{\crys}(Y_{ \bar{S}})$ is a free $\Acrys$-module generated by $d\log t$. By (i) of the proposition in 2.5, $H^1_{\crys}(\Bbb G_{m \bar{K}})\otimes\Bbb Q$ is a free $\Bcrys^+$-module generated by the image of $d\log t$. Let $\kappa$ be the canonical generator  of $H^1_{\et} (\Bbb G_{m \bar{K}}, \Bbb Z_p (1))= H^1_{\et} (\Bbb G_{m \bar{K}}, \Bbb Z_p ) \otimes\Bbb Z_p (1)$; applying $l: \Bbb Z_p (1) \hra  \Acrys$, we get $l(\kappa )\in  H^1_{\et} (\Bbb G_{m \bar{K}}, \Bbb Z_p )\otimes \Acrys $. It remains to prove:
 
 \proclaim{Lemma} One has
 $\rho_{\crys}(d\log t )=l(\kappa ).$
\endproclaim

\demo{Proof of Lemma}
We do mod $p^n$ computation.
Let $\Bbb G\tilde{_m}$ be a copy of $\Bbb G_m$ with parameter $\tilde{t}$, and $\pi: \Bbb G\tilde{_m}\to \Bbb G_m$ be the map $\pi^* (t)=\tilde{t}^{p^n}$. Then $\Bbb G\tilde{_m}_{\bar{K}}/\Bbb G_{m \bar{K}}$ is a $ \Bbb Z/p^n (1)$-torsor, and
 $\kappa_n \in H^1 (\Bbb G_{m \bar{K}\et}, \Bbb Z/p^n (1) )$ is its class. The corresponding  \v Cech hypercovering   is the twist of $ \Bbb G\tilde{_m}_{\bar{K}}$ by the universal $\Bbb Z/p^n (1)$-torsor $\ft$ over the simplicial classifying space $B_\cdot $ of $\Bbb Z/p^n (1)$, and $\kappa_n$ comes from an evident 1-cocycle on $B_\cdot$. 
Now $\pi$ extends to an h-covering $ ( \Bbb G\tilde{_m}_{\bar{K}}, \bar{ \Bbb G}\tilde{_m}_{\bar{S}}   )\to ( \Bbb G_{m\bar{K}} ,  \bar{\Bbb G}_{m\bar{S}} )$ of ss-pairs; its $\ft$-twist $\pi_\cdot : 
  ( \Bbb G\tilde{_m}_{\bar{K}}, \bar{ \Bbb G}\tilde{_m}_{\bar{S}}   )_\cdot 
\to ( \Bbb G_{m\bar{K}} ,  \bar{\Bbb G}_{m\bar{S}} )$ is an h-hypercovering in
$\CV ar^{\ss}_{\bar{K}}$. Changing the notation, we have $\pi : Y\tilde{_{\bar{S}}}\to Y_{\bar{S}}$ and $\pi_\cdot : Y\tilde{_{\bar{S}}}_{\cdot}\to Y_{\bar{S}}$. We want to check that 
 $\pi^*_{\cdot \crys}: R\Gamma_{\!\crys}(Y_{\bar{S}})_{(n)} \to R\Gamma_{\!\crys}(Y\tilde{_{\bar{S}}}_{\cdot}
 )_{(n)}$ sends $d\log t$ to $l(\kappa_n )$. To do this, we extend $\pi_{1\cdot }$ to a map of simplicial pd-smooth thickenings $P\tilde{_\cdot}\to P_\cdot$ over $E_{\crys\, n}$.

Set $P:=  Y_{\crys \, n}$,  $\CG := \Bbb G^\natural_{m E_{\crys \, n}}$(:= the pd-envelope of $\Bbb G_m$ at 1, see 1.2). Then $\CG$ acts on $P$ through $\CG \to \Bbb G_m$. Our $P_\cdot$ is the twist of $P$ by the universal $\CG$-torsor over the
simplicial classifying space $B_{\CG \cdot}$.\footnote{$P_\cdot$ is equal to the simplicial object $P_*$ of $(Y_{\bar{S}_1}/(E_{\crys},\CL_{\crys})_n )_{\crys}$ from 1.6.} Consider now the $\,\tilde{}\,$-copies $P\,\tilde{}$ and $\CG\,\tilde{}$. Together with the $\CG\,\tilde{}$-action,
$P\,\tilde{}$  carries an action of  $p^{-n}\Bbb Z_p$ via the composition $p^{-n}\Bbb Z_p (1)\subset \Bbb Q_p (1) \hra \Acrys^\times \twoheadrightarrow \text{A}_{\crys\, n}^\times = \Bbb G\tilde{_m}(E_{\crys \, n})$, where $\hra$ is the canonical embedding. Its restriction to 
$\Bbb Z_p (1)$ lands in $(1+J_{\crys\, n})^\times $, i.e., we have a homomorphism $\alpha : \Bbb Z_p (1)\to 
 \CG\,\tilde{}\, (E_{\crys \, n})$. Both actions combine into an action  of the group pd-scheme $\CG^+$ which is an extension of $(p^{-n}\Bbb Z_p /\Bbb Z_p )(1)_{E_{\crys \, n}}$ by   $\CG\,\tilde{}\,$ defined as the pushout of   $0\to \Bbb Z_p (1)\to p^{-n}\Bbb Z_p (1)\to (p^{-n}\Bbb Z_p /\Bbb Z_p )(1)\to 0$ by $\alpha$. Our   $P\tilde{_\cdot}$ is the twist of $P\,\tilde{}\,$ by the universal $\CG^+$-torsor over the simplicial classifyling space $B_{\CG^+ \cdot}$. Extension $\CG^+$ splits over $\bar{S}_1$ since $\alpha$ vanishes at
 $\bar{S}_1 \subset E_{\crys\, n}$,  so we get an exact embedding of simplicial log $E_{\crys \, n}$-schemes $Y\tilde{_{\bar{S}_1}}_{\cdot }\hra P\tilde{_\cdot}$. The pd structure on $\CG\,\tilde{}\,$ provides a pd structure on its ideal. Finally, the projection $\pi
 :   P\,\tilde{}\to P$ and an evident ``multiplication by $p^n$" morphism $\pi^+ : \CG^+ \to \CG$ yield a map $\pi_{P\cdot} : 
P\tilde{_\cdot}\to P_\cdot$ of pd-thickenings  that extends $\pi_{1\cdot}$.

The pd-thickenings $P_i$ and $P\tilde{_i}$ are pd-smooth over $E_{\crys\, n}$, so, by (1.8.1), (1.18.1), the map  $\pi^*_{\cdot \crys}: R\Gamma_{\!\crys}(Y_{\bar{S}})_{(n)} \to R\Gamma_{\!\crys}(Y\tilde{_{\bar{S}}}_{\cdot} )_{(n)}$
coincides with the pullback map  $\pi^{*}_{P\cdot} : R\Gamma (P_\cdot ,\Omega^\cdot_{P_\cdot })\to R\Gamma (P\,\tilde{_\cdot}, \Omega^\cdot_{P\,\tilde{_\cdot}})  $  between the total de Rham complexes; here $\Omega^\cdot_? :=\Omega^\cdot_{?/(E_{\crys},\CL_{\crys})_n}$.
Now $d\log t \in \Gamma(P_0 , \Omega^1_{P_0})$ extends to a total 1-cocycle in 
$\Gamma (P_\cdot ,\Omega^\cdot_{P_\cdot })$ by adding the component $\log \chi \in\Gamma (P_1 , \CO_{P_1} )$ which comes from the evident $\CG$-valued 1-cocycle $\chi$ on $B_{\CG \cdot}$. One has $\pi^{*}_{P0} (d\log t)=p^n d\log \tilde{t}=0$, and 
$\pi^{*}_{P1} (\log \chi )$ comes from the 1-cocycle $\log (\chi \pi^+ )$ on 
$B_{\CG \cdot}$, which is 
 $l(\kappa_n )$, q.e.d. \qed \enddemo

(b) {\it Compatibility of $\rho_{\crys}$ with the Gysin maps for codimension 1 closed embeddings of smooth varieties}: Let $i: Y \hra X$ be such an embedding.
For any cohomology theory $R\Gamma_{\!?}$  we deal with,  consider the cohomology with supports  $R\Gamma_{\!?\, Y}(X):= \CC one ( R\Gamma_{\!?}(X)\to
R\Gamma_{\!?}(X\smallsetminus Y))[-1]$. Recall a definition of the Gysin isomorphism
$i_* :R\Gamma_{\! ?}(Y) \iso R\Gamma_{\!?\, Y}(X )[2].$ (Here the Tate twist $(1)$ in the target is canceled due to specifics of the cohomology theories we deal with, see below.)

Let $\CL $ be the normal line bundle, $i_0 : Y\hra \CL$ its zero section. 
There is a canonical identification $\eta :  R\Gamma_{\!?\, Y}(\CL )\iso R\Gamma_{\!?\, Y}(X )$ defined using the deformation to normal cone construction. Namely, we have the diagram  $$\spreadmatrixlines{1\jot}
\matrix
\CL
 &\!\!\!\hra&
 X\tilde{_{\Bbb A^1}} & \hookleftarrow &X   \\  
\,\,\,\,\uparrow i_0 &&\uparrow&& \,\uparrow i \\
Y &\! \!\! \hra & Y_{\Bbb A^1} &\hookleftarrow & Y. 
    \endmatrix
\tag 3.2.1$$
Here $Y_{\Bbb A^1}= Y\times \Bbb A^1$, $X\tilde{_{\Bbb A^1}}$ is $X\times \Bbb A^1$ with $Y\times\{0\}$ blown up,   the bottom embeddings are $y\mapsto (y,0), (y,1)$. The  arrows $R\Gamma_{\!?\, Y}(\CL )\leftarrow R\Gamma_{\!?\, Y_{\Bbb A^1}}(X\tilde{_{\Bbb A^1}})\to R\Gamma_{\!?\, Y}(X )$ are quasi-isomorphisms (this is standard for $R\Gamma_{\!\et}$, $R\Gamma_{\!\dR}$;\footnote{Use purity to identify the arrows with $R\Gamma_? (Y)\leftarrow R\Gamma_? (Y_{\Bbb A^1})\to R\Gamma_? (Y)$.  } the assertion for $R\Gamma_{\!\text{H}\!\text{K}}$ and $R\Gamma_{\!\crys}\otimes\Bbb Q$ is deduced from that for $R\Gamma_{\!\dR}$ using (2.5.1)). Their composition is $\eta$.

One has $i_* := \eta i_{0*}$, so it suffices to define $i_{0*} $. The projection $\CL \to Y$ makes $R\Gamma_{\!?\, Y}(\CL )$ an 
$R\Gamma_{\!?}(Y)$-module, and $ i_{0*} $ is a morphism of $R\Gamma_{\!?}(Y)$-modules. Thus to define $i_{0*}$, we need to specify the orientation class $i_{0*}(1)\in H^2_{? Y}(\CL )\buildrel{\sim}\over\leftarrow
\tau_{\le 0} (R\Gamma_{\!?\, Y}(\CL )[2])$ (here $\buildrel{\sim}\over\leftarrow$ holds due to (ii) of the proposition in 2.5). If $Y$ is connected, then $H^2_{? Y}(\CL )$ is a free module of rank 1 (see loc.~cit.), so, localizing $Y$, we can assume that $\CL$ is trivialized. Then $ H^2_{? Y}(\CL )\buildrel{\sim}\over\leftarrow H^1_? (\Bbb G_{m\bar{K}})$, and
we define  $i_{0*}(1)$ to be $d\log t$  for $R\Gamma_{\!\dR}$, $R\Gamma_{\!\text{H}\!\text{K}}$, $R\Gamma_{\!\crys}\otimes\Bbb Q$, and $l(\kappa )$ for $R\Gamma_{\!\et}(\cdot ,\Bbb Q_p )\otimes\Bcrys$ (see (a) for the notation). We have defined $i_{0*}$, hence $i_*$. It is an isomorphism (standard for $R\Gamma_{\!\et}$, $R\Gamma_{\!\dR}$; the assertion for $R\Gamma_{\!\text{H}\!\text{K}}$ and $R\Gamma_{\!\crys}\otimes\Bbb Q$ is deduced from that for $R\Gamma_{\!\dR}$ using (2.5.1)). Evidently $i_*$ commutes with the maps of (2.5.1) and, by 
lemma in (a), with $\rho_{\crys}$.

(c) {\it The case of smooth projective $X$}: We can assume that $X$ is connected, $\dim X=d$. Then $H^{2d}_? (X)$ is a free module of rank 1 (see (ii) of the proposition in 2.5), and the Poincar\'e duality pairing $H^i_? (X)\times 
H^{2d -i }_? (X)\to H^{2d}_? (X)$ is nondegenerate (standard for $R\Gamma_{\!\et}$, $R\Gamma_{\!\dR}$; the assertion for $R\Gamma_{\!\text{H}\!\text{K}}$ and $R\Gamma_{\!\crys}\otimes\Bbb Q$ is deduced from that for $R\Gamma_{\!\dR}$ using (2.5.1)). Since $\rho_{\crys}$ is a morphism of algebras, it is compatible with the Poincar\'e duality. This, together with (ii) of the proposition in 2.5, implies that  $\rho_{\crys}: H^{i}_{\crys}(X)\otimes  \Bcrys \to H^{i}_{\et}(X,\Bbb Q_p )\otimes  \Bcrys$ is an isomorphism for every $i$ if this is true for $i=2d$. To check the latter assertion, consider the class $c_? \in H^2_? (X)$  of hyperplane section (defined using the corresponding Gysin map).  Then $c^d_{\et}$ is a base in $H^{2d}_{\et}(X,\Bbb Q_p )\otimes  \Bcrys$,
$c^d_{\crys}$ is a base in $H^{2d}_{\crys}(X)\otimes \Bbb Q$ (by (2.5.1), since $c_{\dR}^d$ is a base in $H^{2d}_{\dR}(X)$). Now  $\rho_{\crys}(c_{\crys})=c_{\et}$ by (b), and we are done.

(d) {\it The case when $X$ is the complement to a strict normal crossings divisor in a smooth projective variety}:  Checked exactly as in \cite{B} 3.6,
by induction by the number of the components of the divisor using (c) and (b).

(e) {\it The case of arbitrary $X$}: Checked exactly as in \cite{B} 3.6, using an h-hypercovering of $X$ by varieties as in (d).  
 \qed
\enddemo

\subhead{\rm 3.3}\endsubhead  The   theorem in 3.2 implies the Fontaine-Jannsen
conjecture. To see this, we pull 
 $\rho_{\text{H}\!\text{K}}$ back to  the Fontaine-Hyodo-Kato  $\Bbb G_a$-torsor\footnote{ Recall that $\Bst :=\Bst^+ \otimes_{\Bcrys^+}\Bcrys$.} $\Spec\, \Bst /\Spec\, \Bcrys$ to trivialize the  twist, cf.~(1.18.5). We get a canonical quasi-isomorphism    of $\Bst$-complexes  $$\rho_{\text{H}\!\text{K}} :  R\Gamma_{\! \text{H}\!\text{K}}(X)  
 \otimes_{K_0^{\text{nr}}}\! \Bst
  \iso R\Gamma_{\!\et}(X,\Bbb Q_p )\otimes \Bst \tag 3.3.1$$ 
 compatible with the  $(\varphi ,N)$-action and with the Gal$(\bar{K}/K)$-action on $\CV ar_{\bar{K}}$.
 This is the identification
 asked for in \cite{F2} \S 6.

 Conjectures C$_{\text{pst}}$, C$_{\text{st}}$, and C$_{\crys}$ (see \cite{F2} 6.2.1,  6.2.7, 6.1.4) come as follows:  
 \newline
  (i) Suppose $X$ is defined over a finite extension $K$ of $K_0$, $K\subset \bar{K}$, so we have  $X_{K}$ over $K$ and an identification $X=X_{K}\otimes_{K}\bar{K}$. We get the Gal$(\bar{K}/K)$-action on $R\Gamma_{\! \text{H}\!\text{K}}(X)$. By  \cite{F1} 4.2.4,   $H^n_{\text{H}\!\text{K}}(X)$ coincides with the subspace of those elements in $H^n_{\text{H}\!\text{K}}(X)\otimes_{K_0^{\text{nr}}}\! \Bst
 $ whose stabilizers in Gal$(\bar{K}/K)$ are open, hence, via (3.3.1), with the similar subspace of $H^n_{\et}(X,\Bbb Q_p )\otimes \Bst$.
This is conjecture  C$_{\text{pst}}$.
 \newline
 (ii) Assume we are in the situation of 2.6. Then $R\Gamma_{\! \text{H}\!\text{K}} (Z_1 ,\CM_1 )\otimes_{K_0}K_0^{\text{nr}} \iso R\Gamma_{\! \text{H}\!\text{K}}(X)$  by the proposition in loc.~cit., so Gal$(\bar{K}/K)$ acts trivially on $R\Gamma_{\! \text{H}\!\text{K}}( Z_1 ,\CM_1  )$. We get $\rho_{  \text{H}\!\text{K} }:  R\Gamma_{\! \text{H}\!\text{K}}(Z_1 ,\CM_1 )\otimes_{K_0}\! \Bst  \iso R\Gamma_{\!\et}(X,\Bbb Q_p )\otimes_{\Bbb Q_p}\Bst $, which is  conjecture C$_{\text{st}}$.  \newline
 (iii) Assume we are in the situation of 2.6, and suppose that  $f_0 :(Z_{01} ,\CM_{01} )\to (S_{01}, \CL_{01} )$ from Remark (i) in 1.16 can be realized as the pullback of a log scheme $(Z_{01} ,\CM_{00} )$ over $\Spec\, k$ (with trivial log structure) by the structure map $(S_{01}, \CL_{01} )$ $\to \Spec\, k$. E.g.~this happens if $X$ is smooth proper and  $X_{K}$ has smooth model $Z$: then  the log structure $\CM_{00}$  is trivial. By the base change, $Rf_{0\,\crys *} (\CO_{Z_{01}/\Bbb Z_p})$ is the constant crystal 
with fiber $R\Gamma_{\!\crys}(Z_{01},\CM_{00} )$. Thus, by Remark (i) in 1.16, 
 $ R\Gamma_{\! \text{H}\!\text{K}}(Z_1,\CM_1 ) \iso R\Gamma_{\!\crys}(Z_{01},\CM_{00} )
 \otimes\Bbb Q $ and  $N$ acts trivially on the Hyodo-Kato cohomology. Therefore quasi-isomorphism $\rho_{\text{H}\!\text{K}}$ from the theorem in 3.2 can be rewritten as $R\Gamma_{\!\crys} (Z_{01} ,\CM_{00} )\otimes_{W(k)} \Bcrys \iso R\Gamma_{\!\et}(X,\Bbb Q_p )\otimes \Bcrys$, which is  C$_{\crys}$.

\subhead{\rm 3.4}\endsubhead  In the rest of the section we show that the crystalline period map  is compatible with its derived de Rham cousin
$\rho_{\dR}$ from \cite{B}:

\proclaim{Theorem-construction} For any $X\in\CV ar_{\bar{K}}$ there is a canonical isomorphism $$R\Gamma_{\!\dR}(X)\otimes_{\bar{K}}\BdR^+ \iso R\Gamma_{\!\crys} (X )\otimes_{\Acrys}\! \BdR^+ , \tag 3.4.1$$
compatible with the Galois action  such that $\rho_{\dR}$  from \cite{B} 3.5.4     is
the composition of (3.4.1) with $\rho_{\crys}\otimes \BdR^+ :R\Gamma_{\!\crys} (X )\otimes_{\Acrys}\! \BdR^+ \to R\Gamma_{\!\et}(X,\Bbb Q_p )\otimes \BdR^+$.  \endproclaim

\demo{Proof} Let us fix  notation. As in \cite{B} 3.3, for any  ss-pair    $(V,\bar{V})$ over $\bar{K}$ we denote by $R\Gamma_{\!\dR}^\natural (V,\bar{V})$ its absolute derived de Rham complex
$R\Gamma (\bar{V}, L\Omega^\cdot\!\hat{_{(V,}}{}_{\bar{V})/W(k)})$. Set $R\Gamma_{\!\dR}^\natural (V,\bar{V})_n := R\Gamma_{\!\dR}^\natural (V,\bar{V})
 \otimes^L \Bbb Z/p^n  \iso R\Gamma (\bar{V}, L\Omega^\cdot\!\hat{_{(V,}}{}_{\bar{V})_n/W_n (k)})$, where $\iso$ is the base change identification, $R\Gamma_{\!\dR}^\natural (V,\bar{V})\hotimes\Bbb Z_p $ $ := \text{holim}_n R\Gamma_{\!\dR}^\natural (V,\bar{V})_n$,  
$R\Gamma_{\!\dR}^\natural (V,\bar{V})\hotimes\Bbb Q_p $ $:=
(R\Gamma_{\!\dR}^\natural (V,\bar{V})\hotimes\Bbb Z_p )\otimes\Bbb Q$. These are $F$-filtered E$_\infty$ algebras, where  ``$F$-filtered" means that we view them as mere projective systems of quotients modulo the terms of Hodge filtration $F^m$. Below we consider the homotopy $F$-completions of these complexes, which we denote by $\lim_F$. 
So we have $\lim_F R\Gamma_{\!\dR}^\natural (V,\bar{V})_n = \holim_m 
(R\Gamma_{\!\dR}^\natural (V,\bar{V})_n /F^m )$, $\lim_F R\Gamma_{\!\dR}^\natural (V,\bar{V})\hotimes\Bbb Z_p $ $:= 
\holim_m (R\Gamma_{\!\dR}^\natural (V,\bar{V})\hotimes\Bbb Z_p /F^m ) $ $=
\text{holim}_{m,n} (R\Gamma_{\!\dR}^\natural (V,\bar{V})_n  /F^m )$, etc.

Recall that 
$\AdR =R\Gamma_{\!\dR}^\natural (\Spec\, (\bar{K},O_{\bar{K}}))$ (see the lemma in \cite{B} 3.2). The corresponding $F$-filtered algebras $\text{A}_{\dR\,n}$, $\AdR \hotimes\Bbb Z_p$, $\AdR \hotimes\Bbb Q_p$ are acyclic in degrees $\neq 0$ and their projections $\cdot /F^{m+1} \to \cdot /F^m$ are surjective (see \cite{B} 1.4, 1.5). Thus
 $\text{A}^\ds_{\dR\,n}:= \lim_F \text{A}_{\dR\,n}$ equals $ \limleft_m H^0 (\text{A}_{\dR\,n}/F^m )$, $\AdR^\ds :=\lim_F ( \AdR \hotimes\Bbb Z_p )$ equals $  \limleft_m H^0 (\AdR \hotimes\Bbb Z_p /F^m )$,\footnote{By Example in 3.6 below, $\text{A}^\ds_{\dR}$ is the $J^{[m ]}$-topology
 completion of $\Acrys$.  }
   and  $\lim_F (\AdR \hotimes\Bbb Q_p )= \limleft H^0 (\AdR \hotimes\Bbb Q_p /F^m ) = \BdR^+ $.\footnote{So  $\text{A}^\ds_{\dR}\otimes \Bbb Q$ is dense in $\BdR^+$, but not equal to it.} Notice that $\text{A}^\ds_{\dR\,n}=\AdR^\ds \otimes\Bbb Z/p^n = \AdR^\ds \otimes^L \Bbb Z/p^n$.

 For any $(V,\bar{V})$, the complex $R\Gamma_{\!\dR}^\natural (V,\bar{V})$  is an $F$-filtered E$_\infty$ filtered $\AdR$-algebra, so  $\lim_F R\Gamma_{\!\dR}^\natural (V,\bar{V})_n$ is an 
$\text{A}^\ds_{\dR\,n}$-algebra, $\lim_F (R\Gamma_{\!\dR}^\natural (V,\bar{V})\hotimes\Bbb Q_p )$ is a $\BdR^+$-algebra,
etc. In 3.5--3.7 below, we will construct the next natural quasi-isomorphisms:

\enspace

 (a) In 3.5, we define an $F$-filtered quasi-isomorphism $\gamma : R\Gamma_{\!\dR} (V) \otimes_{\bar{K}} (\AdR \hotimes\Bbb Q_p )$ $\iso R\Gamma^\natural_{\! \dR} (V,\bar{V})\hotimes\Bbb Q_p$. The filtration $F$ on $R\Gamma_{\!\dR} (V)$ is the Hodge-Deligne filtration (see \cite{B} 3.4); it is finite and $\gr^\cdot_F R\Gamma_{\!\dR} (V)$ is a perfect $\bar{K}$-complex. Thus we have $$\lim_F \gamma : R\Gamma_{\!\dR} (V) \otimes_{\bar{K}} \BdR^+ \iso \lim_F (R\Gamma^\natural_{\! \dR} (V,\bar{V})\hotimes\Bbb Q_p ). \tag 3.4.2$$

(b) In 3.6, using the Illusie-Olsson comparison (see 1.9), we construct a natural compatible system of ring homomorphisms $\text{A}_{\crys\, n} \to\text{A}_{\dR\, n}^\ds$ and 
quasi-isomorphisms $\kappa_n :  R\Gamma_{\!\crys}(V,\bar{V})_{(n)} \otimes^L_{\text{A}_{\crys\, n}} \text{A}_{\dR\, n}^\ds \iso \lim_F 
R\Gamma_{\!\dR}^\natural (V,\bar{V})_n$. The proposition in 2.2 yields $
R\Gamma_{\!\crys}(V,\bar{V}) \otimes^L_{\text{A}_{\crys}} \text{A}_{\dR}^\ds \iso
\holim_n \, R\Gamma_{\!\crys}(V,\bar{V})_{(n)} \otimes^L_{\text{A}_{\crys\, n}} \text{A}_{\dR\, n}^\ds $, so  we get
$$\kappa :=\holim_n \kappa_n : R\Gamma_{\!\crys}(V,\bar{V}) \otimes^L_{\text{A}_{\crys}} \text{A}_{\dR}^\ds \iso \lim_F (
R\Gamma_{\!\dR}^\natural (V,\bar{V})\hotimes\Bbb Z_p ).\tag 3.4.3$$

(c) Consider the evident $\AdR \hotimes\Bbb Z_p$-linear map $ R\Gamma_{\! \dR}^\natural  (V,\bar{V})\hotimes\Bbb Z_p \to  R\Gamma_{\! \dR}^\natural (V,\bar{V})\hotimes\Bbb Q_p$.
 In 3.7 we prove that the $\BdR^+$-linear extension of its $\lim_F$ is a quasi-isomorphism: $$(\lim_F (
R\Gamma_{\! \dR}^\natural (V,\bar{V})\hotimes\Bbb Z_p ))\otimes^L_{\AdR^\ds} \BdR^+ \iso \lim_F (R\Gamma_{\! \dR}^\natural (V,\bar{V})\hotimes\Bbb Q_p ). \tag 3.4.4$$

Assuming (a)--(c), let us deduce the theorem. Consider $R\Gamma_{\!\dR} (V) \otimes_{\bar{K}} \BdR^+ \iso R\Gamma_{\!\crys}(V,\bar{V}) \otimes^L_{\text{A}_{\crys}}\BdR^+$ defined as the composition  ((3.4.3)$\otimes \BdR^+ )^{-1}$(3.4.4)${}^{-1}$(3.4.2). This is a quasi-isomorphism of presheaves on $\CV ar^{\ss}_{\bar{K}}$. Its h-sheafification is a quasi-isomorphism of h-sheaves $\CA_{\dR}\otimes_{\bar{K}}\BdR^+ \iso \CA_{\crys}\otimes^L_{\Acrys}\BdR^+$ on $\CV ar_{\bar{K}}$ (we use the notation of \cite{B} 3.4); applying  $R\Gamma (X_{\h}, \cdot )$, we get (3.4.1). 
  The construction is natural, so it commutes with the Galois action. The final property is evident from the constructions of $\rho_{\crys}$ and $\rho_{\dR}$. We are done. \qed
\enddemo

\subhead{\rm 3.5}\endsubhead We consider (a) above.
Recall that we have a filtered quasi-isomorphism $R\Gamma_{\!\dR}^\natural (V,\bar{V})\otimes\Bbb Q \iso R\Gamma_{\!\dR}(V)$, so the evident map $R\Gamma_{\!\dR}^\natural (V,\bar{V})\to R\Gamma_{\!\dR}^\natural (V,\bar{V})\hotimes\Bbb Z_p$ yields the morphism  of $F$-filtered $\bar{K}$-algebras $ R\Gamma_{\!\dR}(V)\to
R\Gamma_{\!\dR}^\natural (V,\bar{V})\hotimes \Bbb Q_p$. Let
  $\gamma : R\Gamma_{\!\dR} (V) \otimes_{\bar{K}} (\AdR \hotimes\Bbb Q_p )\to R\Gamma_{\! \dR}^\natural (V,\bar{V})\hotimes\Bbb Q_p$ be
 its $\AdR \hotimes\Bbb Q_p$-linear extension.

\proclaim{Lemma} $\gamma$  is a filtered quasi-isomorphism.
\endproclaim

\demo{Proof}
By \cite{B} (4.2.1), $\gr^m_F R\Gamma_{\!\dR}^\natural (V,\bar{V}) $ carry a finite filtration $I_\cdot$ and we have an identification  $\gr^I_a \gr^m_F R\Gamma_{\!\dR}^\natural (V,\bar{V})= R\Gamma (\bar{V},\Omega^a_{\langle V,\bar{V} \rangle})[-a] \otimes^L_{O_{\bar{K}}} \gr^{m-a}_F \AdR $. Here
$\Omega^a_{\langle V,\bar{V} \rangle} :=  \Omega^a_{(V,\bar{V})/(\bar{S},\bar{\CL})}  $ are relative differential forms with log singularities. The $I_a$'s are $\gr^\cdot_F \AdR$-submodules of $\gr^\cdot_F R\Gamma_{\!\dR}^\natural (V,\bar{V}) $ and the  identification is $\gr^\cdot_F \AdR$-linear.

Since $ R\Gamma (\bar{V},\Omega^a_{\langle V,\bar{V} \rangle})$ is a perfect $O_{\bar{K}}$-complex, applying $\cdot\hotimes\Bbb Q_p$ yields
a finite filtration $I_\cdot$ on $\gr^m_F R\Gamma_{\!\dR}^\natural (V,\bar{V}) \hotimes\Bbb Q_p$ together with a $\gr_F^\cdot\BdR^+$-linear identification
$\gr^I_a \gr^m_F R\Gamma_{\!\dR}^\natural (V,\bar{V})\hotimes\Bbb Q_p = \mathop\oplus\limits_{m\ge a\ge 0} \!\!
R\Gamma (\bar{V},\Omega^a_{\langle V,\bar{V} \rangle})[-a] \otimes_{O_{\bar{K}}} \!\gr^{m-a}_F \BdR^+$.

One also has an evident $\gr_F^\cdot \BdR^+$-linear identification $\gr^m_F (R\Gamma_{\!\dR} (V) \otimes_{\bar{K}} (\AdR \hotimes\Bbb Q_p ))$ $= \oplus_a \,\gr^a_F (R\Gamma_{\!\dR} (V)) \otimes_{\bar{K}}\gr^{m-a}_F \BdR^+ =\!\!\!\!  \mathop\oplus\limits_{m\ge a\ge 0} \!\!
R\Gamma (\bar{V},\Omega^a_{\langle V,\bar{V} \rangle})[-a] \otimes_{O_{\bar{K}}} \!\gr_F^{m-a}\BdR^+ $.

Since $\gr^\cdot_F \gamma$ sends $\gr_F^a R\Gamma_{\!\dR}(V)$ to $I_a \gr^a_F R\Gamma_{\!\dR}^\natural (V,\bar{V})\hotimes\Bbb Q_p  =\gr^a_F R\Gamma_{\!\dR}^\natural (V,\bar{V}) \hotimes\Bbb Q_p$, it sends
$\gr_F^a R\Gamma_{\!\dR}(V)\otimes\gr^{\cdot }_F \BdR^+ $ to $I_a \gr_F^{\cdot +a} R\Gamma^\natural_{\!\dR}(V,\bar{V})\hotimes\Bbb Q_p$. The lemma follows if we check that $ 
\gr_F^a R\Gamma_{\!\dR}(V)\otimes\gr^{\cdot }_F \BdR^+ \to \gr^I_a \gr_F^{\cdot +a} R\Gamma_{\!\dR}^\natural (V,\bar{V})\hotimes\Bbb Q_p$ is a quasi-isomorphism. By the above, it suffices to do this for $\cdot =0$. Here it is evident: the projection  
$ R\Gamma^\natural_{\!\dR}(V,\bar{V})\to R\Gamma (\bar{V},\Omega^\cdot_{\langle V,\bar{V}\rangle})$ provides, after applying $\gr_F^a \cdot \hotimes\Bbb Q_p$, the inverse map.
\qed \enddemo

\subhead{\rm 3.6}\endsubhead Let us consider (b) in 3.4.
The log schemes $(V,\bar{V})_n /W_n$ satisfy the condition of Remark (i) in 1.9, so we get canonical morphisms $$\kappa_n : R\Gamma_{\! \crys}(V,\bar{V})_{(n)} \to \lim_F R\Gamma_{\!\dR}^\natural (V,\bar{V})_n  .\tag 3.6.1$$ Namely, $\kappa_n$ is projective limit of the maps
$\kappa_{m,n}: R\Gamma_{\! \crys}(V,\bar{V})_{(n)} \to R\Gamma^\natural_{\!\dR} (V,\bar{V})_n /F^m$, where $\kappa_{m,n}$ equals the composition    $R\Gamma_{\! \crys}(V,\bar{V})_{(n)} \iso
R\Gamma (((V,\bar{V})_n /W_n )_{\crys},\CO )\to 
R\Gamma (((V,\bar{V})_n /W_n ),\CO/\CJ^{[m]})\iso  R\Gamma^\natural_{\!\dR}  (V,\bar{V})_n /F^m$, the first $\iso$ comes from  Remark (i) in 1.12, the second $\iso$ is the inverse to the Illusie-Olsson isomorphism (1.9.2). \remark{Example} For $(V,\bar{V})=\Spec\, (\bar{K},O_{\bar{K}})$,  our $\kappa_n$ is a canonical map  $ \text{A}_{crys\, n} \to \text{A}^\ds_{\dR \, n}$. By  Illusie-Olsson (see 1.9), it identifies 
$\text{A}_{\crys\, n} /J^{[m]}_n$ with $\text{A}_{\dR\, n} /F^m$, where 
$J_n$ is the kernel of the projection $\text{A}_{\crys\, n} \twoheadrightarrow \hat{O}_{\bar{K}}/p^n$. Thus $\AdR^\diamondsuit$ is the completion of $\Acrys$ with respect to the $J^{[m]}$-topology.
\endremark

\proclaim{Proposition} The $\text{A}_{\dR\, n}^\ds$-linear extension of $\kappa_n$ is an isomorphism: one has
 $$
 R\Gamma_{\! \crys}(V,\bar{V})_{(n)} \otimes_{\text{A}_{\crys\, n}}^L \text{A}^\ds_{\dR \, n} \iso \text{lim}_F R\Gamma_{\!\dR} (V,\bar{V})_n .\tag 3.6.2$$
\endproclaim

\demo{Proof}  By Proposition in 2.2,  the l.h.s.~in (3.6.2) equals $\text{holim}_m \, R\Gamma_{\! \crys}(V,\bar{V})_{(n)}\otimes_{\text{A}_{\crys\, n}}^L (\text{A}_{\crys\, n}/J^{[m]}_n )\iso  \text{holim}_m R\Gamma (( (V,\bar{V})_n / (\bar{S}_n ,  \text{A}_{\crys\, n}/J^{[m]}_n ))_{\crys} ,\CO_{\bar{V}_n / ( \text{A}_{\crys\, n}/J^{[m]}_n )})$ by base change and (1.18.1). Since $(V,\bar{V})_n $ is log smooth over $(\bar{S},\bar{\CL})_n$, the Hodge-pd filtration on $\Omega^{\cdot}_{\bar{V}_n /  (\text{A}_{\crys\, n}/J^{[m]}_n )}$ is finite (in fact, $F^N \Omega^{\cdot}_{\bar{V}_n / ( \text{A}_{\crys\, n}/J^{[m]}_n )}$ vanishes for $N\ge \dim V + m$), so the above completion equals the completion for the Hodge-pd filtration, which is the r.h.s.~of (3.6.2) by the Illusie-Olsson theorem in 1.9. \qed
\enddemo

\subhead{\rm 3.7}\endsubhead Consider finally (c) in 3.4. Using (3.4.2) and (3.4.3), the map (3.4.4) can be rewritten as $\phi : R\Gamma_{\!\crys}(V,\bar{V})\otimes^L_{\Acrys} \BdR^+ \to R\Gamma_{\!\dR}(V)\otimes_{\bar{K}}\BdR^+$. Both terms are perfect $\BdR^+$-complexes (see Proposition in 2.2) and $\BdR^+$ is a dvr, so to prove that $\phi$ is a quasi-isomorphism it suffices to check that its pullback $\phi_{\Bbb C_p} : R\Gamma_{\!\crys}(V,\bar{V})\otimes^L_{\Acrys}\! \Bbb C_p $ $\to R\Gamma_{\!\dR}(V)\otimes_{\bar{K}}\Bbb C_p$
 to $ \BdR^+ /\fm_{\dR} =\Bbb C_p$ is a quasi-isomorphism. We use 2.2 and the notation in loc.~cit. By (1.18.1) and base change, one has $R\Gamma_{\!\crys}(V,\bar{V})\otimes^L_{\Acrys}\! \hat{O}_{\bar{K}}\iso  R\Gamma (Z, \Omega_{(Z,\CM)/(S,\CL )})\otimes_{O_{K}}^L \!\hat{O}_{\bar{K}}$, hence  $R\Gamma_{\!\crys}(V,\bar{V})\otimes^L_{\Acrys}\! \Bbb C_p \iso R\Gamma_{\!\dR}(V)\otimes_{\bar{K}}\Bbb C_p$. Comparing it with the definition of $\kappa$ in (3.4.2), we see that this identification equals $\phi_{\Bbb C_p}$, and we are done. \qed

\head Index of notation.  \endhead

$\!\!\!\!\!\!\!(Z,\CM )$, $(Z,\CM )^{\text{int}}$, $\CM^{\text{gr}}$, $\Bbb A^{(1)}_{(S,\CL )}$ 1.1; 
$S^\sharp$, $G^\sharp$ 1.2; $\CT_{S^\sharp}$, $\CC_{S^\sharp}$ 1.3; 
$\CO_{Z/S}$, $\CJ_{Z/S}$, $\CM_{Z/S}$, $u^{\log}_{Z/S}$, $f_{\crys}$,
$((Z,\CM )/S^\sharp )_{\crys}$, $(Z/S)^{\log}_{\crys}$ 1.5; $(U, T_* )$, $ \CC^\cdot \CF$ 1.6; $\Omega^\cdot_{(U,T)/S}$, $\Omega^\cdot_{Z/S}$, $F^m$, $\nabla_{(Z,P)}$ 1.7; $L\Omega^\cdot_{(Z,\CM)/(S,\CL )}$ 1.9; $Fr_T$, $Fr'_T$, $C$ 1.10; $D^{\text{pcr}}$ 1.11; 
$i_n$, $\CF_{(n)}$, $(Z,\CM )_{\crys}$, $(Z,\CM )_{\crys (n)}$,
$R\Gamma_{\!\crys}(Z,\CM )_{(n)} $, $R\Gamma_{\!\crys}(Z,\CM )$, $H^i_{\crys}(Z,\CM )$,  $(Z,\CM )_n$
1.12; $R_{\varphi}$-mod, $\Hom^\natural_{R_{\varphi}}\!$, $D_{\varphi}(R)$, $R^{\text{prf}}_{\varphi}$-mod, $D^{\text{prf}}_{\varphi}(R)^{\text{nd}}$ 1.13; $D_{\varphi}((Z/S)^{\log}_{\crys},\CO_{Z/S})$,$D^{\text{pcr}}_{\varphi}(Z/S)$, $D^{\text{pcr}}_{\varphi}(Z/S)^{\nd}$
1.14; $K_0$,  $\Hom^\natural_{\varphi ,N}$,  $(\varphi ,N)$-mod, $(\varphi ,N)^{\text{eff}}$-mod, $D_{\varphi ,N}(K_0 )$, $D_{\varphi ,N}(K_0 )^{\text{eff}}$,  $\epsilon_{\pi}$, $\epsilon_l$, $\epsilon$, $\lambda_l$, $D^{\pcr}_{\varphi}(Y )^{\nd}$, $\tau_{A_{\Bbb Q}}$, $A_{\Bbb Q}^\tau$  1.15; $R\Gamma_{\!\text{H}\!\text{K}}(Z_1,\CM_1 )$, $\iota_{\dR}$,
$R\Gamma_{\!\text{H}\!\text{K}}(Z_1,\CM_1 )^\tau_{A_{\Bbb Q}}$ 1.16; A$_{\crys}$, $J_{\crys}$, $E_{\crys}$, $\CL_{\crys}$, $\CL_{\varphi}$, B$_{\crys}^+$, B$_{\st}^+$  1.17; $R\Gamma_{\!\text{H}\!\text{K}}(Z\bar{_1},\CM\bar{_1} )$, $\iota_{\crys}$, $K_0^{\text{nr}}$ 1.18;  $\cdot \hotimes\Bbb Z_p$, $\cdot \hotimes\Bbb Q_p$ 2.1;
$\CV ar^{\text{ss}}_{\bar{K}}$, $R\Gamma_{\!\crys}(V,\bar{V})$,   $R\Gamma_{\!\crys}(V,\bar{V})_{(n)}$ 2.2;  $\CA_{\crys}$ 2.3; $R\Gamma_{\!\crys}(X)$, $H^n_{\crys} (X )$ 2.4; $\CA_{\text{H}\!\text{K}}$, $R\Gamma_{\!\text{H}\!\text{K}}(X) $, $H^n_{\text{H}\!\text{K}}(X) $, $R\Gamma_{\! \dR}(X)$,
 $R\Gamma_{\!\text{H}\!\text{K}}(X)^\tau_{\Bcrys^+} $, $R\Gamma_{\!\text{H}\!\text{K}}(X)^\tau_{\bar{K}} $ 2.5; $\rho_{\crys}$, $\rho_{\text{H}\!\text{K}}$ 3.1; $R\Gamma_{\!\dR}^\natural (V,\bar{V})$, lim$_F$, A$_{\dR}$, A$_{\dR}^{\diamondsuit}$ 3.4.

\enddocument